\documentclass[11pt]{article}
\usepackage{amsfonts,amsmath,amssymb,amsthm,amscd,epsfig,epic,eepic}
\usepackage{graphics}

 \evensidemargin=\oddsidemargin

\makeatletter
\gdef\thmhead@plain#1#2#3{%
  \thmname{#1}\thmnumber{\@ifnotempty{#1}{ }#2}%
  \thmnote{ {\mdseries#3}}}
\let\thmhead\thmhead@plain
\makeatother
\theoremstyle{plain}
\newtheorem{theorem}{Theorem}[section]
\newtheorem{lemma}[theorem]{Lemma}
\newtheorem{proposition}[theorem]{Proposition}
\newtheorem{corollary}[theorem]{Corollary}
\newtheorem{crr}{Corollary}
\newtheorem*{cri}{Corollary}
\newtheorem*{lemi}{Lemma}
\newtheorem{pro}{Proposition}
\newtheorem{lem}{Lemma}
\newtheorem*{proi}{Proposition}

\newtheorem{assumption}{Assumption}
\newtheorem{thm}{Theorem}

\newtheorem*{thmi}{Theorem}

\newtheorem*{crrend}{Corollary }

\theoremstyle{definition}
\newtheorem{definition}[theorem]{Definition}
\newtheorem{remark}[theorem]{Remark}
\newtheorem{notation}[theorem]{Notation}

\theoremstyle{remark}

\def\alinea#1{\hfill\break%
  \hbox to \parindent{\hss{\upshape{\bf #1)}}\enspace}\ignorespaces}
\def\bul{\hfill\break\hbox to\parindent{\hss$\bullet$\enspace}\ignorespaces}
\def\up{\textup}
\def\from{\colon}
\def\up{\textup}
\def\from{\colon}
\def\mod{\operatorname{mod}}

\def\id{\hbox{Id}}
\def\h{h}
\def\p{p}
\def\H{H}
\def\C{\mathbf{C}}
\def\D{\mathbf{D}}

\def\M{\mathbf{M}}
\def\N{\mathbf{N}}
\def\Q{\mathbf{Q}}
\def\R{\mathbf{R}}
\def\S{\mathbf{S}}

\def\Z{\mathbf{Z}}
\def\a{\mathbf{a}}
\def\AA{\mathcal{A}}
\def\CC{\mathcal{C}}
\def\EE{\mathcal{E}}
\def\DD{\mathcal{D}}
\def\GG{\mathcal{G}}
\def\II{\mathcal{I}}
\def\HH{\mathcal{H}}
\def\Hinf{{\mathcal{H}}_\infty}
\def\LL{\mathcal{L}}
\def\MM{\mathcal{M}}
\def\PP{\mathcal{P}}
\def\QQ{\mathcal{Q}}
\def\RR{\mathcal{R}}
\def\SS{\mathcal{S}}
\def\UU{\mathcal{U}}
\def\VV{\mathcal{V}}
\def\XX{\mathcal{X}}
\def\WW{\mathcal{W}}
\def\Cap_#1{\bigcap\limits_{#1}}
\def\Cup_#1{\bigcup\limits_{#1}}
\def\ol{\overline}

\def\va{f_{\mathbf{a}}(-\mathbf{a})}
\def\wt{\widetilde}

\makeatletter
\def\cqfdsymb{\relax\protect\ifmmode\else\unskip\nobreak\fi
\quad\hfill$\bgroup
\vcenter{\hrule\hbox{\vrule\@height.6em\kern.6em\vrule}\hrule}\egroup$}
\def\cqfd{\cqfdsymb  \endtrivlist}
\gdef\rom#1{\leavevmode\skip@\lastskip\unskip\/%
        \ifdim\skip@=\z@\else\hskip\skip@\fi{\normalshape#1}}
\makeatother
\begin{document}
\title{{\bf Hyperbolic components  of polynomials with a fixed critical point of maximal order}}

\author{{\sc P. Roesch}\thanks{Research partially supported by the Bernouilli Center in Lausanne}}
\date{\today}
\maketitle
\begin{abstract} For the study of the $2$-dimensional space of
cubic polynomials, J.~Milnor considers the complex $1$-dimensional
slice $\mathcal{S}_n$ of the cubic polynomials which have a
super-attracting  orbit of period $n$. He gives in~\cite{M4} a
detailed conjectural picture of  $\mathcal{S}_n$. In this article,
we prove these conjectures for  $\mathcal{S}_1$ and generalize
these results in higher degrees. In particular, this gives a
description  of the closures of the hyperbolic components and of
the Mandelbrot copies sitting in the connectedness locus. We prove
that the closure of hyperbolic components is a Jordan curve, the
points of which are characterized according to their dynamical
behaviour. The global picture of the connectedness locus is a
closed disk together with ``limbs'' sprouting off at the cusps of
Mandelbrot copies  and whose diameter tends to $0$ (which
corresponds to  the Yoccoz inequality in the quadratic case).

\end{abstract}

\section*{Introduction}
In~\cite{BH} Branner and Hubbard  studied the parameter space of
cubic polynomials (in terms of the dynamics). This space has been
intensively studied since then. In this paper we focus on the
$1$-dimensional complex slice $\mathcal{S}_1$ of the cubic maps
which have a fixed critical point. Our goal is to prove the
conjectural picture given by Milnor in~\cite{M4} of the
connectedness locus in $\mathcal{S}_1$, and to generalize these
results to degrees $d \ge 3$.  The first question concerns the
topology of the boundary of the main hyperbolic component. This
leads naturally to the question of characterizing in dynamical
terms the parameters on this boundary. The other problems are to
describe first the
 intersections of the closures of hyperbolic components between
 each others and second  their intersection  with the Mandelbrot copies.
Milnor considered also the limbs which are attached to the main
hyperbolic component and raised the question of the existence of
an analogue of Yoccoz' inequality, namely, the inequality  which
bounds the size of the limbs in the Mandelbrot case.

Let us  consider the families of polynomials of degree $d\ge 3$
having a critical fixed point of maximal multiplicity. When we fix
the degree, $d$, this set of polynomials is described---modulo
affine conjugacy---by the following  family $\{f_{\mathbf{a}},\
\mathbf{a} \in \C\}$, where $0$ is the critical fixed point of
maximal multiplicity\,:
$$ f_{\mathbf{a}}(z)=z^{d-1}\left(z+\frac{d\mathbf{a}}{d-1}  \right).$$

\begin{figure}[!h]
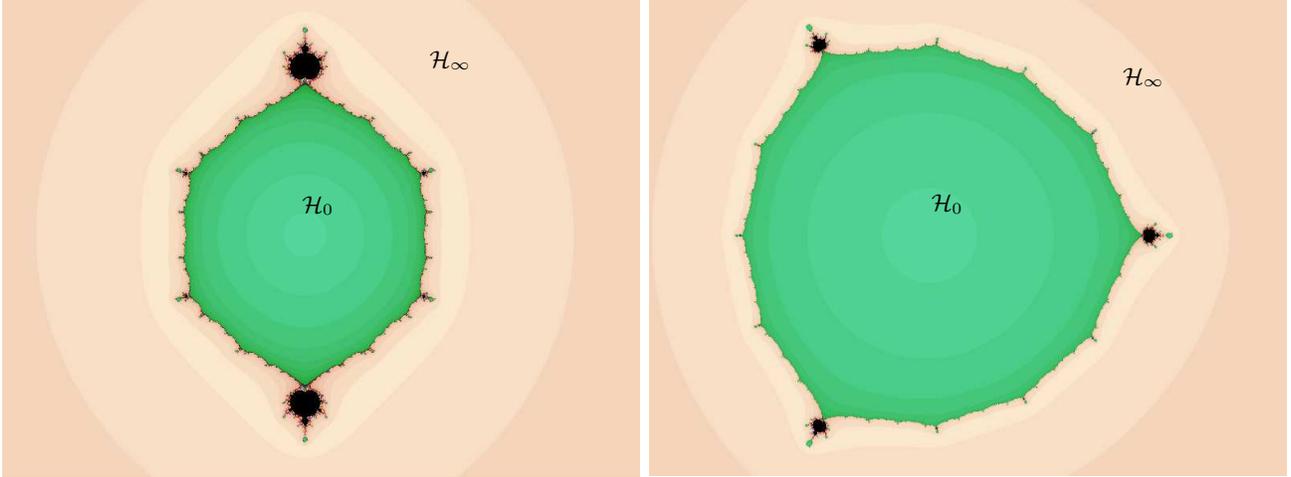
\vskip 0cm
\centerline{
\input{connexite3b.pstex_t}
\input{connexite4.pstex_t}
} \caption{Connectedness locus for $d=3$ and $4$ in dark color.
}\label{f:connexite34}
 \end{figure}

 The set of parameters is partitioned into two loci\,: $\C=\CC\sqcup\HH_\infty$
(after~\cite{BH}). The set $\CC$ denotes the {\it connectedness
locus} {\it i.e.}, the set of  parameters $\mathbf{a}$ such that
the Julia set $J(f_{\mathbf{a}})$ is connected\,;
 $\HH_{\infty}$ consists of the  parameters $\mathbf{a}$ such that
the ``free'' critical point~$-\mathbf{a}$ is attracted by $\infty$
(see~\cite{DH1}). We can continue the partition further,
considering the {\it hyperbolic} parameters
 {\it i.e.}, the parameters such that the
orbit of every critical point converges to an attracting cycle
(see~\cite{M1,M4}). This hyperbolic set is a disjoint union of
open disks called {\it hyperbolic components}. The locus
$\HH_{\infty}$ is the unique unbounded hyperbolic component (see
Lemma~\ref{l:connected} or  \cite{M4} and \cite{BH}). Among the
hyperbolic  components contained in $\CC$,  we focus on the ones
associated to the attracting point $0$. The union of those is
$\mathcal H=\{\mathbf{a}\in \C \mid -\mathbf{a}\in \wt
B_{\mathbf{a}}\}$, where $\wt B_{\mathbf{a}}$ is the basin of
attraction of the fixed point $0$.

\subsection*{Local connectivity of the boundary of hyperbolic components}

\begin{thm}\label{t:Hyp} The boundary of every hyperbolic
component of $\CC$ is a Jordan curve.
\end{thm}
 This Theorem is a consequence of Theorem~\ref{th:H} and
 the renormalization\footnote{Definitions of
``renormalization'' and of ``copies of $\M$'' are given in
section~\ref{s:graploccon} } property of
Proposition~\ref{p:periodicrenormalizable}.
\begin{thm}\label{th:H}
The boundary of every connected component of $\mathcal H$ is a Jordan curve.
\end{thm}
\begin{pro}\label{p:periodicrenormalizable}If the map
$f_{\mathbf{a}}$ has a non-repelling
periodic point $p\neq 0$  then $f_{\mathbf{a}}$ is renormalizable
near $p$ and the parameter $\mathbf{a}$ belongs to a copy of the
Mandelbrot set $\M$.
\end{pro}

Theorem~\ref{th:H} is the analogue in the parameter plane of the
following dynamical result\footnote{We will recall briefly the
proof of it from Yoccoz' Theorem in section~\ref{s:annexe2}}\,:
\begin{thmi}[\cite{F,Ro1}]\label{th:lcdyn} The boundary of every connected
component of $\wt B_{\mathbf{a}}$ is a Jordan curve.
\end{thmi} Let us recall that  D.~Faught gave a
 proof of Theorem~\ref{th:H} in his thesis~\cite{F}.
 This  result  remains  unpublished.
For completeness, we give a proof of this result of local
connectivity\,; our proof is different from that of~\cite{F}, the
argument here is based on an idea of Shishikura that simplifies
the analysis.

Proposition~\ref{p:periodicrenormalizable} has the following two
interesting corollaries\,:
\begin{crr}\label{c:H2} Any  hyperbolic component of  $\CC$ is either
 a connected component  of  $\mathcal H$ or a hyperbolic
component of
 a copy of $\mathbf M$.
\end{crr}

\begin{crr}\label{c:brujno} Assume that $f_{\mathbf{a}}$ possesses a periodic point $p$
with multiplier $\lambda=e^{2i\pi \theta}$, such that $\theta\in
\R\setminus\Q$.
  Then $f_{\mathbf{a}}$ is linearizable near $p$ if and only if $\theta \in \mathcal B$.
  Moreover, if
$\theta\notin \mathcal B$ there exist periodic cycles in any
neighbourhood of $p$.
\end{crr}
Here $\mathcal B$ denotes the set of {\it Brjuno} numbers\,: an
irrational $\theta$ of convergents $p_n/q_n$ (rational
approximations obtained by the continued fraction development) is
a {\it Brjuno} number,   if $\sum_{n=1}^{\infty}(\log
q_{n+1})/q_n$ is finite.

\subsection*{Parameters on the boundary of components of $\partial\HH$}

Let  $\HH_0$ be the connected component of $\mathcal H$ containing
$0$.

\begin{thm}\label{th:comp}Let $\mathbf{a}\in \partial\mathcal{H}_0$. There exists
a unique parameter ray\footnote{Rays and parameter rays are
defined in section~2 and 3}  in $\HH_0$   landing at $\mathbf{a}$,
say $\RR^s_0(t)$. The following dichotomy holds\,{\rm:}
\begin{itemize}

\item either there is a unique external parameter ray converging
to $\mathbf{a}$. In this case $f_{\mathbf{a}}$ is not
renormalizable so $\mathbf{a}$ do not belong to a copy of  $\M$.
Moreover in the dynamical plane, the ray $R_{\mathbf{a}}^0(t)$
lands at the critical value $f_{\mathbf{a}}(-\mathbf{a})\in
\partial B_{\mathbf{a}}$ and there is a unique external ray
converging to $f_{\mathbf{a}}(-\mathbf{a})$\,;

\item or there are   exactly two external parameter rays
 converging to $\mathbf{a}$. In this case  $\mathbf{a}$  is
  the cusp of a copy of $\mathbf M$.
  Furthermore, in the dynamical plane, the  ray
  $R_{\mathbf{a}}^0(t)$ lands at a parabolic point
on $\partial B_{\mathbf{a}}$. The angle $t$ is necessarily
periodic by multiplication by $d-1$.
 \end{itemize}
\end{thm}
\noindent Proposition~\ref{p:convzero} gives a criterium on the
angle $t$ to decide which one of the two cases described above
arises.

\begin{thm}\label{th:U}
Let $\mathbf{a}\in \partial \mathcal U$ where $\mathcal U\neq
\HH_0$ is a connected component of $\mathcal H$. Then $\mathbf{a}$
is the landing point of a unique parameter ray of $\UU$, say
$\RR_\UU(t)$. In the dynamical plane, some iterate
$f_{\mathbf{a}}^k(-\mathbf{a})$ lies in $ \partial B_{\mathbf{a}}$
but $-\mathbf{a}\notin\partial B_{\mathbf{a}}$. Moreover, there
exists a holomorphic function $r$, defined in a neighbourhood of
$\mathbf{a}$, such that the dynamical ray
$R^{r(\mathbf{a})}_{\mathbf{a}}(t)$ converges to the critical
value $\va$. As a consequence, $f_{\mathbf{a}}$ has no parabolic
cycles.
\end{thm}

\begin{figure}[!h]\vskip 0cm
\centerline{
\input{MC5.pstex_t}
} \caption{A copy of $\M$ attached to $\HH_0$ and a component
$\UU$ of $\HH\setminus \HH_0$. } \label{f:copieM}
 \end{figure}

\begin{crr}{\rm(see also}~{\rm\cite{GM})}\label{c:indiff}
For parameters $\mathbf{a}$ on the boundary of a  component of
$\HH$, $f_{\mathbf{a}}$ cannot have an irrational indifferent
periodic point.
\end{crr}

\subsection*{ Intersections between  the closures of hyperbolic components}
\begin{lem}\label{l:inter} Any two distinct components of $\HH$ have disjoint closures.
\end{lem}
Recall that the cusp of $\M$ is the point $c=1/4$, and that the
tips of $\M$ are the parameters $c\in \M$ such that $c$ falls
after some iterations on the repelling fixed point, $\beta_c$ (the
one  that does not disconnect the Julia set). The cusp and the
tips of a copy of $\M$ are the corresponding images by the
homeomorphism defining the copy. \vskip 1em Let $\M_0$ be a copy
of $\M$.
\begin{pro}\label{th:intercopy}
If  $\M_0$ intersects $\partial \UU$ where $\UU$ is a component of
$\HH$, the following dichotomy holds\,:
\begin{itemize} \item If $\UU= \HH_0$, $\M_0\cap \partial \UU$ is reduced
to a single point, which is the cusp of $\M_0$\,; \item If
$\UU\neq \HH_0$, $\M_0\cap \partial \UU$ is reduced to a single
point, which is a tip of $\M_0$. Furthermore,  $\M_0\cap\partial
\HH_0$ is not empty, it reduces to the cusp of $\M_0$.
\end{itemize}
\end{pro}
 Conversely,
\begin{pro}\label{p:tips} If  $\M_0$ intersects $\partial \HH_0$,
then at any of its tips there is a connected component of $\HH
\setminus \HH_0$ attached.
\end{pro}

These results describe all the intersections between the
boundaries of  components of $\HH$ and also with copies of
$\mathbf{M}$, so in particular between all hyperbolic components of $\C$.

\begin{thm}\label{th:conections} The only intersections
between closure of hyperbolic components, and also copies of $\M$
are the following\,:
\begin{itemize}
\item The central component $\HH_0$ has Mandelbrot copies $\M_t$
attached to it at angles $t$ which are  periodic by multiplication
by $d-1$ (a full characterization of these values is given in
Proposition~\ref{p:convzero})\,;

\item At every tip  of such a satellite $\M_t$, there is a
component $\UU$ of $\HH\setminus \HH_0$ attached.
\end{itemize}
Nevertheless, there are infinitely many  copies of $\M$ in $\CC$
and components of $\HH$ not  in the category described above.
\end{thm}

\subsection*{Some global properties of $\CC$}

\begin{thm}\label{th:Clc} $\partial \CC$ is locally connected at every point which
 is not in a  copy of $\M$ and at any point of $\partial \UU$ for
every connected components $\UU$ of $\HH$.
\end{thm}

Concerning the {\it limbs}\footnote{defined in
section~\ref{s:limbs}} of the main component $\HH_0$, we obtain a
qualitative version of Yoccoz' inequality for this family\,:
\begin{thm}\label{th:limbs} For any $\epsilon>0$, only a finite number
of limbs have  diameter greater than $\epsilon$.
\end{thm}

\subsection*{Description of  the content of the article}

In the first section we give some  properties of the polynomials
$f_{\mathbf{a}}$ in dynamical and parameter plane.

The second section is devoted to the parametrization of the
components of $\HH$ and also of $\HH_\infty$. The parametrization
is given by the B\"ottcher coordinate of the critical value and
provides {\it parameter rays} and {\it equipotentials}.

In section $3$, we construct graphs that define puzzles to prove
the local connectivity in the parameter space. They
correspond---via the parametrization---to those used in the
dynamical plane in~\cite{Ro1} for the proof of the local
connectivity of $\partial B_{\mathbf{a}}$ (we will recall the
construction and the results of \cite{Ro1}). Then, the holomorphic
motion of the dynamical graphs allows us to compare the puzzles in
the parameter plane and in the dynamical plane, as pointed out by
Shishikura in the case of quadratic polynomial (see~\cite{Ro2}).

Section $4$ is devoted to the proof of
Theorem~\ref{t:Hyp},~\ref{th:H},~\ref{th:comp} and~\ref{th:U}.
Namely, we prove that when the intersection of the puzzle pieces
(in the parameter plane) is not reduced to a single point, then
this intersection is a copy of the Mandelbrot set.
 
In section~5 we give  the announced  description of $\CC$\, {\it
i.e.} the proof of Theorems~\ref{th:conections},~\ref{th:Clc}
and~\ref{th:limbs}.

Finally, we add in the Appendix (section~\ref{s:annexe2}) the
proof of Theorem~\cite{F,Ro1}.

\vskip 1em
{\it Acknowledgment}\,: I would like to thank Tan Lei, Carsten
Petersen and Curt Mc Mullen for encouraging me to write this down.
I would also like to thank John Milnor for suggesting me
section~\ref{s:annexe1}.

\section{Overview of dynamical and parameter plane }


Through the article we will take angles in $\R / \Z$ but in
general we will have in mind their representant in $[0,1[$. We
will write $dt$ for the image of the angle $t\in \R/\Z$ by the
multiplication by $d$ and $t/d$ for the element whose representant
is in $[0,1/d[$.

\subsection{The dynamics of $f_\mathbf{a}$}\label{s:dyn}
We consider the polynomials $f_{\mathbf{a}}$  for a fixed degree
$d\ge 3$. Note that for $d=2$ only the polynomial $P(z)=z^2$ does
satisfy the condition to have a super-attracting fixed point
(modulo affine conjugacy).

Recall that the {\it filled Julia set} $K(f_{\mathbf{a}})$
consists in the non escaping points and that the {\it Julia set}
$J_{\mathbf{a}}=J(f_{\mathbf{a}})$
 is its boundary\,: $$K_{\mathbf{a}}=K(f_{\mathbf{a}})=\{ z \mid f_\mathbf{a}^{n } (z)\hskip
1em\not \hskip -1em \xrightarrow[n\to\infty] {} \infty\},\quad
 J_\mathbf{a}=\partial K_\mathbf{a}.$$  Note that for every
 $\mathbf{a} \in \C$, $K_\mathbf{a}\neq J_\mathbf{a}$
 since $K_\mathbf{a}$ contains the {\it basin of attraction} of $0$
{\it i.e.}
 $$\wt B_\mathbf{a}=\{z \in \C \mid f_\mathbf{a}^{n }(z) \xrightarrow[n\to\infty] {} 0\}.$$
We denoted by $B_\mathbf{a}$ the {\it immediate basin} of $0$,
that is the connected component  of $\wt B_\mathbf{a}$
containing~$0$. It follows from the  maximum principle that
$B_\mathbf{a}$ is a topological disc. If
 $-\mathbf{a}\notin B_\mathbf{a}$ the map ${f_\mathbf{a}}{\vert _{B_\mathbf{a}}}$
 is conjugated to $z^{d-1}$ on
 $\D$, else $B_\mathbf{a}=\wt B_\mathbf{a}$ (see~\cite{Ro1} and
 the following B\"ottcher's Theorem).

 \begin{thmi}{\rm [B\"ottcher]}\label{th:bottcher}
For $\p= 0$ or $\infty$, there exist  neighborhoods
$V_{\mathbf{a}}^\p, W_{\mathbf{a}}^\p$ of $\p$ such that
$f_{\mathbf{a}}(V_{\mathbf{a}}^\p)\subset V_{\mathbf{a}}^\p$ ,and
conformal isomorphisms $\phi_{\mathbf{a}}^\p \from
V_{\mathbf{a}}^\p \to W_{\mathbf{a}}^\p$ satisfying
$$ \phi_{\mathbf{a}}^\infty\circ
f_{\mathbf{a}}=(\phi_{\mathbf{a}}^\infty)^{d}\quad\hbox{on} \quad
V_{\mathbf{a}}^\infty \quad\hbox{and
}\quad\phi_{\mathbf{a}}^0\circ
f_{\mathbf{a}}=(\phi_{\mathbf{a}}^0)^{d-1}\quad \hbox{on} \quad
V_{\mathbf{a}}^0 \quad (*)$$ with $\phi_{\mathbf{a}}^\infty$
tangent to identity near $\infty$ and $\phi_{\mathbf{a}}^0$
tangent to $z\mapsto \lambda(\mathbf{a})z$ near $0$ where
$\lambda(\mathbf{a})$ is
a   $(d-2)$-th  root of $\frac{d\mathbf{a}}{d-1}$.
 \end{thmi}

  \begin{remark}\label{r:bottcher}  The map
$\phi_{\mathbf{a}}^\infty$ is always unique. Moreover, if we fix
the choice of the $(d-2)$-th root $\lambda(\mathbf{a})$, the map
$\phi_{\mathbf{a}}^0$ is also unique.
  \end{remark}

  \begin{assumption}\label{a:root} Through all the paper  we will
  only consider parameters in
$ \C\setminus \R^-$ {\rm(}because of Remark~\ref{r:symm}{\rm)}.
Thus, for the choice of $\lambda(\mathbf{a})$ we take  the
$(d-2)$-th principal root of $\frac{d\mathbf{a}}{d-1}$, {\it
i.e.}, the one such that $\lambda(\R^+)\subset \R^+$.
  \end{assumption}

 \noindent The Green function $G^{\infty}_\mathbf{a} $ (resp. $G^{0}_\mathbf{a} $)
 associated to $\infty$ (resp.  to $0$)
 is the harmonic map  equal to $\log|\phi^\infty_\mathbf{a}(z)|$ on
$V_\mathbf{a}^\infty$ (resp.  to $-\log|\phi^0_\mathbf{a}(z)|$ on
$V_\mathbf{a}^0$), extended on $\ol \C \setminus K_\mathbf{a}$
through the relation $d
G^{\infty}_\mathbf{a}(z)=G^{\infty}_\mathbf{a}(f_\mathbf{a}(z))$
 (resp. on $\wt B_\mathbf{a}$ through
 $(d-1)G^{0}_\mathbf{a}(z)=G^{0}_\mathbf{a}(f_\mathbf{a}(z))$) and
 vanishing on the complement.

\begin{definition} The \emph {equipotential} of level~$v>0$, $E^{\p}_\mathbf{a}(v)$,
 around $p=0$ or $\infty$ is the curve
 $E^{\p}_\mathbf{a}(v)=\{ z \in \C \mid G^{\p}_\mathbf{a}(z) = v \}$.
A \emph{ray}, $R^{\p}_\mathbf{a}(t)$, of angle $t \in \R/\Z$,
stemming from $p=0$ or $\infty$, is  a gradient line of
$G^{\p}_\mathbf{a}$ that coincides near $p$ with
$(\phi^\p_\mathbf{a})^{-1}(\R^+ e^{2i\pi t})$.
  \end{definition}
Note that if there is no critical point of $G_\mathbf{a}^p$ on a
ray, it is a smooth simple  curve\,; whereas at the critical point
the gradient line divides itself so that several points have the
same potential on this ray and also different points of
$J_\mathbf{a}$ might be on the closure of a unique ray of given
angle. This happens if $\mathbf{a} \in \HH_\infty\cup \HH_0$. In
this case the angle $t$ and the potential  do not define uniquely
points on the ray, we will say that the ray is {\it not well
defined}. In fact the ray is well defined in $\{ z \mid
G^{\p}_\mathbf{a}(z)> G^{\p}_\mathbf{a}(-\mathbf{a}) \}$ since
$\phi_\mathbf{a}^\p$ extends
 to this set (via ($*$)).
 Note  that it  is possible
to  extend $\phi^0_\mathbf{a}$ continuously at the critical point
$-\mathbf{a}$, but not $\phi^\infty_\mathbf{a}$ since  there are
two external rays crashing on $-\mathbf{a}$.

Finally, if a ray $R_\mathbf{a}^p(t)$ is well defined, it
accumulates on the Julia set. We say that it {\it lands }if its
accumulation set is reduced to one point, the landing point is in
$J_\mathbf{a}$.

We have the following behaviour for rational rays
(see~\cite{DH1,M1,Pe})\,:
\begin{lemma}\label{l:douady1}  Let $\mathbf{a}_0\in \C$, $p=0$ or
$\infty$,
 and $t\in \Q/\Z$. If  the ray $R_{\mathbf{a}_0}^\p(t)$ is well defined, it lands
at an eventually periodic point which is repelling or parabolic
\footnote{A point $x$ of period $p$ is \emph{repelling},
\emph{attracting}, \emph {parabolic} or \emph{indifferent
irrational}, respectively, if $|(f^p)'(x)|>1$, $|(f^p)'(x)|<1$,
$(f^p)'(x) = e^{2i\pi \theta}$ and  $\theta \in \Q/\Z$, or $\theta
\in (\R\setminus \Q)/\Z$.}.
\end{lemma}

\begin{lemma}\label{l:douady}
Under the assumptions of lemma~\ref{l:douady1}, if the landing
point is repelling and not eventually a critical point, there
exists a neighbourhood $A$ of $\mathbf{a}_0$ such that for all
$\mathbf{a}\in A$ the ray $R_{\mathbf{a}}^\p(t)$ lands at a
repelling point. Moreover, the map $(\mathbf{a},s)\mapsto
\psi^p_{\mathbf{a},t}(s)$ is continuous on $A\times [0,\infty]$
and holomorphic in $\mathbf{a}$, where $\psi^p_{\mathbf{a},t}(s)$
is the point on $R_{\mathbf{a}}^\p(t)$ of potential $s$.
\end{lemma}

\begin{proposition}{\rm[Yoccoz]}\label{p:yocc}
For every eventually periodic point of $f_\mathbf{a}$ that is
repelling or parabolic, there exists a rational angle $t$ such
that $R^\infty_\mathbf{a}(t)$ lands at this point if
$J(f_\mathbf{a})$ is connected.
\end{proposition}

We will not need the analogue result in the non connected case
since we will start from the rays obtained when the Julia set is
connected and proceed to a holomorphic motion.

\begin{lemma}\label{l:aboutdistray0} If two rays $R^0_\mathbf{a}(t)$ and $R^0_\mathbf{a}(t')$ land,
their landing points are distinct when $t\neq t'$.
\end{lemma}
\proof This follows, by a contradiction argument, from the maximum
principle applied to the iterates of $f_\mathbf{a}$ on the domain
bounded by the closed curve
$\ol{R^0_\mathbf{a}(t)}\cup\ol{R_\mathbf{a}^0(t')}$.\cqfd

\subsection{Parameter plane }
\begin{remark}\label{r:symm}The rotation $\tau(z)=\tau z$,
where $\tau=e^{\frac{2i\pi}{d-1}}$, is  the only possible
conformal conjugacy between  polynomials
 $f_\mathbf{a}$ and $f_{\mathbf{a}'}$\,;  the relation is
 $f_{\tau \mathbf{a}}(\tau z)=\tau f_\mathbf{a}(z)$.
Besides this, $f_\mathbf{a}$ is conjugated to
$f_{\ol{\mathbf{a}}}$ by the complex conjugacy $\sigma(z)= \ol z$.
\end{remark}

Hence a ``fundamental domain'' for the study of the family
$f_\mathbf{a}$ is
 $$\mathcal{S}=\left\{\mathbf{a} \in \C \mid 0\le \arg(\mathbf{a})\le \frac{1}{2(d-1)}\right\}.$$

\begin{figure}[!h]
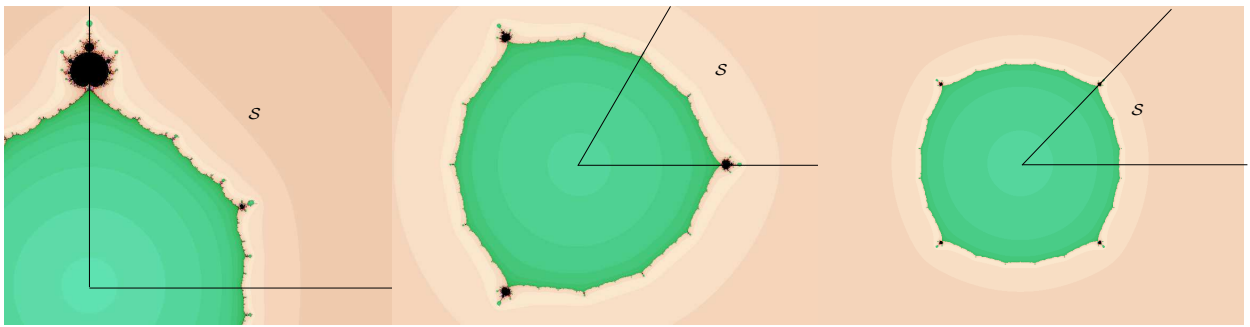
\vskip 0cm
\centerline{
\input{domaineC3.pstex_t}\hskip -2em\input{domaineC4.pstex_t}\input{domaineC5.pstex_t}}
\caption{Fundamental domain in degree $d=3,4$ and $5$.
}\label{f:connexite345}
  \end{figure}

 The  connectedness locus  $\mathcal C$,
 {\it i.e.}, the set of parameters $\mathbf{a}$ such that $K_\mathbf{a}$ (or equivalently $J_\mathbf{a}$)
 is connected, admits the following classical  characterization (see~\cite{DH1})\,:
$$\mathcal C=\{\mathbf{a} \in \C \mid f_\mathbf{a}^{n }(-\mathbf{a})\hskip 1em\not \hskip -1em
\xrightarrow[n\to\infty] {} \infty\}.$$

\begin{remark}The sets $\CC$, $\HH$ and $\HH_\infty$ admit $\sigma$
and $\tau$ as symmetries.
\end{remark}
\begin{lemma}\label{l:connected}
 The set  $\HH_\infty=\C\setminus \CC$ is a connected component
  of the set of  hyperbolic parameters.
  Similarly, the connected component $\HH_0$ is exactly
 the set $\{\mathbf{a}\mid -\mathbf{a}\in B_\mathbf{a}\}$.
  \end{lemma}
\proof Clearly $\HH_\infty$ contains a neighbourhood of $\infty$.
So if $\HH_\infty$ is not connected there is a  bounded connected
component, $\UU\subset\HH_\infty$. The boundary of $\UU$ is in
$\CC$ so there exists some $M\in \R$ such that
$|f_\mathbf{a}^n(-\mathbf{a})|\le M$ for all $n\ge 0$ and $
\mathbf{a} \in
\partial \UU$. For $\mathbf{a}_0$ a parameter in $\UU$, there exists some
$N$ such that $|f_{\mathbf{a}_0}^n(-\mathbf{a}_0)|\ge 2M$ for
$n\ge N$. This contradicts the maximum principle for the function
$\mathbf{a}\mapsto f_\mathbf{a}^N(-\mathbf{a})$.

The proof for $\HH_0$ goes with the same arguments exchanging $0$
and $\infty$. Assume (by  contradiction) that there is a connected
component   $\UU\subset\C$ of $\{\mathbf{a}\mid -\mathbf{a}\in
B_\mathbf{a}\}$ which is different from $\HH_0$.
There  exists
 some $\epsilon>0$ such that on $\partial \UU$,
 $|f_\mathbf{a}^n(-\mathbf{a})|>\epsilon$ and for $\mathbf{a}_0 \in \UU$ there exists some
 $N$ such that  $|f_{\mathbf{a}_0}^n(-\mathbf{a}_0)|<\epsilon/2$ for $n\ge N$.
This contradicts the maximum principle (as before) for the map
$\mathbf{a}\mapsto 1/f_\mathbf{a}^N(-\mathbf{a})$ which is well
defined on a neighbourhood of $\ol \UU$ since
$f_\mathbf{a}^N(-\mathbf{a})\neq 0$ for $\mathbf{a} \in\UU$
  because the sequence
 $f_\mathbf{a}^n(-\mathbf{a})$ tends geometrically to $0$ in $B_\mathbf{a}$.
 \cqfd

\begin{definition} The so-called {\it capture
 components} of depth $i\ge 1$ are the connected components of
 $\HH_i$, where
 $ \mathcal H_i=\{\mathbf{a}\in \C \mid
f_\mathbf{a}^i(f_\mathbf{a}(-\mathbf{a}))\in B_\mathbf{a} \hbox{
and } f^{i-1}(f_\mathbf{a}(-\mathbf{a}))\notin B_\mathbf{a}\}$.
\end{definition}
We have the following  decomposition of the hyperbolic components
of $\HH$\,:\begin{remark}
 $\HH=\Cup_{i\ge 0}\HH_i$.
\end{remark}
 \proof For $\mathbf{a}\in \HH$, the critical point $-\mathbf{a}$ is attracted by
 $0$, so there exists $k\ge 0$ such that $f_\mathbf{a}^k(-\mathbf{a}) \in B_\mathbf{a}$. If
$\mathbf{a}\notin \HH_i$ for any $i\ge 0$, necessarily
$f_\mathbf{a}(-\mathbf{a})\in B_\mathbf{a}$ with
$-\mathbf{a}\notin B_\mathbf{a}$. This is not possible since any
point near $f_\mathbf{a}(-\mathbf{a})$ in $B_\mathbf{a}$ would
have $d-1$ preimages in $B_\mathbf{a}$ (see B\"ottcher's Theorem)
plus two near $-\mathbf{a}$ (the critical point), which exceeds
the degree~$d$ of $f_\mathbf{a}$.\cqfd

 A rough picture of the dynamics of $f_\mathbf{a}$ for $\mathbf{a}\in\HH\cup\HH_\infty$ is the following.
 For parameters $\mathbf{a}$ in $\mathcal H_\infty=\C \setminus \mathcal C$,
the filled Julia set $K_\mathbf{a}$ is not connected but not
totally disconnected since it contains the closed disc $\ol
B_\mathbf{a}$.
 More precisely, $K_\mathbf{a}$ is the disjoint union of all the inverse
  images of $\ol B_\mathbf{a}$, the dynamics of
$f_\mathbf{a}$ on  $K_\mathbf{a}=\cup f_\mathbf{a}^{-i}(\ol
B_\mathbf{a})$ is the "shift"
 and on $\ol B_\mathbf{a}$ it is conjugated to $z \mapsto z^{d-1}$ on $\ol{\D}$.

 For parameters $\mathbf{a}$ in $\HH_0$ the critical point $-\mathbf{a}$ is in the immediate basin $
B_\mathbf{a}$. Indeed, for the center $\mathbf{a}=0$ this is clear
and the situation is stable.
 Thus for  $ \mathbf{a} \in \HH_0$, the dynamics is very simple\,: $K_\mathbf{a}=\ol B_\mathbf{a}$, $B_\mathbf{a}=\wt B_\mathbf{a}$,
    $J_\mathbf{a}$ is a quasi-circle and $f_\mathbf{a}|_{J_\mathbf{a}}$ is
quasi-conformally conjugated to $z\mapsto z^d$ on $\S^1$.

For parameters $\mathbf{a} \in \HH_i$ with $i\ge 1$,
$K_\mathbf{a}=\bigcup_{j\ge 0} f_\mathbf{a}^{-j} (\ol
B_\mathbf{a})$, the map $f_\mathbf{a}$ is conjugated to $z \mapsto
z^{d-1}$ on $\ol B_\mathbf{a}$ and corresponds to the "shift" on
the components of $(f_\mathbf{a}^{-j} (\ol B_\mathbf{a}))_{j\ge
0}$ not containing the critical point.

  \begin{lemma}\label{l:1connected}
  Any  connected component of $\HH$, as well as
   $\HH_{\infty} \cup \{ \infty\}$,  is simply connected.
  \end{lemma}
\proof Once more this is an application of the maximum principle.
The proof goes by contradiction. For seek of simplicity we give
the proof for a connected component $\UU$ of $\HH_n$ with $n\ge
0$. The argument follows for $\HH_{\infty} \cup \{ \infty\}$ by
exchanging  $0$ and $\infty$.

Assume by contradiction that there exists a bounded connected
component  $K$ of $\C\setminus \UU$.
 Then there exist points of $\partial \UU$ in $K$ and also a
 simple closed curve $\gamma\subset \UU$ surrounding $K$ since
 $\UU$ is arcwise connected. In $\UU$, the iterates of the  critical point $-\mathbf{a}$
 converge to $0$. Thus, for every $\epsilon>0$, there exists an $N\ge n$
such that for every $j\ge N$ and every parameter $\mathbf{a}\in
\gamma$, $|f_\mathbf{a}^j(-\mathbf{a})|<\epsilon$. Now let $x$ be
a point in $K\cap
\partial \UU$. Since $x\notin \HH$, there exists an $r>0$
(depending on $x$) such that the iterates $f^j_x(-x)\cap
B(0,r)=\emptyset$ for every $j\ge 0$. Taking $\epsilon=r/2$, this
contradicts the maximum principle for the holomorphic function
$g(z)=f^N_z(-z)$, on the bounded open set delimited by~$\gamma$.
 \cqfd

  \begin{notation}Let  $\Upsilon_{0} \from \HH_0 \to \D$,
  resp.   $\Upsilon_{\infty} \from \HH_\infty \to\C\setminus\overline{\D}$,
   be  the conformal representation
tangent  to the  identity at  $0$, resp. at $\infty$.
  \end{notation}

  \begin{remark}\label{r:symetry}Then, for $p=0$ or
$\infty$, $\Upsilon_p (\sigma \mathbf{a})=\sigma{ \Upsilon_p (
\mathbf{a})}$ and $\Upsilon_p (\tau \mathbf{a})=\tau\Upsilon_p (
\mathbf{a})$ with $\tau=e^{\frac{2i \pi}{d-1}}$ and
$\sigma(\mathbf{a})=\ol{\mathbf{a}}$. In other words  $\HH_p$
admits $\sigma$ and $\tau$ as symmetries.
  \end{remark}
\proof Since $\HH_p$ is invariant by the complex involution
$\sigma$ and the rotation $\tau$ (see Remark~\ref{r:symm}), the
maps $\sigma {\Upsilon_p(\sigma z )}$ and $\Upsilon_p(\tau
z)/\tau$ are
  conformal representations of $\HH_p$ onto $\D$, or $\C \setminus \ol \D$,
which are tangent to the identity at $0$, or at $\infty$. Hence
$\Upsilon_p(\sigma z)=\sigma {\Upsilon_p(z)}$ and $\Upsilon_p(\tau
z)=\tau \Upsilon_p(z)$. \cqfd

  \begin{corollary}\label{c:domaineimage} Let $\rho=e^{\frac{i
\pi}{d-1}}$, for any $k\in \N$ the line $\rho^k\ \R^+$ cut $\HH_0$
and $\HH_\infty$ under a connected set. As a consequence,
$\Upsilon_0(\rho^k\ \R^+)= \rho^k [0,1[$ and
$\Upsilon_{\infty}(\rho^k\ \R^+) =\rho^k]1,+\infty]$.
  \end{corollary}
This does not imply that $\R^+$ crosses only $\HH_0$, $\HH_\infty$
and $\partial \CC$ (see corollary~\ref{c:R+}). \proof Fix $p\in
\{0,\infty\}$. By remark~\ref{r:symetry}, for $\mathbf{a} \in
\R^+$, $\Upsilon_p (\mathbf{a}) \in \R^+$ and
$\Upsilon_p(\mathbf{a}\rho^k) \in \rho^k\ \R^+$ (since $\rho^k
\mathbf{a}=\rho^{2k}\sigma(\rho^k
\sigma(\mathbf{a}))=\rho^{2k}\sigma(\rho^k \mathbf{a})$ for
$\mathbf{a}\in\R$, where $\sigma$ denotes the complex conjugacy).
To prove that $\R^+\cap \HH_p$ is connected we apply the maximum
principle to a loop that we construct now.
 Assume by contradiction that there exist  $x_0<x<x_1$ with $x\notin\HH$,
 $x_0,x_1\in \HH_p$\,;
then there is  a simple arc  $\gamma_0\subset\HH_p$ with endpoints
$x_0$ and $x_1$  such that $\gamma_0\setminus \{x_0,x_1\}$ avoids
$\R$ (otherwise we exchange $x_0$ and $x_1$). The desired loop in
$\HH_p$ surrounding~$x$ is $\gamma=\gamma_0\cup\sigma(\gamma_0)$
(obtained by adding the conjugate). Hence $\R^+\cap \HH_p$ is
connected by the same argument as in Lemma~\ref{l:1connected}.
Using the symmetries $\sigma$ and $\tau$ (Remark~\ref{r:symetry})
we deduce that for $k\ge 1$, the set $\rho^k\R^+\cap \HH_p$ is
also connected. \cqfd

The following Lemma will be useful for describing the domains of
parametrization of the $\HH_i$ for $i\in\N\cup\{\infty\}$.  It
gives some symmetry properties of the rays. These properties  are
specific to the family under consideration.
\begin{lemma}\label{l:symray}
Fix $p\in \{0,\infty\}$.
  If the parameters $\mathbf{a}$ and $\tau \mathbf{a}$ are in $\C
 \setminus \R^{-}$, the  B\"ottcher maps are related by
  some constants $\kappa_p$ as follows\,:
   $ \sigma({\phi^p_{\sigma (\mathbf{a})}(\sigma
(z))})=\phi^p_\mathbf{a}(z)=\kappa_p(\mathbf{a}) \phi^p_{\tau
\mathbf{a}}(\tau z)$, with $\kappa_\infty(\mathbf{a})=\tau^{-1}$
and $\kappa_0(\mathbf{a})=\frac{\lambda(\mathbf{a})}{\tau
\lambda(\tau \mathbf{a})}$.
 Then the rays at parameters $\mathbf{a}$, $\tau \mathbf{a}$
 and $\sigma(\mathbf{a})$ satisfy   the following  relations,
 where $t_p(\mathbf{a})=\arg\left(\kappa_p(\mathbf{a})\right)$\,:
  $$  R^\p_{
\sigma(\mathbf{a})}(t)= \sigma \left(R^\p_{
\mathbf{a}}\left(-t\right)\right)\hskip 1em \hbox{and} \hskip 1em
 R^\p_{\tau \mathbf{a}}(t)= \tau R^\p_{\mathbf{a}}\left(t+t_p(\mathbf{a})\right ).$$
  \end{lemma}
 \proof
Since $\tau^{-1}f_{\tau \mathbf{a}}(\tau z)=f_\mathbf{a}(z)$, the
map $\tau^{-1} \phi^{\infty}_{\tau \mathbf{a}}(\tau z)$ conjugates
$f_\mathbf{a}$ to $z\mapsto z^d$ near $\infty$. Since it is
tangent to identity at $\infty$, $\tau^{-1} \phi^{\infty}_{\tau
\mathbf{a}}(\tau z)= \phi^{\infty}_\mathbf{a}(z)$. Applying  the
same argument to the maps $\phi^{0}_{\tau \mathbf{a}}(\tau z)$ and
$\sigma(\phi^p_{\sigma (\mathbf{a})}(\sigma (z)))$, we obtain that
$\phi^{\infty}_{\sigma (\mathbf{a})}(\sigma (z))=\sigma(
\phi^{\infty}_{\mathbf{a}}(z))$ and that $\kappa'_0(\mathbf{a})
\sigma(\phi^0_{\sigma (\mathbf{a})}(\sigma
(z)))=\phi^0_\mathbf{a}(z)=\kappa_0(\mathbf{a}) \phi^0_{\tau
\mathbf{a}}(\tau z)$ where $\kappa_0(\mathbf{a})$ and
$\kappa'_0(\mathbf{a})$ are appropriate constants. Taking the
derivatives at $0$, we obtain
$\kappa'_0(\mathbf{a})\sigma(\lambda(\sigma(
\mathbf{a})))=\lambda(\mathbf{a})=\kappa_0(\mathbf{a})\tau{\lambda(\tau
\mathbf{a})}$. Note that $\kappa'_0(\mathbf{a})=1$ since
$\lambda=\sigma\circ\lambda\circ \sigma$.
 \cqfd

\begin{figure}[!h]
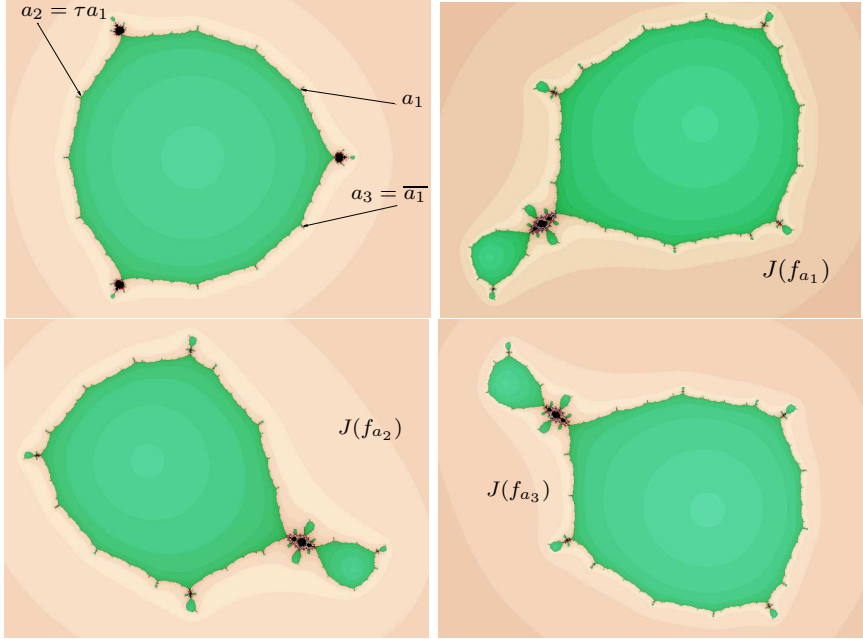
\vskip 0cm
\centerline{
\input{symmetries.pstex_t}  \input{a1.pstex_t} }
\centerline{\ \input{a2.pstex_t}  \input{a3.pstex_t} }
\caption{Symmetric parameters
$\mathbf{a}_1,\mathbf{a}_2,\mathbf{a}_3$ in $\CC_4$ and
$J(f_{\mathbf{a}_1})$, $J(f_{\mathbf{a}_2})$,
$J(f_{\mathbf{a}_3})$ . } \label{f:symmconnexite4}
\end{figure}
\begin{notation}\label{n:domainS+}
Let $\SS^+$ denote the connected component  of $ \C \setminus(
\tau^{-1}\R^{-}\cup\R^-)$ containing $\R^+$ for $d>3$, and
$\SS^+=\{z\mid \Im m(z)<0\}$ for $d=3$. Note that $\SS\subset
\SS^+$.
\end{notation}

\begin{remark}\label{r:t0(a)}
If $\mathbf{a}$ belongs to $\SS^+$, then
$\kappa_0(\mathbf{a})=\frac{1}{\tau
\lambda(\tau)}=e^{\frac{-2i\pi}{d-2}}$ and thus
$t_0(\mathbf{a})=-\frac{1}{d-2}$.
\end{remark}

\section{Coordinates in the parameter plane}\label{s:coordonees}

 The conformal representations  $\Upsilon_0$ and $\Upsilon_\infty$
are ``a priori'' independent of the dynamics.
 In this section we  define a {\it dynamical parametrization } of
 $\HH \cup \HH_\infty$ as well as {\it parameter rays} and {\it equipotentials}.

\subsection{``Dynamical'' parametrization of $\mathcal H_0$ and
$\HH_\infty$.}\label{s:H0}


 As usual, this
parametrization is given by the ``position'' of the critical
value. It is not  defined everywhere, but can be extended by
symmetry (see also~\cite{M4}).
  \begin{proposition}\label{p:conju} The following  map is
a holomorphic covering of degree~$d$\,: $$  \Phi_\infty \from \left\{
  \begin{aligned}
  \HH_\infty   & \longrightarrow \C\setminus \ol \D \\
   \mathbf{a} 
 & \longmapsto \Phi_\infty(\mathbf{a})= \phi^\infty_\mathbf{a}(f_\mathbf{a}(-\mathbf{a}))
  \end{aligned} \right. $$
Its restriction to $ \HH_\infty \cap \SS$  is a homeomorphism onto
$ \Delta_d$ where
$$\Delta_{d}=(-1)^{d-1}\left\{re^{i\theta}\ \left| r>1, \ \ 0\le
\theta\le\frac{1}{2}+\frac{1}{2(d-1)}\right.\right\}.$$
  \end{proposition}

 \proof   The B\"ottcher coordinate $\phi^\infty_\mathbf{a}(z)$ is holomorphic in
$(\mathbf{a},z)$  where it is defined (see~\cite{Bl}) so
$\Phi_\infty (\mathbf{a})$ is  holomorphic   on $\HH_\infty$ with
values
 in $\C \setminus\ol \D$. It extends
by ${\Phi_\infty(\infty)=\infty}$ and
  $\Phi_\infty \from\HH_\infty\cup \{\infty\} \to \ol\C\setminus \ol \D$
  is proper. Indeed,  for any sequence
  ${\mathbf{a}_n \in \HH_\infty}$  tending to
$\mathbf{a}\in\partial \mathcal{C}$,  the modulus
$\left|\Phi_\infty(\mathbf{a}_n)\right|
=e^{-G^\infty_{\mathbf{a}_n}(f_{\mathbf{a}_n}(-\mathbf{a}_n))}$
  tends to  $e^{-G^\infty_{\mathbf{a}}(f_{\mathbf{a}}(-\mathbf{a}))}$ since
the map $(\mathbf{a},x)\mapsto G^\infty_\mathbf{a}(x)$ is
continuous on $(\C\setminus\{0\})\times \C$.  Moreover,
$G^\infty_{\mathbf{a}}(f_{\mathbf{a}}(-\mathbf{a}))=0$ since
$f_\mathbf{a}(-\mathbf{a})\in K_\mathbf{a}$ for $\mathbf{a} \in
\mathcal{C}$, so $\left|\Phi_\infty(\mathbf{a}_n)\right|\to 1$.

The map $\Phi_\infty \from\HH_\infty\cup \{\infty\} \to
\ol\C\setminus \ol \D$ is holomorphic  and proper, hence it is  a
ramified covering whose degree is determined by the number of
preimages of $\infty$. Since $\infty$ is its sole preimage we have
to check the local degree at $\infty$. The following computations
show that $\Phi_\infty(\mathbf{a})\sim
\frac{-(-\mathbf{a})^d}{d-1}$ at $\infty$, so that the degree is
$d$. The function $\displaystyle
F_\mathbf{a}(z)=f_\mathbf{a}(z)/z^d=1+\frac{d \mathbf{a}}{(d-1)z}$
is greater than $1/2$ for $z \in D_\mathbf{a}=\{z \mid
|z|>|\mathbf{a}|^2\}$ and large~$\mathbf{a}$. Thus
$|f_\mathbf{a}(z)|\geq \frac{1}{2} |z|^d\ge |z|^2$ and
$f_\mathbf{a}$ sends
   $ D_\mathbf{a}$ into itself.
Hence, since $F_\mathbf{a}(D_\mathbf{a})$ is a small neighbourhood
of $1$, for every $k\ge 1$, the quantity $\displaystyle
(F_\mathbf{a}(f^k_\mathbf{a}(z)))^{\frac{1}{d^{k+1}}}$ is
well-defined on $D_\mathbf{a}$ by taking the principal
determination (since $\phi_\mathbf{a}^\infty(z)$ is tangent to
identity at $\infty$). The B\"ottcher coordinate
$\phi_\mathbf{a}^\infty$ is the limit of the functions
$$\phi^\infty_{\mathbf{a},k}(z) \left(f^k_\mathbf{a}(z)\right)^{1\over d^k}= z\,
\bigl(F_\mathbf{a}(z)\bigr)^{1\over d}\,
   \bigl(F_\mathbf{a}(f_\mathbf{a}(z))\bigr)^{1\over d^2} \,\cdots\,
   \bigl(F_\mathbf{a}(f^k_\mathbf{a}(z))\bigr)^{1\over d^{k+1}}.$$
For $k\ge 1$,
$\left|1-F_\mathbf{a}\bigl(f^k_\mathbf{a}(-\mathbf{a})\bigr)
\right| \leq \frac{d}{d-1} |\mathbf{a}|^{1-2^k}$ since for large
$\mathbf{a}$ the critical value $f_\mathbf{a}(-\mathbf{a})$
belongs to $D_\mathbf{a}$
 and $|f_\mathbf{a}^k(-\mathbf{a})|\ge |\mathbf{a}|^{2^k}$. Thus,
  $$ \sum_{k\ge 0} \frac{\log \left(F_\mathbf{a}(f^{k+1}_\mathbf{a}(-\mathbf{a}))\right)}{d^{k+1}
}\xrightarrow[n\to\infty] {} 0
 \quad \hbox{ and }
 \quad \displaystyle\Phi_\infty(\mathbf{a})\sim_\infty f_\mathbf{a}(-\mathbf{a})=\frac{-(-\mathbf{a})^d}{d-1}.$$

Hence $\Phi_\infty \from \HH_\infty \cup \{\infty\}\to \ol\C
\setminus \ol\D$ is a covering of degree $d$ which is ramified
only at $\infty$ (Riemann Hurwitz formula). We can
 lift it through the covering $\ol\C \setminus\ol \D \xrightarrow {z \mapsto z^d}
\ol\C \setminus\ol \D$ to a map $\upsilon \from\HH_\infty\cup
\{\infty\} \to \ol\C \setminus\ol \D$ satisfying
$\upsilon(\mathbf{a})^d \sim \frac{-(-\mathbf{a})^d}{d-1}$  at
$\infty$, so that we can choose $\upsilon$ tangent to $-\alpha
e^{\frac{i\pi}{d}}\hbox{Id}$  at $\infty$ with $\alpha>0$.
Therefore  $\upsilon(\mathbf{a})=-e^{i
\pi/d}\Upsilon_{\infty}(\mathbf{a})/(d-1)^{1/d}$  (they are
conformal representations from $\HH_\infty\cup \{\infty\}$ to
$\ol\C\setminus \ol \D$ tangent  at $\infty$, so they coincide).
Hence $$\Phi_\infty
(\mathbf{a})=(-1)^{d-1}\Upsilon_{\infty}(\mathbf{a})^d/(d-1).$$
This determines the image of $\SS \cap \HH_\infty$ by
$\Phi_\infty$ using Corollary~\ref{c:domaineimage}. Indeed,
$\Phi_\infty(\R^+\cap \HH_\infty)=(-1)^{d-1}\R^+\setminus \ol \D$
and $\Phi_\infty(\rho \R^+\cap \HH_\infty)=(-1)^d\rho\R^+\setminus
\ol \D$ so that $\Phi_\infty(\HH_\infty\cap \SS)= \Delta_d$
 since $\Phi_\infty$ preserves the cyclic order at
$\infty$. Finally  $\Phi_\infty$ is  injective on $\SS \cap
\HH_\infty $ because the opening  of
$\Upsilon_\infty(\mathcal{S})$ is less than $1/d$ and $\Phi_\infty
(\mathbf{a})=(-1)^{d-1}\Upsilon_{\infty}(\mathbf{a})^d/(d-1)$.
 \cqfd

\begin{figure}[!h]
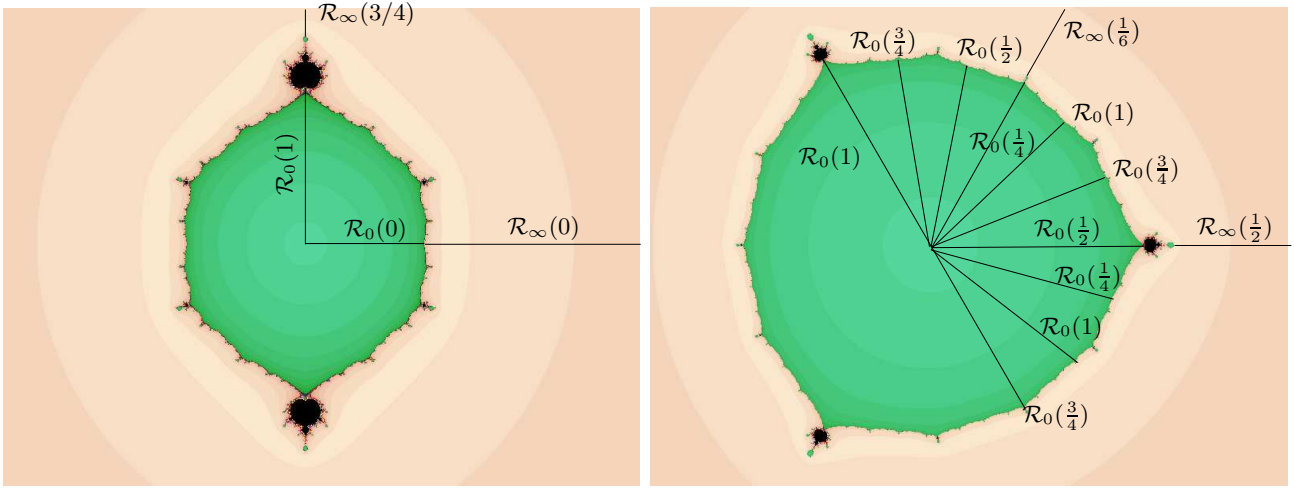
\vskip 0cm
\centerline{\input{rayonsC3.pstex_t}
\input{rayonsC4.pstex_t}}
\caption{Some rays in $\HH_0$ and $\HH_\infty$ for $d=3$ and $d=4$.}\label{f:rays34}
  \end{figure}

  \begin{remark}For $d=3,4$ one has $$\Delta_3=\left\{]1,\infty[e^{2i\pi\theta}\mid 0\le
\theta\le\frac{3}{4}\right\},
\Delta_4=\left\{]1,\infty[e^{2i\pi\theta}\mid \frac{1}{2}\le
\theta\le \frac{1}{6}+1\right\}.$$
  \end{remark}
\begin{remark}\label{r:symrayinf}
>From the proof above and Remark~\ref{r:symetry}, the following
symmetries hold\,:
$\Phi_\infty(\sigma(\mathbf{a}))=\sigma(\Phi_\infty(\mathbf{a}))$
and $\Phi_\infty(\tau \mathbf{a})=\tau \Phi_\infty(\mathbf{a})$
where $\tau=e^\frac{2i\pi}{d-1}$.
\end{remark}

\begin{proposition}\label{p:H0}
The map $\Phi_0(\mathbf{a})=
\phi^0_\mathbf{a}(f_\mathbf{a}(-\mathbf{a}))$ is well-defined on
$\HH_0\setminus \R^-$.
\begin{itemize}
\item For $d>3$, its restriction to $ \HH_0\cap \mathcal S $ is a
holomorphic homeomorphism onto
$$\Lambda_{d}=(-1)^{d-1}\left\{re^{i\theta}
 \mid r \in [0,1[,\  0\le
 \theta\le\frac{1}{2}+\frac{1}{2(d-2)}\right\}$$
 \item For $d=3$, the restriction of  $\Phi_0$ to $\HH_0\cap \dot \SS$
 is a holomorphic  homeomorphism onto 
 $\D\setminus \R^+$,
 where $\dot \SS$ denotes  the interior of $\SS$.
 Moreover,  $\Phi_0$ maps bijectively   each of the two boundary lines
 $\HH_0\cap\R^+$ and  $\HH_0\cap i\R^+$ onto $[0,1[$   {\rm(}with $i^2=-1${\rm)}.
\end{itemize} \end{proposition}

 \proof  The proof goes exactly as in Proposition~\ref{p:conju}\,;
 we just have to take care since $\phi_\mathbf{a}^0$ is not defined for
 $\mathbf{a}\in \R^-$ (see B\"ottcher's Theorem, section~\ref{s:dyn}). We do it  for $d>3$, the
 proof for $d=3$ follows  the same arguments\,; the only difference is that
 the boundaries of $\Lambda_d$ glue together to give the disk for $d=3$,
 so if we forget the components of the boundary  the arguments go through for $d=3$.
  
The map $\Phi_0$ is holomorphic on $\mathcal H_0  \setminus \R^-$.
Indeed, for $\mathbf{a} \in \mathcal H_0  \setminus (\R^-\cup\{
0\})$
 $\Phi_0(\mathbf{a})=\phi_\mathbf{a}(\lambda(\mathbf{a})f_\mathbf{a}(-\mathbf{a}))$,
 where $\phi_\mathbf{a}$ is the B\"ottcher
 coordinate near $0$ of the map
$
g_\mathbf{a}(z)=\lambda(\mathbf{a})f_\mathbf{a}\left(z/\lambda(\mathbf{a})\right)=
z^{d-1}+\lambda(\mathbf{a})^{1-d}\ z^d$ (remark that
$\lambda(\mathbf{a})$ is a non-vanishing holomorphic function on
$\HH_0  \setminus (\R^-\cup\{ 0\})$). As in
Proposition~\ref{p:conju} above, $\phi_\mathbf{a}$ is the limit of
$$\phi_{\mathbf{a},k}(z) = \left(g^k_\mathbf{a}(z)\right)^{1\over
(d-1)^k}= z\, \bigl(G_\mathbf{a}(z)\bigr)^{1\over d-1}\,
   \,\cdots\,
   \bigl(G_\mathbf{a}(g^k_\mathbf{a}(z))\bigr)^{1\over (d-1)^{k+1}}$$
   where  $\displaystyle G_\mathbf{a}(z)=g_\mathbf{a}(z)/z^{d-1}=1+\lambda(\mathbf{a})^{1-d}\ z$.
The extension of $\phi_\mathbf{a}$, using $\phi_\mathbf{a}\circ
g_\mathbf{a}=(\phi_\mathbf{a})^{d-1}$
 gives a holomorphic function
$(\mathbf{a},z)\to \phi_\mathbf{a}(z)$  with
$(\phi_\mathbf{a})'(0)=1$.

It remains to prove that $\Phi_0$ is proper. As in the proof of
Proposition~\ref{p:conju}, if  $\mathbf{a}\in \HH_0  \setminus
(\R^-\cup\{ 0\})$ tends to $\partial \HH_0$, $\Phi_0(\mathbf{a})$
tends to $\partial \D$\,; an analogous computation shows that
 $\Phi_0(\mathbf{a})\sim -\lambda(\mathbf{a})(-\mathbf{a})^d /(d-1)$  near $0$. So we can extend
 $\Phi_0$  by $\Phi_0(0)=0$.

 In order to determine the images of $\R^+\cap \HH_0$ and $\rho\R^+\cap
  \HH_0$,
 with $\rho=e^{i\pi/(d-1)}$, we can use  that $\phi_\mathbf{a}^0$ is defined
 by continuity at the critical point $-\mathbf{a}$ and satisfies
 $\Phi_0(\mathbf{a})=\phi^0_\mathbf{a}(f_\mathbf{a}(-\mathbf{a}))=(\phi^0_\mathbf{a}(-\mathbf{a}))^{d-1}$.
 From the proof of  Remark~\ref{l:symray},
 $\phi^0_\mathbf{a}(\R^-)\subset \R^-$  for $\mathbf{a} \in \R^+$ since
  $(\phi^0_\mathbf{a})'(0)\in \R^+$ and $\phi^0_\mathbf{a}(\sigma(z))=\sigma(\phi^0_\mathbf{a}( z))$.
 Hence $\Phi_0(\R^+)\subset(-1)^{d-1}\R^+$.
 Now for $\mathbf{a} \in \rho \R^+ \cap \HH_0$, we  determine $t$ such that $-\mathbf{a}$ ``is''
  on the  ray $R_\mathbf{a}^0(t)$, using the symmetries of Lemma~\ref{l:symray}
   with $\tau=\rho^2$.
 Since $\tau\sigma(\mathbf{a})=\mathbf{a}$, the critical point $-\mathbf{a}$ is on
 $R_\mathbf{a}^0(t)=R_{\tau \sigma (\mathbf{a})}^0(t)$, thus
 $-\sigma(\mathbf{a})=-\tau^{-1}\mathbf{a}\in \tau^{-1}R_{\tau \sigma (\mathbf{a})}^0(t)$.
 Since $\tau^{-1}R_{\tau \sigma (\mathbf{a})}^0(t)= R^0_{\sigma \mathbf{a}}(t+t_0)
 =\sigma(R^0_\mathbf{a}(-t-t_0))$,
 where $t_0=t_0(\sigma(\mathbf{a}))=\frac{-1}{d-2}$
 (Remark~\ref{r:t0(a)}),  we deduce that the critical point $-\mathbf{a}$
 is  also on  $R^0_\mathbf{a}\left(-t+\frac{1}{d-2}\right)$\,; so   $t=-t+\frac{1}{d-2} \mod 1$
 because there is a unique ray stemming from $0$ which contains the critical
 point (the case of bounded Fatou components).
 Hence, $2t=\frac{1}{d-2} \mod 1$ and the critical point belongs to
$R^0_\mathbf{a}\left(\frac{1}{2(d-2)}\right)$ or to
$R_\mathbf{a}^0\left(\frac{1}{2(d-2)}+\frac{1}{2}\right)$. The map
$\phi_\mathbf{a}^0 (z)$ is asymptotic to $\lambda(\mathbf{a})z$
near $0$\,; so any point $z\in \rho \R^-$ near $0$ is sent inside
$e^{i\pi/(d-2)}\R^-$, since $\mathbf{a} \in \rho \R^+$. Finally,
$\phi_\mathbf{a}^0(-\mathbf{a})\in e^{i\pi/(d-2)}\R^-$ and
$\Phi_0(\mathbf{a})\in(-1)^{d}e^{i\pi/(d-2)}\R^+$.

The image $\Phi_0(\SS\cap \HH_0)$ is exactly $\Lambda_d$. Indeed,
 the points of $\SS$ near $0$ are sent by $\Phi_0$ inside $\Lambda_d$.
Moreover, if $\Phi_0(\SS\cap \HH_0)$ were bigger than $
\Lambda_d$, there would be another connected component of
$\Phi_0^{-1}(\partial \Lambda_d)$ in $\SS$, but then
$\Phi_0^{-1}(0)\cap \SS\neq 0$ which is impossible
($\phi^0_\mathbf{a}(f_\mathbf{a}(-\mathbf{a}))=0 \implies
\mathbf{a}=0$).
 Finally $\Phi_0$ is a proper holomorphic map from $ \SS \cap \HH_0$ onto
 $\Lambda_d $. It is a ramified covering since $\HH_0\cap \SS$ is
 simply connected (by Lemma~\ref{l:1connected} and
 Corollary~\ref{c:domaineimage}).
The covering  $\Phi_0\circ \Upsilon_0^{-1}\from \Upsilon_0(\SS\cap
\HH_0) \to \Lambda_d $ extends
 to all the boundary and its degree is  the number of  preimages  of $0$.
Thus $\Phi_0\from \SS \cap \HH_0 \to \Lambda_d$ is a
 holomorphic homeomorphism.
 \cqfd

\begin{remark}\label{r:symray0}
>From  Lemma~\ref{l:symray}, Remark~\ref{r:t0(a)} and the proof
above, the following symmetries hold\,:
$\Phi_0(\sigma(\mathbf{a}))=\sigma(\Phi_0(\mathbf{a}))$ for
$\mathbf{a}\notin \R^-$ and, for $\mathbf{a} \in \SS^+$,
$\Phi_0(\tau \mathbf{a})=e^\frac{2i\pi}{d-2} \Phi_0(\mathbf{a})$,
where $\tau=e^\frac{2i\pi}{d-1}$ and
$\sigma(\mathbf{a})=\ol{\mathbf{a}}$.
\end{remark}

\begin{remark}
The difference (in the angles) between $\Lambda _d$ and $\Delta_d$
comes from the fact that the  boundary of  $\HH_\infty$ will also touch
the boundary of the components of $\HH_n$ with $n>0$.
\end{remark}

\subsection{Parametrization of $\mathcal H_n$. }\label{s:Hi}

\begin{lemma}\label{l:UinterS} If $\mathcal U$ is a connected component of
$\mathcal H_n$ with $n>0$, it cannot contain at the same time a
point $\mathbf{a}_0$ and  its symmetric
 $\tau \mathbf{a}_0$. Therefore, either $\mathcal U$ or $\tau \UU$
  is included in $\C \setminus \R^-$.
\end{lemma}

\proof Suppose by contradiction that $\UU$ contains a point
$\mathbf{a}_0$ and its symmetric $\tau \mathbf{a}_0$. Then, we can
construct
 a curve $\gamma$ surrounding $0$ on which $f_\mathbf{a}^n$ is uniformly bounded for every
 $n$\,: $\gamma$  is the union of an arc in $\UU$ joining $\mathbf{a}_0$
 to $\tau \mathbf{a}_0$ and
 all the symmetric images by  $\tau^k$. This curve  $\gamma$ is contained in $\UU$
 since  $\UU\cap \tau\UU\neq \emptyset $, so that $\UU=\tau \UU=...=\tau^k\UU$.
Then there are points  of $\HH_\infty$ surrounded by $\gamma$
since $\ol\HH_0$ is surrounded by $\gamma$ ($0\in \HH_0$ and
$\gamma\subset \UU\neq \HH_0$) and $\partial \HH_0\subset\partial
\HH_\infty$. This contradicts the maximum principle exactly as in
Lemma~\ref{l:1connected}.

 Suppose now that $\UU$ crosses $\R^-$ and also $\tau^{-1}\R^-$\,; then it will
necessary  cross $\rho^{-1}\R^{-}$ or $\rho\R^{-}$where
$\rho=e^{i\pi/{d-1}}$.
 Moreover $\UU=\sigma(\UU)$ since both have a common point on $\R^-$. Let
 $\mathbf{a}_0$ be some point of $\UU\cap \rho^{-1}\R^-$.
 Then $\sigma(\mathbf{a}_0)$ belongs to $\UU$.
 But $\sigma(\mathbf{a}_0)=\tau \mathbf{a}_0$  since $\tau\rho^{-1}=\sigma
 (\rho^{-1})$.  By the first part of this Lemma, this is again a contradiction.
 \cqfd

\begin{lemma}\label{l:parametrizationU}
Let $\mathcal U$ be a connected component of $\mathcal H_n$ with
$n>0$ included in $\C \setminus \R^-$. The map $$ \Phi_{\mathcal
U} \from \left\{
\begin{aligned}
 \mathcal{U} & \longrightarrow  \D  \\
 \mathbf{a} 
& \longmapsto \Phi_{\mathcal U}(\mathbf{a})=
\phi^0_\mathbf{a}(f_\mathbf{a}^n(f_\mathbf{a}(-\mathbf{a})))
 \end{aligned} \right. $$
is a conformal homeomorphism.
\end{lemma}
\proof For $\mathbf{a}\in \mathcal U$, the map $\phi_\mathbf{a}^0$
is well-defined, and is holomorphic in $(\mathbf{a},z)$ for $z \in
B_\mathbf{a}$ (see Remark~\ref{r:bottcher}). Hence since $\UU$ is
simply connected $\Phi_{\mathcal U}$ is a ramified covering from
$\UU$ to $\D$ (it is holomorphic and proper\,: the proof is
similar to that of Proposition~\ref{p:H0}). It remains to show
that it has degree one. For this we will prove that $\Phi_{\UU}$
is a local homeomorphism near every point of $\mathcal U$. Let
$\mathbf{a}_0\in \mathcal U$ and $z_0=\Phi_{\mathcal
U}(\mathbf{a}_0)$. We will construct by surgery a local inverse
$\mathbf{a}(z)$ to $z=\Phi_{\mathcal U}(\mathbf{a})$ in a
neighbourhood of $z_0$. We first modify $f_{\mathbf{a}_0}$ near
the point $f_{\mathbf{a}_0}^{n}(-\mathbf{a}_0)$, in order to send
$f_{\mathbf{a}_0}^{n}(-\mathbf{a}_0)$ to
$(\phi_{\mathbf{a}_0}^0)^{-1}(z)$ and we prove that this new map
is quasi-conformally conjugated to $f_{\mathbf{a}(z)}$, so that
$\Phi_{\mathcal U}(\mathbf{a}(z))=z$.

Let $B'=f_{\mathbf{a}_0}^{-1}(B_{\mathbf{a}_0})$, $B'$ is a
topological disc. For $\epsilon$ small enough and for every $z\in
D(z_0,\epsilon)$ one can construct a map $\delta_z\from \ol\C \to
\ol\C$ smooth in the variable $z$ and satisfying\,:

$\bullet$ $\delta_{z_0}=f_{\mathbf{a}_0} $,  

$\bullet$
$\delta_z(f_{\mathbf{a}_0}^{n}(-\mathbf{a}_0))=(\phi_{\mathbf{a}_0}^0)^{-1}(z)$,

$\bullet$  $\delta_z$ coincides with $f_{\mathbf{a}_0}$ outside
$V$  which is a small  neighbourhood of
$f_{\mathbf{a}_0}^{n}(-\mathbf{a}_0)$ compactly included in $B'$
(and independent of $z$),

$\bullet$ $\delta_z$ is a diffeomorphism  from $B'$ to
$B_{\mathbf{a}_0}$.

\noindent  We denote by $\sigma_z$ the complex structure which
coincides with the standard complex structure on
$B_{\mathbf{a}_0}\cup (\C\setminus \wt B_{\mathbf{a}_0})$ and
which is invariant by $\delta_z$. This complex structure has
bounded dilatation. Let $g_z$ be the homeomorphism that integrates
$\sigma_z$, given by Ahlfors-Bers' Theorem which is normalized to
fix $0$ and to  be tangent  to identity at $\infty$. Then the map
$h_z=g_z\circ \delta_z\circ g_z^{-1}$ is holomorphic. Moreover it
fixes $\infty$ with local degree $d$ and $0$ with local degree
$d-1$. Therefore,
$h_z(u)=u^{d-1}(u+d\mathbf{a}(z)/(d-1))=f_{\mathbf{a}(z)}(u)$
 with  $\mathbf{a}(z)$ a continuous function of $z$ and $\mathbf{a}(z_0)=\mathbf{a}_0$.
 On $B_{\mathbf{a}_0}$, $g_z$ is  holomorphic,  conjugates
 $f_{\mathbf{a}_0}$ to  $f_{\mathbf{a}(z)}$\,; so  $\phi_{\mathbf{a}_0}^0$ and $\phi_{\mathbf{a}(z)}^0\circ g_z$
differ from a $d-2$ root of unity (since they both conjugate
$f_{\mathbf{a}_0}$ to $z\mapsto z^{d-1}$). For $z=z_0$ this root
of unity is $1$\,; so, by continuity,
$\phi_{\mathbf{a}_0}^0=\phi_{\mathbf{a}(z)}^0\circ g_z$. Moreover
$g_z(-\mathbf{a}_0)=-\mathbf{a}(z)$ since $g_z$ preserves the
critical points (looking at the local degree). This implies that
$\Phi_{\mathcal
U}(\mathbf{a}(z))=\phi_{\mathbf{a}(z)}^0(f^{n+1}_{\mathbf{a}(z)}(-\mathbf{a}(z)))=z$
since $f^{n+1}_{\mathbf{a}(z)}(-\mathbf{a}(z))=g_z\circ
\delta^{n+1}_z\circ g_z^{-1}(-\mathbf{a}(z))=g_z\circ
\delta_z^{n+1}(-\mathbf{a}_0)= g_z\circ
\delta_z(f_{\mathbf{a}_0}^{n}(-\mathbf{a}_0))$ and $g_z\circ
\delta_z(f_{\mathbf{a}_0}^{n}(-\mathbf{a}_0))=g_z\circ
(\phi_{\mathbf{a}_0}^0)^{-1}(z)=(\phi^0_{\mathbf{a}(z)})^{-1}(z).$
\cqfd

\begin{remark}From Lemma~\ref{l:symray} and Remark~\ref{r:t0(a)},
the following symmetries hold\,: for $\mathbf{a}\in \UU$,
$\Phi_{\sigma(
\UU)}(\sigma(\mathbf{a}))=\sigma(\Phi_\UU(\mathbf{a}))$ and,
$\Phi_{\tau\UU}(\tau \mathbf{a})=e^\frac{2i\pi}{d-2}
\Phi_\UU(\mathbf{a})$ if  $\UU \subset \SS^+$, where
$\tau=e^\frac{2i\pi}{d-1}$.
\end{remark}


\subsection{Parameter rays and equipotentials in $\C\setminus \R^-$.}

 \begin{definition}
   We define the  {\it equipotential } of level $v>0$, in $\HH_\infty$ and  in $\HH_0$,
   by\,:
   $$\mathcal{E}_\infty(v)=\Phi_\infty^{-1}\left(\left\{e^{ v+2 i \pi t},
   \ t\in [0,1]\right\}\right),
   \quad  \mathcal{E}_0(v)=\Phi_0^{-1}\left(\left\{e^{-v+2 i \pi t},
   \ t\in [0,1]\right\}\right).$$
\end{definition}
Note that an equipotential in $\HH_\infty$ is a closed curve
surrounding $\CC$ and that  an equipotential in $\HH_0 $ is never
closed because the point in $\R^-$ is missing. One can close it
however by adding this point.
  
 \begin{definition}\label{d:pararays}Let   $p \in \{ 0, \infty\}$, we define the union
 of the rays of angle $t$ in $ \HH_p$,  by $$U\RR_p(t)=(\Phi_p)^{-1}(\R^+ e^{2i\pi t}).$$

\end{definition}

 \begin{remark}
With this definition, there is no intersection between  $\R^-$ and
$U\RR_0(t)$ for any $t$.
\end{remark}

 \begin{definition}
Let  $p \in\{0,\infty\}$ and  $s\in \{\hbox{Id},
\tau,\sigma,\tau\sigma,\cdots\}$ the result of  any composition of
the symmetries
 $\sigma$ and $\tau$. Denote by $\Phi^{s}_p$ the
 restriction to $s(\SS)$ of the map $\Phi_p$
 and define the  {\it ray}
 of angle $t$  in $\HH_p\cap s(\SS)$ by
  $$ \mathcal{R}^s_p(t)=\left(\Phi^{s}_p\right)^{-1}\left(\left\{re^{2 i \pi
t},\ r\in \R^+\right\}\right)=U\RR_p(t)\cap s(\SS) .$$
\end{definition}

\begin{remark}\
\begin{itemize}
\item

The set $\RR^s_p(t)\setminus\{0\}$ is connected except when $d=3$,
$t=0$ and $s=\id$\,; \item For $d=3$ and $t=0$,
$\RR^{\hbox{Id}}_0(0)\setminus\{0\}$ has  two connected
components. Therefore we will make the convention that
$\RR_0^{\id}(0)$ is the ray $\R^+\cap \HH_0$ and $i \R^+\cap
\HH_0$ is the ray $\RR_0^{\tau \sigma}(0)$\,;

\item We will get rid of these notations when we will work in a
para-puzzle piece $\PP_0(\mathbf{a}_0)$ since from
Remark~\ref{r:inclusdsS}, for every $t$, $U\RR_p(t)\cap
\PP_0(\mathbf{a}_0)$ has  only  one connected component that we
will call $\RR_p(t)$.
\end{itemize}\end{remark}

\begin{lemma}\label{l:raysym} We have
$\sigma(\RR^{s}_p(t))=\RR^{\sigma s}_p(-t), \
\tau\RR^{s}_\infty(t)=\RR^{\tau s}_\infty(t+\frac{1}{d-1})$
 and if
 $\RR^{s}_0(t) \subset  \SS^+$, $\tau\RR^{s}_0(t)=\RR^{\tau s}_0(t+\frac{1}{d-2})$ .
 \end{lemma}
\proof This follows from Remarks~\ref{r:symray0}
and~\ref{r:symrayinf} \cqfd
 
 It follows from Lemma~\ref{l:raysym} that  $U\RR_p(t)$ is not invariant  in general
 by $\sigma$ neither by $\tau$.
\begin{remark}\label{r:rayons0} We have the following
correspondences for parameters $\mathbf{a}\in \C \setminus \R^-$
and $p\in\{0,\infty\}$\,:
 $$(\mathbf{a}\in U\RR_p(t)\iff f_\mathbf{a}(-\mathbf{a})\in R^p_\mathbf{a}(t))\ \hbox{ and
 }\  (\mathbf{a}\in \mathcal{E}_p(v) \iff f_\mathbf{a}(-\mathbf{a})\in E^p_\mathbf{a}(v)).
 $$
\end{remark}

\begin{remark}\label{r:rayons}
Let $p\in\{0,\infty\}$. \begin{itemize} \item The line $\R^+\cap
\HH_p$ is the ray $\RR^{\id}_p(0)$ if $d$ is odd, and
$\RR^{\id}_p(1/2)$  if $d$ is  even\,; \item The line
$\rho\R^+\cap \HH_0$ and $\rho\R^+\cap \HH_\infty$ are the rays
$\RR^{\id}_0\left(\frac{1}{2}+\frac{1}{2(d-2)}\right)$ and
$\RR^{\id}_\infty\left(\frac{1}{2}+\frac{1}{2(d-1)}\right)$
respectively if $d>3$ is odd, and
$\RR^{\id}_0\left(\frac{1}{2(d-2)}\right)$ and
$\RR^{\id}_\infty\left(\frac{1}{2(d-1)}\right)$ respectively if
$d$ is even.
\end{itemize}\end{remark}

\proof This follows from Proposition~\ref{p:H0}
and Proposition~\ref{p:conju} (those angles are given by the
boundaries of $\Lambda_d$ and $\Delta_d$).

\begin{definition}
Let $\UU\subset \C \setminus \R^{-}$ be a connected component of
$\HH_n$, with $n>0$.  We define the {\it center} of $\UU$ by
$c_{\UU}=\Phi_\UU^{-1}(0)$ and  the {\it internal ray} of angle
$t\in \R/\Z$, resp. the {\it equipotential} of level $v$ by\,:
$$ \mathcal{R}_\UU(t)=\Phi_\UU^{-1}\left(\left\{re^{2 i \pi t},\
r<1\right\}\right),\ \hbox { resp. }\
\mathcal{E}_\UU(v)=\Phi_\UU^{-1}\left(\left\{e^{-v+2 i \pi t},\
t\in [0,1]\right\}\right).$$
\end{definition}
Note that $\RR_\UU(t)\setminus\{c_\UU\}$ is connected  since
$\Phi_\UU$ is a homeomorphism.
 
\begin{remark}For $\mathbf{a} \in \HH_n\setminus \R^-$ with $n>0$, the critical value
$f_\mathbf{a}(-\mathbf{a})$ belongs  to a connected component
$U_\mathbf{a}$ of $\wt B_\mathbf{a}$. The map $f^n_\mathbf{a}\from
U_\mathbf{a}\to B_\mathbf{a}$ is a homeomorphism. Thus we can pull
back the B\"ottcher coordinate to get coordinates
on~$U_\mathbf{a}$.
\end{remark}
\begin{notation}
We denote by $r$ and call the {\it center} of $U_\mathbf{a}$ the
unique point of $U_\mathbf{a}$ which is sent by $f^n_\mathbf{a}$
to $0$ (the center of $B_\mathbf{a}$). We denote by
$R^r_\mathbf{a}(t)$ the ray stemming from $r$ with B\"ottcher
coordinate $t$, {\it i.e.} the preimage
$(f^n_\mathbf{a}|_{U_\mathbf{a}})^{-1}(R^0_\mathbf{a}(t))$ which
contains $r$ in its closure.
\end{notation}
As in Lemma~\ref{l:aboutdistray0} we get now the criterium  for
connected components in the dynamical plane  to have a common
boundary point which is   the landing point of rays.
\begin{lemma}\label{l:aboutdistray1}
Let $\mathbf{a}\in \C\setminus \R^-$ and $U,V$  be two connected
components of $\wt B_\mathbf{a}$ with center $r,r'$ respectively.
If the rays $R_\mathbf{a}^r(t)$ and $R_\mathbf{a}^{r'}(t')$ land
at the same point $x$ then either $r=r'$ and $t=t'$ or the landing
point $x$ is eventually critical.
\end{lemma}
\proof We assume that the rays $R^r_\mathbf{a}(t)$ and
$R^{r'}_{\mathbf{a}}(t')$ are different. After several iterations
by $f_\mathbf{a}$, the image of the two rays in $B_\mathbf{a}$
should coincide by Lemma~\ref{l:aboutdistray0}. Therefore we have,
at some  step of the iterations, the situation of two rays landing
at a common point and  having the same image under the map
$f_\mathbf{a}$. Then the common landing point is the critical
point $-\mathbf{a}$. \cqfd

 Note that when $\mathbf{a}$ is the center of a
connected component
 $\UU$ of $\HH_n$, the critical value is the center of $U_\mathbf{a}$.

\begin{lemma}\label{l:defcentre}
Let $\UU$  be a connected component of $\HH_n$ with $n>0$. Let
$\Omega$ be a simply connected neighbourhood of $\ol\UU$ which
avoids the centers of all the components of  $\HH_j$ with $j< n$.
There exists on $\Omega$ a holomorphic map $r=r_\UU\from \Omega\to
\C$ such that for $\mathbf{a}\in \UU$, $r(\mathbf{a})$ is the
center of the connected component that contains the critical
value.
\end{lemma}
\proof We apply to $F(\mathbf{a},z)=f_{\mathbf{a}}^n(z)$ the
Implicit Function Theorem  in a neighbourhood of the point
$(\mathbf{a},z)=(c_\UU, f_{c_\UU}(-c_\UU))$. The only point where
it is not possible to apply the Theorem are the centers of $\HH_j$
for $j< n$ since then the critical point is sent after $n$
iterations to~$0$. \cqfd

\begin{corollary}\label{c:corresp}
Let $\UU\subset \C\setminus \R^-$ be a connected component of
$\HH_n$ with $n>0$. We have the following equivalence\,: $
\mathbf{a} \in \RR_\UU(t) \iff f_\mathbf{a}(-\mathbf{a})\in R^{r
(\mathbf{a})}_\mathbf{a}(t).$
\end{corollary}


\subsection{Landing properties}

Most of the results in this subsection  follow from the classical
case of quadratic polynomials, see~\cite{DH1, M1} and
also~\cite{M4}. Recall that for $\UU=\HH_\infty$ or $\HH_0$, the
ray $\RR^s_\UU(t)$
 is the one in $s(\SS)$ where $s$ is any composition of $\sigma$ and $\tau$.
  Thus, it is not defined for  any  $t$\,; for instance for $s=\id$, the  angle
  has to be in  $\Delta_d$ or in $\Lambda_d$.
 For a more detailed description  and an other proof of the following Lemma,
see section~\ref{s:annexe1}.
\begin{lemma}\label{l:misiurpreim1} Let $\UU\subset \C \setminus \R^-$
be a connected component of $\HH_n$ with $n\in \N\cup \{\infty\}$.
For $t$ rational, the ray   $\RR^s_\UU(t)$ converges. Let
$\mathbf{a}_0$ denote the landing point. If $\mathbf{a}_0\notin
\R^-$, the ray $R^{r(\mathbf{a}_0)}_{\mathbf{a}_0}(t)$ is periodic
{\rm(}resp. eventually periodic{\rm)}
 and lands at a
 parabolic periodic  {\rm(}resp. eventually
periodic{\rm)} point or at $f_{\mathbf{a}_0}(-\mathbf{a}_0)$ that
is a
 repelling periodic {\rm(}resp. eventually
periodic{\rm)} point.
\end{lemma}
In the last sentence, ``eventually'' depends on $t$, meaning  that
the number of iterates after which
$R_{r(\mathbf{a}_0)}^{\mathbf{a}_0}(t)$ becomes periodic, and the
period, both depend on $t$.

\proof Assume that $\UU=\HH_\infty$, the proof being easier for
the other components. Let $\mathbf{a}_0$ be an accumulation point
of $\RR^s_\infty(t)$. Since $\mathbf{a}_0\in \CC$ and $t \in \Q$
the ray $R^\infty_{\mathbf{a}_0}(t)$ is well-defined and
converges. The landing point is (eventually) periodic, either
parabolic or repelling. If it is repelling, and not eventually the
critical value, by lemma~\ref{l:douady1} we should have the
stability of the rays $R^\infty_{\mathbf{a}_0}(t/d+k/d)$. But for
parameters $\mathbf{a}$ on $\RR^s_\infty(t)$ near $\mathbf{a}_0$
the critical value is on $R^\infty_{\mathbf{a}}(t)$ (by
definition), so at least two of the previous   rays crash on the
critical point. So there exist $p,l$ depending only on $t$ such
that
$f_{\mathbf{a}_0}^{p+l}(-\mathbf{a}_0)=f_{\mathbf{a}_0}^l(-\mathbf{a}_0)$.
If the landing point is eventually parabolic, the resultant of the
two polynomials $f_{\mathbf{a}_0}^p (y) -y$ and
$(f_{\mathbf{a}_0}^p)' (y) -1$ vanishes since they have at least a
common root. In both cases $\mathbf{a}_0$ is a root of a
polynomial. Hence the accumulation set is a finite set and, since
it is connected, it reduces to a point.

In the case where the landing point of
$R^\infty_{\mathbf{a}_0}(t)$ is eventually repelling, let $k$ be
the first integer such that the critical value lies on
$\ol{R^\infty_{\mathbf{a}_0}(d^kt)}$. Thus the compact set $\ol
{R^\infty_{\mathbf{a}_0}(d^kt)}$ moves continuously (by $\ol
{R^\infty_{\mathbf{a}}(d^kt)}$) for $\mathbf{a}$ in a
neighbourhood $\Omega$ of $\mathbf{a}_0$. In particular, for
$\mathbf{a}\in (\RR^s_\infty(t)\cup \{\mathbf{a}_0\})\cap \Omega$,
the compact set $\ol{R^\infty_{\mathbf{a}}(t)}$ is a continuous
image of $\ol {R^\infty_{\mathbf{a}}(d^kt)}$  since
$f_\mathbf{a}^k$ restricts
 on ${R^\infty_{\mathbf{a}}(t)}$ to a  homeomorphism onto  $R^\infty_{\mathbf{a}}(d^kt)$
 (the holomorphic motion of the closure is obtained using  the $\lambda$-Lemma).
Therefore the critical value, which  for $\mathbf{a} \in
\RR^s_\infty(t)\cap\Omega$ is in $R^\infty_{\mathbf{a}}(t)$,  is
for $\mathbf{a}=\mathbf{a}_0$ in
$\ol{R^\infty_{\mathbf{a}_0}(t)}$. Therefore
$R^\infty_{\mathbf{a}_0}(t)$ lands at
$f_{\mathbf{a}_0}(-\mathbf{a}_0)$. \cqfd

\begin{definition}
A parameter $\mathbf{a}$ is  {\it Misiurewicz} (or of  {\it
Misiurewicz type})
 if for some $l\ge 1$, $z=f_\mathbf{a}^l(-\mathbf{a})$ is a
periodic point of  $f_\mathbf{a}$.
\end{definition}
Note  that if $f_\mathbf{a}^l(-\mathbf{a})$ is periodic (with
$l\ge 1$), it is necessarily  a repelling point. Indeed, if it is
 attracting or parabolic  it would
 attract a critical point
and  there is no other ``free'' critical point that can converge
to it. Note also that all the Misiurewicz points are in $\CC$.

\begin{lemma}\label{l:misiurewicz}
Let $\mathbf{a}\in \C \setminus  \R^-$  be a Misiurewicz point.
There exists $t\in \Q$  such that $R_\mathbf{a}^{\infty}(t)$ lands
at $ f_\mathbf{a}(- \mathbf{a}) $. Moreover,  the ray
$\RR^s_\infty(t)$ lands at $\mathbf{a}$, for $s$ such that
$\mathbf{a}\in s(\SS)$.
\end{lemma}
\proof  The proof is exactly the same as in~\cite{DH1}.
\cqfd

\begin{corollary}\label{c:R+}
If $d$ is odd, $\R^+=\RR^{\id}_\infty(0)\cup \RR^{\id}_0(0)\cup
\{*\}$ where $*$ is a Misiurewicz point.
 \end{corollary}
 \proof  If $d=2l+1$, $\R^+ \cap \HH_0=\RR^{\id}_0(0)$
(Remark~\ref{r:rayons}). Let $\mathbf{a}_0$ be the landing point
of the ray $\RR^{\id}_0(0)$, $\mathbf{a}_0\in \R^{+*}$. The ray
$R_{\mathbf{a}_0}^0(0)\subset \R^+$ lands at a fixed point, say
$x_0\in \R^{+*}$. If this fixed point is parabolic, the critical
point $-\mathbf{a}_0$ is in a Fatou component attached to $x_0$.
Thus, this Fatou component contains a curve which joins
$-\mathbf{a}_0$ and $x_0$ and avoids $0$. By symmetry ($\sigma$)
this Fatou component contains a curve surrounding $0$\,: this
contradicts the fact that Fatou components are simply connected
for polynomials.
 Therefore $x_0$ is repelling and $\mathbf{a}_0$ is a Misiurewicz
 parameter,  so $x_0=f_{\mathbf{a}_0}(-\mathbf{a}_0)$ (Lemma~\ref{l:misiurpreim1}).

 The fixed ray $R^\infty_{\mathbf{a}_0}(0)\subset \R^+$
also converges to a positive fixed point, say $x_1$. Assume that
$x_1\neq x_0$. Then, from the shape of the graph  of
$f_{\mathbf{a}_0}|_{\R^+}$, it is easy to see that since $x_0$ is
repelling, either there are two other fixed points (one attracting
and one repelling) or there is a parabolic fixed point of
multiplier~$1$. This implies that, including $0$, there are at
least $d+1$ fixed points in $\C$ counted with multiplicity. This
is not possible for a polynomial of degree $d$.

Therefore $x_0=x_1$ and the ray $R^\infty_{\mathbf{a}_0}(0)$ lands
at $f_{\mathbf{a}_0}(-\mathbf{a}_0)$. So  $\RR^{\id}_\infty(0)$
lands at $\mathbf{a}_0$ (by Lemma~\ref{l:misiurewicz}).\cqfd

Proposition~\ref{p:convzero} and Proposition~\ref{p:UcapM} give
the precise dynamical behaviour of $f_\mathbf{a}$ for  parameters
on $\partial \HH_0$.

The parameters on $\R^-$ excluded in all the results are obtained
by symmetry.

\begin{lemma}\label{l:aboutdistpararay}
Two different rays in $\UU$, where $\UU$ is a connected component
of $\HH$, cannot converge to the same parameter.
 \end{lemma}
\proof  The proof is the same in any $\HH_i$ so we do it for
$\UU=\HH_0$. Assume, to get a contradiction, that two rays of
$\HH_0$ converge to the same point $\mathbf{a}_0$. One can suppose
 (up to changing the rays) that they belong to the same $s(\SS)$,
 so that  it is enough to consider the
case $s=\id$. Let  $\RR_0(t)$ and $\RR_0(t')$ be the two rays
under consideration. Let $\gamma$ be the curve
$\RR_0(t)\cup\RR_0(t')\cup\{\mathbf{a}_0\}\cup\{0\}$. There are
infinitely many angles of the form $\frac{p}{q((d-1)^k-1)}$
between $t$ and $t'$, and infinitely many of them give rays
landing to Misiurewicz parameters. Indeed, for such an angle
$\theta$, the  ray $\RR_0(\theta)$ converges to a parameter
$\mathbf{a}$ which is either of Misiurewicz type or such that the
map $f_\mathbf{a}$ has a parabolic point of period $k$ with
multiplier $1$ since it is the landing point of a ray in
$B_\mathbf{a}$ (see Lemma~\ref{l:misiurpreim1}). As there is only
a finite number of parameters $\mathbf{a}$ satisfying the second
alternative (they are solutions of a polynomial equation of degree
at most $d^k$), we know that infinitely many of these landing
parameters   are Misiurewicz points (on $\partial
\HH_0\setminus\{\mathbf{a}_0\}$) lying in the bounded component of
$\C\setminus \gamma$. This contradicts the fact that
 Misiurewicz parameters are landing points of external parameter rays
(Lemma~\ref{l:misiurewicz}) because such external rays will have
to cross $\gamma$ to enter the bounded component of $\C \setminus
\gamma$.
\cqfd

\subsection{Description of the dynamical position of the critical
value.}\label{s:annexe1}

 Note that in degree $d\ge 3$, the position of the critical value does
 not give directly  the position of the critical point like in degree $2$.
We will now give the B\"ottcher coordinate of the critical point
for any parameter $\mathbf{a} \in \SS\cap(\HH_\infty\cup\HH_0)$,
so in this subsection we forget the exponent specifying the sector
for the parameter rays.
  \begin{lemma}\label{l:localisation} For $\mathbf{a}\in \mathcal{R}_\infty(t)$,
the rays
$R^\infty_\mathbf{a}\left(\frac{t}{d}+\frac{\lfloor\frac{d-1}{2}\rfloor}{d}\right)$,
$R_\mathbf{a}^\infty
\left(\frac{t}{d}+\frac{\lfloor\frac{d+1}{2}\rfloor}{d}\right)$
crash on the critical point  $-\mathbf{a}$. If $\mathbf{a}\in
\mathcal{R}_0(t)$, the critical point $-\mathbf{a}$ belongs to
$R_\mathbf{a}^0\left(\frac{t}{d-1}+\frac{\lfloor\frac{d-1}{2}\rfloor}{d-1}\right)$.
  \end{lemma}
\proof  For $\mathbf{a}\in \mathcal{R}_\infty(t)$,
$f_\mathbf{a}(-\mathbf{a})$ belongs to $R^\infty_\mathbf{a}(t)$,
so the two rays crashing on $-\mathbf{a}$ belong to the set of
rays
$\daleth=\left\{R^\infty_\mathbf{a}\left(\frac{t}{d}+\frac{k}{d}\right)\mid
0\le k\le d-1\right\}$. If $\mathbf{a}\in \mathcal{R}_0(t)$, the
critical point  $-\mathbf{a}$ belongs to a unique  ray of the set
$\beth=\left\{R^0_\mathbf{a}\left(\frac{t}{d-1}+\frac{k}{d-1}\right)\mid
0\le k\le d-2 \right\}$.

We first describe  the case of parameters  $\mathbf{a}\in \R^+$\,:
this case is more visual because of the symmetry
$R_\mathbf{a}^p(-\theta)=\sigma({R_\mathbf{a}^p(\theta)})$ (where
$\sigma$ is the complex conjugacy). Then we conclude by moving
$\mathbf{a}$ trough $\SS$.

1) For $0<\mathbf{a}\ll 1$,
 the critical point is
on $R_\mathbf{a}^0(\frac{1}{2})\subset \R^-$, since $\mathbf{a}
\in \HH_0\cap \R^+$. For $d=2l+1$, $\mathbf{a}\in \RR_0(0)$ so
$t=0$\,; we verify then that
$R_\mathbf{a}^0(\frac{1}{2})=R^0_\mathbf{a}(0+\frac{l}{d-1})$. For
$d=2l+2$, $\mathbf{a}\in \RR_0(\frac{1}{2})$, so
$t=\frac{1}{2}$\,; we verify in that case that
$R_\mathbf{a}^0(\frac{1}{2})=R^0_\mathbf{a}\left(\frac{1}{2(d-1)}+\frac{l}{d-1}\right)$.

For $\mathbf{a} \in \SS \cap \HH_0$,
$\phi_\mathbf{a}^0(-\mathbf{a})$ is well-defined and continuous\;
 it belongs to $e^{2i \pi(\frac{t}{d-1}+\frac{k}{d-1})}\R^+$ for $\mathbf{a}\in \RR_0(t)$.
 The integer $k$ is a continuous function of $\mathbf{a}$, so it is  constant and
 therefore equal to
 $\lfloor\frac{d-1}{2}\rfloor$.

 2) We consider now the case  $\mathbf{a}
\gg 1$. The fixed rays $R_\mathbf{a}^0(\frac{k}{d-2})$,
$R_\mathbf{a}^\infty(\frac{k}{d-1})$ with $k\in \N$ are well
defined. Indeed, $-\mathbf{a} \notin K_\mathbf{a}$ so every
rational ray in $B_\mathbf{a}$ converges. Moreover, the only fixed
rays in $\daleth$ are $R_\mathbf{a}^\infty(0)\subset \R^+$ and, if
$d$ is odd, $R_\mathbf{a}^\infty(\frac{1}{2})=]-\infty,p[$ where
$p$ is the unique negative fixed point of $f_\mathbf{a}$\,; note
that  $p<-\mathbf{a}$ since $f_\mathbf{a}(-\mathbf{a})>0$. By
Lemma~\ref{l:aboutdistray0}, distinct internal rays of
$R_\mathbf{a}^0(\pm\frac{k}{d-2})$ ($k\in\N$) converge to distinct
fixed points, named $x_{\pm k}$.
 Those $x_k$ are repelling points ($\mathbf{a}\notin \mathcal C$).
They are the landing points of  external rays (see~\cite{LP}),
which are also fixed rays because of the cyclic order at $x_k$
(see also~\cite{Pe}). Those rays belong to
 $\{R^\infty_\mathbf{a}(\frac{p}{d-1}),\ 0\le p\le d-2\}$.
Because of the symmetry, $R^\infty_\mathbf{a}(\frac{\pm k}{d-1})$
and $R^0_\mathbf{a}(\frac{\pm k}{d-2})$ converge to $x_{\pm k}$
for $0\le k<\lfloor \frac{d}{2}\rfloor$. Thus
 $\gamma_\mathbf{a}= \ol {R^\infty_\mathbf{a}}\left(\pm\frac{ l-1}{d-1}\right)\cup
 \ol{R^0_\mathbf{a}}\left(\pm\frac{ l-1}{d-2}\right)$  is a curve "separating"
  $\C$ into two connected components\,; let  $U_\mathbf{a}$ be the one which  contains $-\mathbf{a}$. The only
rays  of $\daleth$ entering $U_\mathbf{a}$ are
$R_\mathbf{a}^\infty(\frac{d-1}{2d})$, $R_\mathbf{a}^\infty
(\frac{d+1}{2d})$ if $d=2l+1$ and
$R_\mathbf{a}^\infty(\frac{1}{2d}+\frac{l}{d})$,
$R_\mathbf{a}^\infty (\frac{1}{2d}+\frac{l+1}{d})$ if $d=2l+2$\,;
so they crash on $-\mathbf{a}$ for $\mathbf{a} \gg 1$.

By Lemma~\ref{l:douady}, $\gamma_\mathbf{a}$ admits a holomorphic
motion
  parameterized by $\SS \cap \HH_\infty$. Indeed, if $-\mathbf{a}$
  belongs to $\gamma_\mathbf{a}$ the critical value  would describe on the  external
  rays  a sector of opening more than $\frac{d\pi}{d-1}$
 which is  impossible since it is exactly the opening of $\Phi_\infty(\SS)$.
Therefore  $-\mathbf{a}$ stays in $U_\mathbf{a}$ and the only
candidates of $\daleth$ in $U_\mathbf{a}$ are
$R_\mathbf{a}^\infty(\frac{t}{d}+\frac{\lfloor\frac{d-1}{2}\rfloor}{d})$
and $R_\mathbf{a}^\infty
(\frac{t}{d}+\frac{\lfloor\frac{d+1}{2}\rfloor}{d})$, so they
crash on~$-\mathbf{a}$. \cqfd

  \begin{corollary}
For $\mathbf{a} \in \SS$, if the critical point $-\mathbf{a}$
belongs to $R_\mathbf{a}^\infty(\theta)$ then $\theta\in
[\frac{1}{2}-\frac{1}{d},\frac{1}{2}+\frac{1}{d}]$, if
$-\mathbf{a} \in R_\mathbf{a}^0(\theta)$ then
$\theta\in[\frac{1}{2}-\frac{1}{2(d-1)},\frac{1}{2}+\frac{1}{d-1}]$.
  \end{corollary}

   \begin{proposition}\label{p:convinfini}
 For $t\in \Q$, the ray
$\mathcal{R}_\infty(t)$ converges to a parameter
$\mathbf{a}_\infty(t)$. For $l=\lfloor\frac{d-1}{2}\rfloor$\,:
  \begin{enumerate}
\item if  $\frac{t}{d}+\frac{l}{d}$
 {\rm(}resp. $\frac{t}{d}+\frac{l+1}{d}${\rm)} is  periodic by multiplication by $d$,
$\mathbf{a}_\infty(t)$ is a parabolic parameter. The ray
$R^\infty_{\mathbf{a}_\infty(t)}\left(\frac{t}{d}+\frac{l}{d}\right)$
{\rm (}resp.
$R^\infty_{\mathbf{a}_\infty(t)}\left(\frac{t}{d}+\frac{l+1}{d}\right)${\rm)}
lands at a parabolic point $p$, root of the Fatou component $P_t$
containing the critical point $-\mathbf{a}_\infty(t)$\ and  the
ray
$R^\infty_{\mathbf{a}_\infty(t)}\left(\frac{t}{d}+\frac{l+1}{d}\right)$
(resp.
$R^\infty_{\mathbf{a}_\infty(t)}\left(\frac{t}{d}+\frac{l}{d}\right)$)
lands at the preimage of $f_{\mathbf{a}_\infty(t)}(p)$ on
$\partial P_t$; \item otherwise $\mathbf{a}_\infty(t)$ is a
Misiurewicz parameter,
$R^\infty_{\mathbf{a}_\infty(t)}\left(\frac{t}{d}+\frac{l}{d}\right)$
and
$R^\infty_{\mathbf{a}_\infty(t)}\left(\frac{t}{d}+\frac{l+1}{d}\right)$
land at $-\mathbf{a}_\infty(t)$.
  \end{enumerate}
  \end{proposition}

\begin{figure}[!h]\vskip 0cm
\centerline{\input{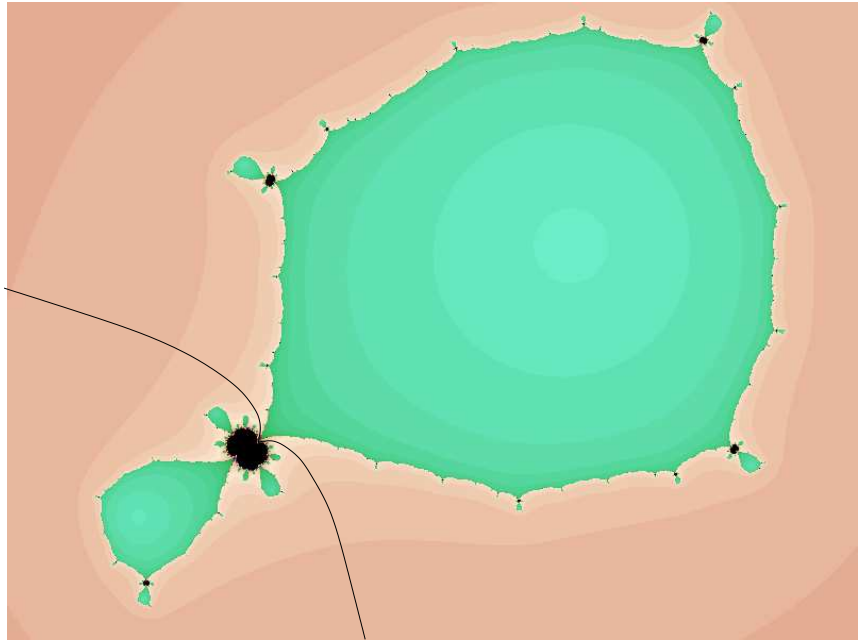}} \caption{Julia set with a
parabolic point in $\CC_4$.}\label{f:para}
\end{figure}

\begin{figure}[!h]\vskip 0cm
\centerline{\input{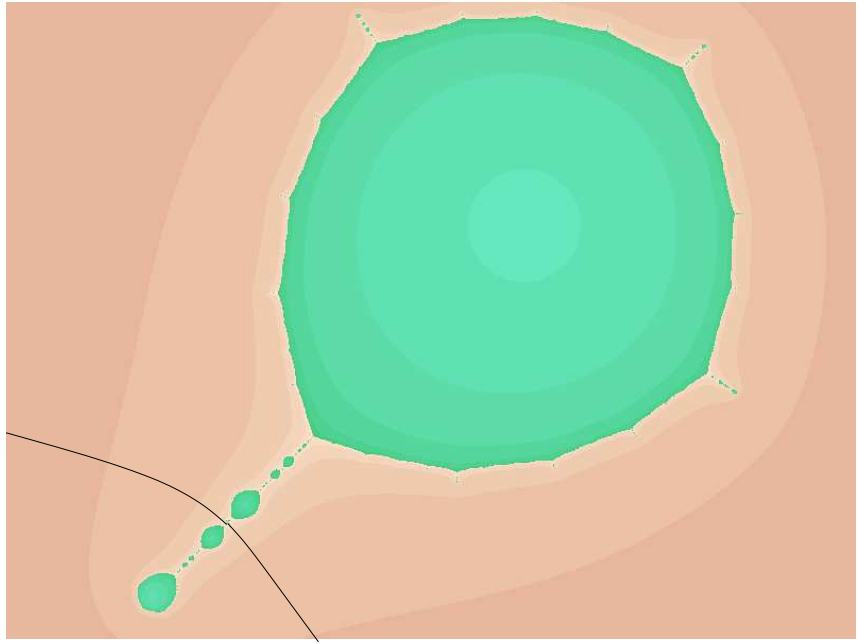}} \caption{Two rays converging
to the critical point.}\label{f:renorm4}
  \end{figure}

\proof We will use the notations\,:
$t_0=\frac{t}{d}+\frac{\lfloor\frac{d-1}{2}\rfloor}{d}$,
$t_1=\frac{t}{d}+\frac{\lfloor\frac{d+1}{2}\rfloor}{d}$ and
$x_0,x_1$ for the landing points of
$R_{\mathbf{a}_\infty}^\infty(t_0)$,
 $R_{\mathbf{a}_\infty}^\infty(t_1)$ respectively.
We distinguish two cases according to wether $t_0$, $t_1$ are
periodic or not.

 \noindent 1) {\it  $t_0$ is
periodic }: The point $x_0$, which is eventually critical or
parabolic, is now periodic and in the Julia set. So it is
necessarily parabolic. If $t_0$ is fixed the result follows.
Assume now that $t_0$ is not fixed. Suppose, to get a
contradiction, that $x_0$ is not the root of $P_t$. So some ray
with angle $t_2=d^i t_0\neq t_0$ converges  to the root of the
Fatou component $P_t$ which contains the critical point, that is
to the image $f_\mathbf{a}^i(x_0)$ belonging to $\partial P_t$.
Since the angles $\{\frac{k}{d^j(d-1)}, k\in\N,\ j \in\N\}$ are
dense in $\R/\Z$ (the distance between two consecutive terms tends
to 0), there exist three angles of this form called $\theta_j\neq
0$ separating $t_0$, $t_1$ and $t_2$. Hence the external rays of
angles $\theta_j$ with the internal rays which have the same end
points form a curve $\delta$ that separates $\C$ into three
connected components\,: one contains
$R^\infty_{\mathbf{a}_\infty(t)}(t_2)$ and the critical point, the
other ones contain $R^\infty_{\mathbf{a}_\infty(t)}(t_0)$ and
$R^\infty_{\mathbf{a}_\infty(t)}(t_1)$.
 By Lemma~\ref{l:douady}, $\delta$ varies continuously for $\mathbf{a} \in
\mathcal{R}_\infty(t)$ since the critical point cannot be on
$\gamma$ (Lemma~\ref{l:localisation}). Here we get a contradiction
since the critical point varies continuously, cannot cross
$\delta$ and has to break the rays
$R^\infty_{\mathbf{a}_\infty(t)}(t_0)$ and
$R^\infty_{\mathbf{a}_\infty(t)}(t_1)$
(Lemma~\ref{l:localisation}). Hence $x_0$ is the root of $P_t$.
Let $x'$ be the other preimage on $\partial P_t$ of
$f_{\mathbf{a}_\infty(t)}(x_0)$. We want to prove now that
$R^\infty_{\mathbf{a}_\infty(t)}(t_1)$ converges to $x'$. The
curve $\gamma$ formed by all the fixed  rays (defined in
Lemma~\ref{l:localisation}) varies continuously with $\a$ on
$\RR_\infty(t)$ until $\mathbf{a}_\infty(t)$ since at this
parameter the critical point is not on $\gamma$. So the rays
$R_{\mathbf{a}_\infty(t)}^\infty(t_0)$ and
$R_{\mathbf{a}_\infty(t)}^\infty(t_1)$ are the only preimages of
$R_{\mathbf{a}_\infty(t)}^\infty(dt_0)$ in the same connected
component of $\C \setminus \gamma$ as $P_t$. So
$R^\infty_{\mathbf{a}_\infty(t)}(t_1)$ converges to $x'$.

If $t_1$ is periodic instead of $t_0$ the proof is the same.

 \noindent 2) {\it Neither $t_0$ nor $t_1$ is periodic}: Assume by contradiction
 that  $x_0$ is
eventually parabolic, by the same argument using $\gamma$ as
before, $R_{\mathbf{a}_\infty(t)}^\infty (t_0)$ lands at the root
 $p$ of $P_t$. But $t_0$ is not periodic, so there
is another ray $R^\infty_{\mathbf{a}_\infty(t)}(d^it_0)$ landing
at $p$. Let $k$ be the first integer such that $d^{i+k}t_0\equiv
d^k t_0\mod 1$\,; then $f_\mathbf{a}^{k-1}(p)$ must be critical
since two different rays land at this point\,:
$R^\infty_{\mathbf{a}_\infty(t)}( d^{i+k-1}t_0)$ and
$R_{\mathbf{a}_\infty(t)}^\infty(d^{k-1} t_0)$ which have the same
image $R_{\mathbf{a}_\infty(t)}^\infty(d^{k} t_0)$. Finally, by
Lemma~\ref{l:douady}, $x_0$ is eventually critical.\cqfd

  \begin{proposition}\label{p:convzero}
For $t\in \Q$, let $\mathbf{a}_0(t)$ be the landing point of
$\mathcal{R}_0(t)$. We have\,:
  \begin{enumerate}
\item if $ \
\frac{t}{d-1}+\frac{\lfloor\frac{d-1}{2}\rfloor}{d-1}$ is periodic
by multiplication by $d-1$, then $\
R^0_{\mathbf{a}_0(t)}\left(\frac{t}{d-1}+\frac{\lfloor\frac{d-1}{2}\rfloor}{d-1}\right)$
lands at a parabolic point which is the  root of the Fatou
component $P_t$ which contains $-\mathbf{a}_0(t)$; \item
otherwise, $\
R^0_{\mathbf{a}_0(t)}\left(\frac{t}{d-1}+\frac{\lfloor\frac{d-1}{2}\rfloor}{d-1}\right)$
lands at the critical point  $-\mathbf{a}_0(t)$.
  \end{enumerate}
  \end{proposition}
\proof The proof goes exactly as the one of
Proposition~\ref{p:convinfini}.\cqfd

  \begin{lemma}\label{l:misiurewicz2} Let $\mathbf{a}$  be a Misiurewicz point on
  $\partial \HH_j$ with
$j\in\{0,\infty\}$, $$R_\mathbf{a}^{r(\mathbf{a})}(t)\quad
\hbox{lands at } \quad f_\mathbf{a}(- \mathbf{a}) \iff
\RR_j(t)\quad \hbox{lands at } \quad \mathbf{a}.$$
  \end{lemma}
\proof The proof of the implication $\Rightarrow$ is exactly the
same as in degree $2$ (see \cite{DH1}). The proof of $\Leftarrow$
is just the case 2 of Proposition~\ref{p:convinfini} and
Proposition~\ref{p:convzero} since $-\mathbf{a}$ cannot be at the
same time eventually periodic and attracted by a parabolic point.

We can now give another proof of the following corollary\,:
  \begin{crrend}\hskip -0.75em{\bf \ref{c:R+}}
  If $d$ is odd, $\R^+=\RR^{\id}_\infty(0)\cup \RR^{\id}_0(0)\cup
\{*\}$ where $*$ is a Misiurewicz point.
  \end{crrend}
\proof If $d=2l+1$, $\R^+ \cap \HH_0=\RR_0(0)$  so $t=0$ and
$\frac{t}{d-1}+\frac{\lfloor\frac{d-1}{2}\rfloor}{d-1}=\frac{l}{d-1}=\frac{1}{2}$\,;
this angle is not periodic by multiplication by $d-1$, neither the
angle $\frac{l+1}{d-1}$ since their images are $l$ and $l+1$. Thus
$\RR_0(0)$ converges to a parameter $\mathbf{a}_0$ such that
$R_{\mathbf{a}_0}^0(\frac{1}{2})$ lands $-\mathbf{a}_0$. So
$\mathbf{a}_0$ is a Misiurewicz point. The unique positive fixed
point is $f_{\mathbf{a}_0}(-\mathbf{a}_0)$ so the fixed ray
$R^\infty_{\mathbf{a}_0}(0)\subset \R^+$ converges to
$f_{\mathbf{a}_0}(-\mathbf{a}_0)$ and by
Lemma~\ref{l:misiurewicz2}, $\RR_\infty(0)$ lands at
$\mathbf{a}_0$.\cqfd

\section{Graphs, puzzles, para-graphs and para-puzzles}\label{s:graphes}
This section is devoted to the construction of the puzzles and the
para-puzzles. We recall in section~\ref{s:dyngraph} the graphs
used in \cite{Ro1} and we construct in section~\ref{s:paramgraph}
the analoguous graphs in the parameter plane, then we establish
the relations between graphs and  para-graphs
(section~\ref{s:graphmotion} and~\ref{s:realtiong-pg}) and show
how to use them for the question of local connectivity
(section~\ref{s:graploccon}).
 

\subsection{Dynamical puzzles and graphs }\label{s:dyngraph}


We define the puzzles and the graphs (as in~\cite{Ro1}) when the
Julia set is connected, {\it i.e.} for parameters in $\CC\setminus
\R^-$, and obtain the puzzles  in the other cases through  a
holomorphic motion of these graphs.

 Fix  $\mathbf{a}\in \CC \setminus \R^-$. For  large  $l$ and  $\theta_\pm
=\pm\frac{ 1}{(d-1)^l-1}$, the ray $R_\mathbf{a}^0(\theta_\pm)$ is
well defined and  converges to a point~$x_\pm$ which is repelling
($f_\mathbf{a}$ possesses at most one parabolic orbit). Let
$\eta_\pm$ be the angle of some external ray
$R^\infty_\mathbf{a}(\eta_\pm)$ landing at~$x_\pm$ (given by
Proposition~\ref{p:yocc}).
Since the internal ray $R_\mathbf{a}^0(\theta_\pm)$ is fixed by
$f_\mathbf{a}^l$, the external ray $R^\infty_\mathbf{a}(\eta_\pm)$
will also be fixed by $f_\mathbf{a}^l$ so $\eta_\pm$ is periodic
(see Remark~\ref{r:aboutdistrayext}). Using  these rays we
construct the graphs in $\ol X^\mathbf{a}$, where $X^\mathbf{a} =
\{z \in\C \mid G^0_\mathbf{a}(z)<1 \quad\hbox{and}\quad
G^\infty_\mathbf{a}(z)<1\}$, as follows.

  \begin{definition} \label{d:puzzle}
Let $\mathbf{a} \in  \CC \setminus  \R^-$, $\theta_\pm =
\pm\frac{1}{(d-1)^l-1}$ with $l$ large (as before). A \emph{
puzzle} for~$f_\mathbf{a}$ is defined by the following graph\,:
 $$
I^\mathbf{a}_0(\theta_\pm) = \partial X^\mathbf{a} \cup \Bigl(
X^\mathbf{a} \cap \Bigl( \Cup_{i\ge0} \left(\ol
{R^0_\mathbf{a}((d-1)^i\theta_\pm)}\cup \ol{ R^\infty_\mathbf{a}
\bigl( d^i\eta_\pm \bigr)}\right) \Bigr)\Bigr).
 $$
  The \emph{puzzle pieces of depth~$n\ge0$} are the connected
components of $$f_\mathbf{a}^{-n}(X^\mathbf{a}) \setminus
I^\mathbf{a}_n = f_\mathbf{a}^{-n} (X^\mathbf{a} \setminus
I_0^\mathbf{a}),\hbox{ where } I^\mathbf{a}_n = f_\mathbf{a}^{-n}
(I^\mathbf{a}_0) \quad\text{for all }n\ge1.
$$ The puzzle piece containing a given point~$z$
is denoted by $P^\mathbf{a}_n(z)$. The puzzle pieces containing
the critical value~$f_\mathbf{a}(-\mathbf{a})$ are denoted simply
by $P_0^\mathbf{a}, \dots, P_n^\mathbf{a}, \dotsc$ if there is no
ambiguity.

The \emph{puzzle } is the union of all the puzzle pieces.
  \end{definition}
\begin{remark}\label{r:aboutdistrayext}
The ray $R^\infty_\mathbf{a}(\eta_\pm)$ is the only external ray
of the cycle $R^\infty_\mathbf{a}(d^j\eta_\pm)$, $j\ge 0$,  to
converge to $x_\pm$.
\end{remark}
\proof  Assume (to get a contradiction) that
$R^\infty_\mathbf{a}(d^j\eta_\pm)$ with $d^j\eta_\pm\neq \eta_\pm
\hbox{ mod } 1$ converges to $x_\pm$. Since there is a finite
number of rays in the cycle converging to $x_{\pm}$ we can assume
(up to changing the notations) that the angles  are all in the
interval $(\eta_{\pm}, d^j\eta_\pm)$.
  Since  the map $f^{l-j}_\mathbf{a}$ is
conformal at this point, it preserves the ``cyclic order'' of the
rays at $x_\pm$. But it maps $R^\infty_\mathbf{a}(d^j\eta_\pm)$ to
$R^\infty_\mathbf{a}(\eta_\pm)$  and
$R^\infty_\mathbf{a}(\eta_\pm)$ to
$R^\infty_\mathbf{a}(d^{l-j}\eta_\pm)$. These rays land at $x_\pm$
but (because of the cyclic order) $d^{l-j}\eta_\pm$ will not be
 in the interval  $(\eta_{\pm}, d^j\eta_\pm)$.
So the two rays cannot be in the same cycle.
\cqfd

If we fix some $\theta$ as in definition~\ref{d:puzzle} but vary
the parameter $\mathbf{a}$ in $\C\setminus\R^-$, we will see that
for the graph $I^\mathbf{a}_0(\theta)$  the landing points of the
rays, $x_\pm$, can become parabolic, the rays
$R^\mathbf{a}_\infty(\eta_\pm)$ and $R_\mathbf{a}^0(\theta_\pm)$
can land at different points, the rays can crash on critical
points and no more be well-defined etc... For these reasons we
should restrict the domain (in the parameter space) on which we
consider the graph at each depth. The para-puzzle pieces defined
in section~\ref{s:paramgraph} correspond exactly to the region
were the dynamical pieces are defined by the same rays.


\subsection{Para-graphs and para-puzzles}\label{s:paramgraph}

 The para-graph are just the copy of dynamical graphs in the
parameter plane so depend from a preferred parameter.
 Let $\mathbf{a}_0\in  \CC  \setminus  \R^-$ and
 $I^{\mathbf{a}_0}_0(\theta)$ be  the graph associated to this parameter by
Definition~\ref{d:puzzle}.

 \begin{definition}
For $n\in \N$, let  $\kappa_n$  be the set of all the pairs
$(\mathcal U,v)$ where\,:
\begin{itemize} \item
$\mathcal U$ is a connected component of $\HH_i$ with $i\in
\{\infty,0,\cdots,n\}$\,;
 \item $(d-1)^{n-i}v=1$ if
 $0\le i \le n$ and  $d^nv=1$ if $i=\infty$.
\end{itemize}
Let $\XX_n$  be the connected component containing $\mathbf{a}_0$
of $ \C  \setminus ( \Cup_{(\mathcal U,v)\in \kappa_n}\ol{\mathcal
E_{\mathcal U}(v)})$.
 \end{definition}

  \begin{definition}
For $n\in \N$, let  $\Theta_n$  be the set of all the pairs
$(\mathcal U,t)$ where\,:
\begin{itemize} \item
$\mathcal U$ is a connected component of $\HH_i$ with $i\in
\{\infty,0,\cdots,n\}$\,;
 \item ${(d-1)^{n-i}}t \in
\{\theta,\ldots,(d-1)^{l-1}\theta\}$ if $0\le i\le n$\,;
 \item $d^nt\in \{d^{j}\eta, \ j\ge 0\}$ if $i=\infty$.
\end{itemize}
 \end{definition}

\begin{figure}[!h]\vskip 0cm
\centerline{\input{graph.pstex_t}} \caption{Schematic
representation of a para-graph $\II_0$ in $\CC_4$.}\label{f:graph}
  \end{figure}
\begin{figure}[!h]\vskip 0cm
\centerline{\input{graph1.pstex_t}} \caption{Schematic
representation of  para-graphs $\II_0$ and $\II_1$  in
$\CC_4$.}\label{f:graph1}
  \end{figure}

  \begin{definition}\label{d:ppuzzle}
The \emph{para-puzzle} is the union of the \emph{para-puzzle
pieces}.  The \emph{para-puzzle pieces of depth~$n$} are the
connected components of $\XX_n \setminus \II_n$ where
$$ \II_n (\theta)=\partial\XX_n \cup
    \Cup_{ (\mathcal U,t)\in \Theta_n}
   \left( \ol {U\RR_{\mathcal U}(t)} \cap \XX_n \right) .$$
   The para-puzzle piece  containing a given parameter~$\mathbf{a}$ will be
denoted by $\PP_n(\mathbf{a})$. For the given parameter
$\mathbf{a}_0$ we will simply write $\PP_n$ for
$\PP_n(\mathbf{a}_0)$.
  \end{definition}
  The points of $\R^-$ can be added or not to the graph. We only
  care of parameters in $\SS$ and for these parameters, all
the para-puzzle pieces are compactly contained in $\C\setminus
\R^-$ by the following Remark.
\begin{remark}\label{r:inclusdsS}
 For parameters $\mathbf{a}_0\in \SS$, the para-puzzle piece $\PP_0(\mathbf{a}_0)$
 is compactly contained in $\C\setminus \R^-$.
 Moreover there is
 only one connected component of $U\RR_p(t)\setminus\{p\}$ in $\PP_0(\mathbf{a}_0)$.
 Thus for the simplicity of the exposition we will forget the exponent in $\RR^s_p(t)$
 and call  $\RR_p(t)$ the parameter ray that belongs to $\PP_0(\mathbf{a}_0)$.
 \end{remark}
 \proof By definition
 of the para-graphs, any piece of  depth $0$ is bounded only by rays in
 $\HH_0$ and in $\HH_\infty$ as well as equipotentials. If
  $\PP_0(\mathbf{a}_0)$ intersects $\R^+$ then it is compactly contained
  in the interior of $\SS\cup \sigma \SS\cup \tau^{-1} \SS$.
  Indeed, the sector $\sigma \SS\cup \tau^{-1}
  \SS=\rho^{-1}(\SS\cup\sigma\SS)$ with $\rho=e^\frac{i\pi}{d-1}$
  contains all the rays of $\HH_0$ in a sector of angles of width
  greater than $\displaystyle 1+\frac{1}{d-2}$ by proposition~\ref{p:H0}
   and  using  the  coordinate $\Phi_0$.
Therefore,  $\PP_0(\mathbf{a}_0)$ is contained in the interior of
$\SS\cup \sigma \SS\cup \tau^{-1} \SS$ since   by definition
$\PP_0(\mathbf{a}_0)\cap \HH_0$ corresponds to angles in a sector
of width less than $1$. Using the same argument, if
$\PP_0(\mathbf{a}_0)$ intersects $\rho\R^+$ then it is compactly
contained
  in the interior of $\SS\cup \tau\sigma \SS\cup \tau \SS$
  which is $\SS\cup \tau(\SS\cup \sigma \SS$).
 Hence, $\PP_0(\mathbf{a}_0)$ is compactly contained in
$\C\setminus \R^-$  since the interior of $\SS\cup \rho \SS \cup
\sigma \SS$ is included in $\C\setminus \R^-$.

Since   by definition $\PP_0(\mathbf{a}_0)\cap \HH_0$ corresponds
to angles in a sector of width less than $1$,
$U\RR_0(t)\setminus\{0\}$ has only one connected component in
$\PP_0(\mathbf{a}_0)$. The same holds for
$U\RR_\infty(t)\setminus\{\infty\}$. \cqfd

 \subsection{Holomorphic motion of the dynamical graphs}~\label{s:graphmotion}

  \begin{definition} Let $\Lambda$ be a $\C$-analytic variety. Let
  $\lambda_0\in \Lambda$.
A holomorphic motion of a subset $\Gamma \subset \ol \C$
parameterized by $\Lambda$  is a map $\Psi \from \Lambda\times
\Gamma \to \ol \C$ such that $\Psi(.,z)$ is holomorphic on
$\Lambda$, $\Psi^\lambda=\Psi(\lambda,.)$ is injective on $\Gamma$
and $\Psi^{\lambda_0}=\id$.
   \end{definition}

For the given parameter $\mathbf{a}_0\in \CC\setminus \R^-$, we
define now the set of parameters for which the graph
$I_0^{\mathbf{a}_0}(\theta)$ admits a holomorphic motion. Let
$\eta$ be the angle of the external ray converging to the same
point as $R_{\mathbf{a}_0}^0(\theta)$ in
$I_0^{\mathbf{a}_0}(\theta)$.
 
  \begin{lemma}\label{l:holom}
  Let $\Omega_{\eta}$ be the set of parameters $\mathbf{a}\in \C\setminus \R^-$ such that
  for all  $i\ge 0$  the ray $R_\mathbf{a}^\infty(d^i\eta)$ is well-defined
  and converges to a  repelling periodic point.
\begin{enumerate}
\item $\Omega_{\eta}$ is a non empty  open set\,; \item
  The set  $\Gamma_{\mathbf{a}_0}^\infty(\eta)=\Cup_{i\geq 0}
  \ol{R_{\mathbf{a}_0}^\infty(d^i\eta)}$
  admits a holomorphic motion
parameterized by $\Omega_\eta$\,; \item The boundary $
\partial \Omega_{\eta}$ is a subset of $\R^-\cup \Cup_{i\ge
0}\ol{U\RR_\infty(d^i\eta)}$ without  isolated points.
\end{enumerate}
  \end{lemma}
\proof 1. and 2.\,:  Since $\mathbf{a}_0\in \Omega_{\eta}$, it is
clearly not empty. To prove that $\Omega_{\eta}$ is open, take
$\mathbf{a}_1\in\Omega_{\eta}$. For each $i\ge 0$, the landing
point of $R_{\mathbf{a}_1}^\infty(d^i\eta)$  is periodic and
repelling, thus not eventually  critical. Therefore the compact
set $\ol{ R_{\mathbf{a}_1}^\infty(d^i\eta)}$ admits a holomorphic
motion in some neighbourhood of $\mathbf{a}_1$
(lemma~\ref{l:douady}). The set
$\Gamma_{\mathbf{a}_1}^\infty(\eta)$ is  a finite union of such
compact sets, so it admits a holomorphic motion in a neighbourhood
of $\mathbf{a}_1$. Hence $\Omega_{\eta}$ is open.

3. a) Assume by  contradiction that there is a point
$\mathbf{a}_1$ of $\C\setminus\R^-$  isolated in $\partial
\Omega_{\eta}=\ol{ \Omega_{\eta}}\setminus \Omega_{\eta}$. Then
there exists an open neighbourhood $\mathcal O$ of $\mathbf{a}_1$
such that $\mathcal O \setminus \{\mathbf{a}_1\}\subset
\Omega_{\eta}$. Since the parameter $\mathbf{a}_1$ is not in $
\Omega_{\eta}$, one of the rays either is not well defined or
converges to a parabolic point. If the ray
$R_{\mathbf{a}_1}^\infty(d^i\eta)$ crashes on $-\mathbf{a}_1$,
then $f_{\mathbf{a}_1}(-\mathbf{a}_1)$ belongs to
$R_{\mathbf{a}_1}^\infty(d^{i+1}\eta)$ and $\mathbf{a}_1 \in
\RR^s_\infty(d^{i+1}\eta)$ for some $s$ composition of $\tau$ and
$\sigma$.
 Then for parameters $\mathbf{a}\in \RR^s_\infty(d^{i+1}\eta)\cap
\mathcal O$ near $\mathbf{a}_1$, the ray
$R^{\infty}_{\mathbf{a}}(d^i\eta)$ also crashes. This contradicts
the fact that $\mathcal O \setminus \{\mathbf{a}_1\}\subset
\Omega_{\eta}$. Consider now the case where all the rays
$R_{\mathbf{a}_1}^\infty(d^i\eta)$ are well-defined but  converge
to a parabolic periodic cycle ($d^i\eta$ is periodic). For
$\mathbf{a}\in \mathcal O$, the landing point $x_i(\mathbf{a})$ of
the ray $R_\mathbf{a}^\infty(d^i\eta)$  defines a holomorphic map
(Lemma~\ref{l:douady}). It is repelling for $\mathbf{a} \in
\mathcal O \setminus \{\mathbf{a}_1\}$ and parabolic at
$\mathbf{a}_1$. So $(f_\mathbf{a}^k)'(x(\mathbf{a}))$ can be
extended to a holomorphic map from $\mathcal O$ to $\C \setminus
\D$ (where $k$ denotes the period of the cycle). Its modulus
reaches its minimum at $\mathbf{a}_1$\,; this contradicts the
maximum principle for the map $\mathbf{a}\mapsto
1/(f_\mathbf{a}^k)'(x(\mathbf{a}))$.

 b) For $\mathbf{a} \in \partial \Omega_{\eta}\setminus\R^-$, either one
ray $R_\mathbf{a}^\infty(d^i\eta)$ crashes on the critical point
$-\mathbf{a}$ and so $\mathbf{a}\in U\RR_\infty(d^{i+1}\eta)$, or
the rays $R_\mathbf{a}^\infty(d^i\eta)$ converge to a parabolic
periodic cycle and so $\mathbf{a} \in \mathcal{PA}=\{\mathbf{a}
\mid \exists\ x\hbox{ such that } f_\mathbf{a}^k(x)=x \hbox{ and }
(f_\mathbf{a}^k)'(x)=1\}$ (where $k$ is the period of $\eta$). The
set $\mathcal{PA}$ is finite  since every $\mathbf{a}\in
\mathcal{PA}$
 is a root of the discriminant of the two polynomials $f_\mathbf{a}^{k}(z)-z$ and
$(f_\mathbf{a}^{k})'(z)-1$. Since there is no isolated points in
$\partial \Omega_{\eta}$  those parameters of
$\mathcal{PA}\cap\partial \Omega_{\eta}$ are in the closure $
U\RR_\infty(d^i\eta)$ for $i\geq 0$. \cqfd

\noindent The same result holds for internal rays\,:
\begin{lemma}\label{l:holom2}
  Let $\Omega'_{\theta}$ be the set of parameters
   $\mathbf{a}\in \C\setminus \R^-$ such that
  for all  $i\ge 0$  the ray $R_\mathbf{a}^0((d-1)^i\theta)$
  is well-defined and
converges to a  repelling periodic point.
\begin{enumerate}
\item $\Omega'_{\theta}$ is a non empty  open set\,; \item
  The set  $\Gamma^0_{\mathbf{a}_0}(\theta)=\Cup_{i\geq 0} \ol{R_{\mathbf{a}_0}^0((d-1)^i\theta)}$
  admits a holomorphic motion
parameterized by $\Omega'_\theta$ \,; \item The boundary $
\partial \Omega'_{\theta}$ is a subset of $\R^-\cup \Cup_{i\ge
0}\ol{U\RR_0((d-1)^i\theta)}$. It   has no isolated points.
\end{enumerate}
  \end{lemma}
 \begin{corollary}\label{c:stable} In the connected component containing $\mathbf{a}_0$
 of $\Omega_{\eta}\cap\Omega'_{\theta}$,
 the rays $R_\mathbf{a}^0(\theta)$ and $R_{\mathbf{a}}^\infty(\eta)$  land at
a common point.
  \end{corollary}
\proof The landing points $x_0(\mathbf{a})$ of
$R_{\mathbf{a}}^0(\theta)$, $x_\infty(\mathbf{a})$ of
$R_{\mathbf{a}}^\infty(\eta)$, are both repelling periodic points.
The period is determined by the angles $\theta$, and $\eta$, and
is at most say  $k$. At the parameter $\mathbf{a}_0$ the points
coincide by definition of the graph\,:
$x=x_0(\mathbf{a}_0)=x_\infty(\mathbf{a}_0)$. Since $x$ is
repelling, by Rouch\'e's Theorem, on some neighbourhood $U$ of $x$
there is exactly one point of period less than $k$, for
$\mathbf{a}$ in a neighbourhood $\UU\subset
\Omega_{\eta}\cap\Omega'_{\theta}$ of $\mathbf{a}_0$. Moreover,
the points $x_0(\mathbf{a})$ and $x_\infty(\mathbf{a})$ vary
continuously for $\mathbf{a}\in
 \Omega_{\eta}\cap\Omega'_{\theta}$ (Lemma~\ref{l:holom} and
Lemma~\ref{l:holom2}). Therefore, they coincide on $\UU$ and
finally  on the connected component containing $\mathbf{a}_0$ of
$\Omega_{\eta}\cap\Omega'_{\theta}$ (since they are holomorphic
maps).\cqfd

  \begin{corollary}\label{c:holom} The para-puzzle piece $\PP_0$
  is contained in the connected component of
  $\Omega_{\eta } \cap \Omega'_{\theta}$ containing $\mathbf{a}_0$.
  Therefore the  graph $I_0^{\mathbf{a}_0}(\theta)$  admits a holomorphic
  motion defined on $\mathcal P_0(\mathbf{a}_0)$ so that $I_0^\mathbf{a}(\theta)$
  is well defined {\rm(}{\it i.e.} the rays of
the graph  are well defined{\rm)}.
  \end{corollary}
\proof The boundary $\partial\Omega'_\theta$ is included in
$\Omega_\eta\cup\R^-$, except for the landing points of the rays,
since $\partial\Omega'_\theta\cap\partial\Omega_\eta\subset
\R^-\cup(\ol{U\RR_0((d-1)^i\theta)}\cap\ol{U\RR_\infty(d^i\eta)})$.
The same holds for $\partial\Omega_\eta$, so the boundary of
$\Omega_\eta\cap \Omega'_\theta$ is simply the union $\partial
\Omega_\eta\cup
\partial\Omega'_\theta$. Thus  $(\partial \Omega_\eta\cup
\partial\Omega'_\theta)\cap \XX_0$ is included in $\II_0(
\theta)\cup\R^-$. Therefore  $\PP_0$ is contained in the connected
component of $\Omega_\eta\cap \Omega'_\theta$ containing
$\mathbf{a}_0$ (since $\PP_0=\PP_0(\mathbf{a}_0)$ is a connected
component of $\C \setminus \II_0(\theta)$ in $\C\setminus\R^-$).

Hence, the graph $I_0^\mathbf{a}(\theta)=\partial X_\mathbf{a}\cup
((\Gamma^\infty_\mathbf{a}(\eta)\cup
\Gamma^0_\mathbf{a}(\theta))\cap X_\mathbf{a})$ is well defined
for $\mathbf{a}$ in $\PP_0(\mathbf{a}_0)$, since
$\PP_0(\mathbf{a}_0)$ is included in $ \XX_0 $. The holomorphic
motion of the graph $I_0^{\mathbf{a}_0}$ follows from
Lemma~\ref{l:holom}, Lemma~\ref{l:holom2} and the fact that the
map $(\mathbf{a},z)\to
\phi^p_\mathbf{a}\circ(\phi^p_{\mathbf{a}_0})^{-1}(z)$ defines a
holomorphic motion of the equipotentials $E^p_{\mathbf{a}_0}(1)$
for $\mathbf{a}$ in $\XX_0$ and $p\in \{0,\infty\}$. \cqfd

>From now on through the rest of the paper, we restrict ourself to
parameters $\mathbf{a}_0 \in \SS$ and study para-puzzles inside
the open region $\PP_0(\mathbf{a}_0)$. Hence by
Remark~\ref{r:inclusdsS} we don't need assumptions on the sector
containing the parameters considered.

  \begin{corollary}\label{c:misiurgraph}
For $n\ge 1$, the points which are  in $\PP_0$ of $\II_n(
\theta)\cap\CC$  are of Misiurewicz type.
  \end{corollary}
\proof A parameter $\mathbf{a} \in  \II_n( \theta)\cap \partial
\CC$ is necessarily the landing point of a ray  $\RR_\infty(t)$
with $d^nt\in \{d^j\eta, j\ge 0\}$ (by definition of $\II_n(
\theta)$). Thus the ray $R_\mathbf{a}^\infty(t)$ belongs to
$I_n^{\mathbf{a}}$ since its image by $f_\mathbf{a}^n$,
$R_\mathbf{a}^\infty(d^n t)$, belongs to $I_0^\mathbf{a}$ (by
definition of $t$ and of $I_0^\mathbf{a}$). Since they are in
$\PP_0$, the landing points of rays in $I_0^\mathbf{a}$ are
repelling periodic points (Lemma~\ref{l:holom} and
Corollary~\ref{c:stable}). Therefore $\mathbf{a}$ is a Misiurewicz
point since we are in the second alternative of
Lemma~\ref{l:misiurpreim1}\,: $R_\mathbf{a}^\infty(d^nt)$ lands at
a repelling periodic point.\cqfd

  \begin{lemma}\label{l:graph}
For parameters  $\mathbf{a} \in \mathcal P_0(\mathbf{a}_0)$, the
following equivalence holds\,: $$\mathbf{a}\in \II_n \iff
f_\mathbf{a}(-\mathbf{a})\in I^\mathbf{a}_{n}.$$
  \end{lemma}

  \begin{proof}
By construction of $\II_n$, the rays and equipotentials involved
in $\II_n$ and $I^\mathbf{a}_n$  correspond to each other via the
change of coordinates (Remark~\ref{r:rayons0} and
Corollary~\ref{c:corresp}). From Corollary~\ref{c:misiurgraph} and
its proof, the  points in $\II_n\cap \partial \CC$ are Misiurewicz
points and $f_\mathbf{a}(-\mathbf{a})$ is the landing point of the
corresponding ray in $I_n^\mathbf{a}$. Conversely, if
$f_\mathbf{a}(-\mathbf{a})\in I_n^\mathbf{a}$ is in the Julia set,
it is the landing point of some external ray
$R^\infty_\mathbf{a}(t)$ of $I_n^\mathbf{a}$, so $d^nt\in \{
d^j\eta\quad j\ge 0\}$. Since $f_\mathbf{a}(-\mathbf{a})$ is
eventually periodic, $\mathbf{a}$ is a Misiurewicz point and by
Lemma~\ref{l:misiurewicz} the external ray $\RR_\infty(t)$ lands
at $\mathbf{a}$. Hence, the parameter $\mathbf{a}$ belongs to
$\II_n$ (by definition of this para-graph).
  \end{proof}

\begin{corollary}\label{c:crit} For $\mathbf{a}\in \PP_n(\mathbf{a}_0)$, the $n$-th para-puzzle piece,
the critical point $-\mathbf{a}$ is not  on any of  the graphs
$I^\mathbf{a}_0, \dots, I^\mathbf{a}_{n+1}$.
\end{corollary}
\begin{corollary}\label{c:1conx}
The para-puzzle pieces are simply connected.
\end{corollary}
\proof It is equivalent  to prove that the graph $\II_n$ is
connected. Any part of an equipotential involved in $\II_n\cap
\HH_i$ ($i\in \N$) is connected to $\partial \CC$ by a ray in
$\II_n$. By Corollary~\ref{c:misiurgraph}, this ray converges to a
Misiurewicz parameter, say $\mathbf{a}_1$. At this parameter,  in
the dynamical graph $I_n^{\mathbf{a}_1}$, some external ray
$R_{\mathbf{a}_1}^\infty(t')$  converges to
$f_{\mathbf{a}_1}(-\mathbf{a}_1)$. Then the external parameter ray
$\RR_\infty(t')$ of $\II_n$ (by Lemma~\ref{l:graph}) converges to
$\mathbf{a}_1$ (Lemma~\ref{l:misiurewicz}). Finally all these
external rays are connected to the external equipotential of the
graph $\II_n$. \cqfd
 
Let $n\ge1$ and
$\PP_n=\PP_n(\mathbf{a}_0)\subset\PP_{n-1}=\PP_{n-1}(\mathbf{a}_0)$.
  \begin{lemma} \label{l:motion}
There exists a holomorphic motion $h_n\from \PP_n\times
I_{n+1}^{\mathbf{a}_0}\to \ol \C$ such that\,:
\begin{itemize} \item $I_{n+1}^\mathbf{a} = \h_n^\mathbf{a} (I_{n+1}^{\mathbf{a}_0})$ for all
$ \mathbf{a} \in \PP_n$\,;\item  $\h_n$ coincides with $\h_{n-1}$
on $\PP_n \times I^{\mathbf{a}_0}_n$\,; \item
 for every $\mathbf{a} \in
\PP_n$ the following diagram is commutative\,: $$   \begin{CD}
I_{n+1}^{\mathbf{a}_0} @>{\h^\mathbf{a}_n}>> I_{n+1}^\mathbf{a}
\\
              @V{f_{\mathbf{a}_0}}VV               @VV{f_\mathbf{a}}V \\
              I_n^{\mathbf{a}_0}  @>>{\h^\mathbf{a}_{n-1}}> I_n^\mathbf{a}   \end{CD} $$
\end{itemize}  \end{lemma}

  \begin{proof}
By Corollary~\ref{c:holom}, the graph $I_0^{\mathbf{a}_0}$ admits
a holomorphic motion on~$\PP_0(\mathbf{a}_0)$. For
$\mathbf{a}\in\PP_n$ the critical point is not on any of the
graphs $I_k^\mathbf{a}$ for $k\le n+1$ (Corollary~\ref{c:crit}),
so we can pull-back the holomorphic motion of $I_0^{\mathbf{a}_0}$
(by $f^{j}_\mathbf{a}$ with $j\le n+1$) to get  the sequence of
holomorphic motions of the graphs $I_{j}^{\mathbf{a}_0}$ on the
restricted domain $\PP_{j-1}$. By construction they  satisfy the
announced properties.
  \end{proof}

\subsection{Relation between graphs and para-graphs}~\label{s:realtiong-pg}
\begin{lemma} \label{l:homeo}
The following map $H_n$ is a homeomorphism. $$ \H_n \from \left\{
  \begin{aligned}
   \PP_n \cap \II_{n+1} & \longrightarrow P^{\mathbf{a}_0}_n \cap I^{\mathbf{a}_0}_{n+1} \\
   \mathbf{a} 
 & \longmapsto (\h^\mathbf{a}_n)_{}^{-1} (f_\mathbf{a}(-\mathbf{a}))
  \end{aligned} \right. $$
  \end{lemma}

  \begin{proof}
  For $\mathbf{a}$ in
    $\PP_n\cap \II_{n+1}$, $(\h^\mathbf{a}_n)^{-1} (f_\mathbf{a}(-\mathbf{a}))$ is well-defined
      by Lemma~\ref{l:motion} and Lemma~\ref{l:graph}.
  The  image $H_n(\PP_n\cap\II_{n+1})$ is clearly included in $I_{n+1}^{\mathbf{a}_0}$.
Moreover for $\mathbf{a}=\mathbf{a}_0$,  the critical value
$f_{\mathbf{a}_0}(-\mathbf{a}_0)$ belongs to the puzzle piece
$P_n^{\mathbf{a}_0}$. Therefore $f_{\mathbf{a}}(-\mathbf{a})$
belongs to the (open) puzzle piece bounded by
$h_n^\mathbf{a}(\partial P_n^{\mathbf{a}_0})$, since
$f_{\mathbf{a}}(-\mathbf{a})$ and
$I_{n}^{\mathbf{a}}=h_n^\mathbf{a}(I_{n}^{\mathbf{a}_0})$ move
continuously and never meet when $\mathbf{a} \in \PP_n$
(Corollary~\ref{c:crit}). Hence, $H_n(\PP_n\cap \II_{n+1})\subset
P_n^{\mathbf{a}_0}$ since $(\h^\mathbf{a}_n)_{}^{-1}$  is
injective on $I^\mathbf{a}_{n+1}$.

By construction, the map $H_n$ is clearly a homeomorphism on the
rays and equipotentials of $\II_{n+1}\cap \PP_n$ that are in
$\HH_\infty$. We prove now that it is injective in $\HH$.
 Assume by  contradiction that $\UU_1, \UU_2$ are two
 connected components of $\HH$ and that there exist parameters $\mathbf{a}_1$, $\mathbf{a}_2$
 on two rays $\RR_{\UU_1}(t_1)$, $\RR_{\UU_2}(t_2)$ respectively such that
 $H_n(\mathbf{a}_1)=H_n(\mathbf{a}_2)$.
Since $\PP_n$ is a simply connected region
(corollary~\ref{c:1conx}) that avoids the center of all the
components of $\HH_i$ for $0\le i\le n$, we can define functions
$r_{\UU_1}(\mathbf{a})$ and $r_{\UU_2}(\mathbf{a})$ on $\PP_n$ by
Lemma~\ref{l:defcentre}. Since for $j=1,2$ the critical value
$f_{\mathbf{a}_j}(-\mathbf{a}_j)$ belongs to
$R^{r_{\UU_j}(\mathbf{a}_j)}_{\mathbf{a}_j}(t_j)$
(corollary~\ref{c:corresp}), $H_n(\mathbf{a}_j)$ then  belongs to
$(h_n^{\mathbf{a}_j})^{-1}\left(R^{r_{\UU_j}(\mathbf{a}_j)}_{\mathbf{a}_j}(t_j)\right)=R^{r_{\UU_j}(\mathbf{a}_0)}_{\mathbf{a}_0}(t_j)$.
Since $H_n(\mathbf{a}_1)=H_n(\mathbf{a}_2)$, the two rays have a
common point so coincide, and  $\UU_1=\UU_2$, $t_1=t_2$. The same
arguments work (simpler) for the injectivity on the equipotentials
in $\HH\cap \II_{n+1}$.

 To achieve the proof of injectivity,
 it is enough to  show that $H_n$ is injective
  on $\II_{n+1}\cap \partial \CC$.
 Thus, we consider two distinct rays
$\RR_{\UU_1}(t_1), \RR_{\UU_2}(t_2) \subset \II_{n+1}$ landing at
points $\mathbf{a}_1, \mathbf{a}_2 \in \PP_n$ such that
$\H_n(\mathbf{a}_1)=\H_n(\mathbf{a}_2)$, with $\UU_j$ connected
components of $\HH\cup \HH_\infty$. As before, the corresponding
dynamical rays $R^{r_{\UU_j}(\mathbf{a}_j)}_{\mathbf{a}_j}(t_j)$
for $j=1,2$
 are pulled back  by the holomorphic motion  to the rays
$R^{r_{\UU_j}(\mathbf{a}_0)}_{\mathbf{a}_0}(t_j)$.
  These rays land at a common point\,:
$\H_n(\mathbf{a}_1)=\H_n(\mathbf{a}_2)$. Since this point is
eventually repelling and not eventually critical,
 this situation is possible only if  one (at least)
 of the centers  $r_{\UU_1}(\mathbf{a}_0)$ or $r_{\UU_2}(\mathbf{a}_0)$ is at $\infty$
 (Lemma~\ref{l:aboutdistray1} and Lemma~\ref{l:aboutdistray0}),
  or in the trivial case where $\UU_1=\UU_2$, $t_1=t_2$.
 Say $r_{\UU_2}(\mathbf{a}_0)=\infty$, so $\UU_2=\HH_\infty$ and
  $r_{\UU_2}(\mathbf{a})=\infty$  for every $\mathbf{a}\in \PP_n$.
Now, pulled back to the dynamical plane of $\mathbf{a}_1$ (through
the holomorphic motion) the rays $R^{\infty}_{\mathbf{a}_1}(t_2)$
and $R^{r_{\UU_1}(\mathbf{a}_1)}_{\mathbf{a}_1}(t_1)$ still land
at a common point,
 by  Corollary~\ref{c:stable}. This common point is
 $h_n^{\mathbf{a}_1}(H_n(\mathbf{a}_1))=f_{\mathbf{a}_1}(-\mathbf{a}_1)$.
This implies that $ \RR_\infty(t_2)$ lands at $\mathbf{a}_1$
(Lemma~\ref{l:misiurewicz}) and therefore that
$\mathbf{a}_1=\mathbf{a}_2$.

The surjectivity follows from the same kind of arguments. The map
is clearly surjective on the part of  $I^{\mathbf{a}_0}_{n+1}\cap
P^{\mathbf{a}_0}_n$ which is in $\ol{B_{\mathbf{a}_0}(\infty)}$
(the closure of the basin of $\infty$) by construction and by
Lemma~\ref{l:misiurewicz}. Now let $z$ be a point in
$P^{\mathbf{a}_0}_n$ on $I_{n+1}^{\mathbf{a}_0}\cap \wt
B_{\mathbf{a}_0}$, so $z\in R^r_{\mathbf{a}_0}(t)$ for some center
$r$ and angle $t$. Let $z_1$  be the landing point of this ray. By
Remark~\ref{r:aboutdistrayext} and Corollary~\ref{c:crit} there is
only one external ray in $I_{n+1}^{\mathbf{a}_0}$ that also lands
at $z_1$, say $R_{\mathbf{a}_0}^\infty(t')$. Thus, the ray
$\RR_\infty(t')$ belongs to the para-graph $\II_{n+1}$ (by the
surjectivity). It lands at a parameter $\mathbf{a}_1$ and there is
only one ray in $\II_{n+1}$ that also lands at $\mathbf{a}_1$ (by
the injectivity of $H_n$). Consider the simple arc formed by the
union of  this ray, $\{\mathbf{a}_1\}$ and $\RR_\infty(t')$. Its
image by $H_n$ must contain $z$, by the injectivity of $H_n$ and
since there is no other branch of $I_{n+1}^{\mathbf{a}_0}$ at
$z_1$. This reasoning extends to the parameters of $\II_{n+1}$
lying on equipotentials.
  \end{proof}

  \begin{corollary} \label{c:critpiece}
 Let $\mathbf{a} \in \PP_{n-1}$. The parameter $\mathbf{a}$ belongs
to the annulus $\PP_{n-1} \setminus \ol\PP_{n}$ if and only if the
critical value~$f_\mathbf{a}(-\mathbf{a})$ lies in
$P^\mathbf{a}_{n-1} \setminus \ol{ C_n^\mathbf{a}}$ where
$C_n^\mathbf{a} $ is the puzzle piece bounded by $
{h^\mathbf{a}_{n-1}(\partial P_n^{\mathbf{a}_0})}$.
  \end{corollary}
  \begin{proof} By definition, $\mathbf{a}\in \PP_{n-1}$ implies that
  $f_\mathbf{a}(-\mathbf{a})\in P_{n-1}^\mathbf{a}$.
In the proof of Lemma~\ref{l:homeo} we showed  that the piece
$P_{n}^\mathbf{a}$ is bounded by $h^\mathbf{a}_{n-1}(\partial
P_n^{\mathbf{a}_0})$, for $\mathbf{a} \in \PP_{n}$. Thus  if
$f_\mathbf{a}(-\mathbf{a})$ lies in $P^\mathbf{a}_{n-1} \setminus
\ol{ C_n^\mathbf{a}}$, then the parameter $\mathbf{a}$ is in
$\PP_{n-1} \setminus \ol\PP_{n}$.

 Conversely, suppose $\mathbf{a}$ belongs to
$\PP_{n-1} \setminus \ol\PP_{n}$ and  denote by
$\PP_n(\mathbf{a})$ the new para-puzzle piece of depth $n$
containing $\mathbf{a}$ (if $\mathbf{a}\in
\XX_{n-1}\setminus\XX_n$ the result is clear).
 We construct a continuous  path  $\mathbf{a}_t\subset \PP_{n-1}$
 joining $\mathbf{a}_0$ to~$\mathbf{a}_1 = \mathbf{a}$,  crossing
$\partial \PP_{n}$ (resp. $\partial \PP_n(\mathbf{a})$)  at
exactly one point~$\mathbf{a}_{t_0}$ (resp. $\mathbf{a}_{t_1}$) on
equipotentials of $\II_n$ and avoiding $\II_{n} \setminus
\{\mathbf{a}_{t_0},\mathbf{a}_{t_1}\}$. For this, we connect
$\mathbf{a}_0$ by a path in $\PP_n$ to a point $\mathbf{a}'\in
\HH_\infty\cap \PP_n$ then we follow the ray containing
$\mathbf{a}'$ and cross $\partial \PP_n$ at an equipotential, we
then take an equipotential contained in
$\XX_{n-1}\setminus\ol{\XX_n}$ and join $\mathbf{a}_1$ by a ray
entering $\PP_n(\mathbf{a})$ and a piece of path, as before,
inside $\PP_n(\mathbf{a})$.

Thus, for $t<t_0$, the critical value
$f_{\mathbf{a}_t}(-\mathbf{a}_t)$ belongs to the puzzle piece
bounded by $h_{n-1}^{\mathbf{a}_t}(\partial P_n^{\mathbf{a}_0})$.
Moreover, since the parameter path~$\mathbf{a}_t$ crosses
$\II_{n}$ at $t=t_0$ on a ray (changing so the value of its
potential) the critical value~$f_{\mathbf{a}_t}(-\mathbf{a}_t)$
goes out of the piece bounded by $ {h^\mathbf{a}_{n-1}(\partial
P_n^{\mathbf{a}_0})}$ when $t$ passes over~$t_0$. Then
Lemma~\ref{l:homeo} insures that, for $t_0<t<t_1$, the critical
value does not cross $ I_{n}^{\mathbf{a}_t}$ again. Hence the
critical value $f_{\mathbf{a}_t}(-\mathbf{a}_t)$ is outside the
piece bounded by $ {h^{\mathbf{a}_t}_{n-1}(\partial
P_n^{\mathbf{a}_0})}$. Now at $t_1$, the critical value
$f_{\mathbf{a}_{t_1}}(-\mathbf{a}_{t_1})$ belongs to $
I_{n}^{\mathbf{a}_{t_1}}$, but not to
${h^{\mathbf{a}_{t_1}}_{n-1}(\partial P_n^{\mathbf{a}_0})}$ since
$H_{n-1}(\mathbf{a}_{t_1})\notin \partial P_n^{\mathbf{a}_0}$ by
the injectivity of $H_n$ (see Lemma~\ref{l:homeo}). Thus as before
going inside $\PP_n(\mathbf{a})$ along a ray, the critical value
enters a new puzzle piece which  is not bounded by
${h^{\mathbf{a}_t}_{n-1}(\partial P_n^{\mathbf{a}_0})}$.

  \end{proof}

  \begin{corollary} \label{c:nondeg}
If $\ol{P^{\mathbf{a}_0}_{n}} \subset P^{\mathbf{a}_0}_{n-1}$ then
$\ol\PP_{n} \subset \PP_{n-1}$.
  \end{corollary}

  \begin{proof}
Assume, to get a contradiction, that $\partial \PP_{n} \cap
\partial \PP_{n-1}\neq\emptyset$. In the graph $I^{\mathbf{a}_0}_0(\theta)$
two rays never converge to the same point, they don't neither  in
$I^{\mathbf{a}_0}_n(\theta)$ by pullback nor in $\II_n$ by
isomorphism (or Lemma~\ref{l:misiurewicz}). Therefore, the
intersection $\partial \PP_{n} \cap\partial \PP_{n-1}$ contains at
least a part of external rays. Let $\mathbf{a}_1$ be such an
intersection point contained in some ray $\RR_\infty(t)$. So
$d^nt\in \{d^j\eta, j\ge 0\}$ but also $d^{n-1}t\in \{d^j\eta,
j\ge 0\}$. Through a path in $\PP_n$ one can go from
$\mathbf{a}_1$ to the ``center'' $\mathbf{a}_0$ of $\PP_n$ without
crossing the graph $\II_n$ (outside $\mathbf{a}_1$). Since the
graph admits a holomorphic motion in $\PP_n$, the critical value
$f_\mathbf{a}(-\mathbf{a})$ enters the piece $P_n^\mathbf{a}$
which is bounded by $h_{n-1}^\mathbf{a}(\partial
P_n^{\mathbf{a}_0})$ (this is clear taking a path which starts by
some part of equipotential). On the path and by holomorphic
motion, the ray $R_\mathbf{a}^\infty(t)\cap X_n^\mathbf{a}$ is in
the boundary $\partial P_n^\mathbf{a}=h_{n-1}^\mathbf{a}(\partial
P_n^{\mathbf{a}_0})$. Therefore $R^\infty_{\mathbf{a}_0}(t)$
belongs to $I_n^{\mathbf{a}_0}$. But, since $d^{n-1}t\in
\{d^j\eta, j\ge 0\}$, the ray $R^\infty_{\mathbf{a}_0}(t)$ is also
in $I_{n-1}^{\mathbf{a}_0}$. Therefore, $\partial
P_n^{\mathbf{a}_0}\cap\partial P_n^{\mathbf{a}_0}\neq \emptyset$
which contradicts the hypothesis.
\end{proof}

\begin{corollary}\label{c:homeobord} If $\partial \PP_{n+1} \subset
\PP_n$, then $ H_n$ induces a homeomorphism between
$\partial\PP_{n+1}$ and $\partial P_{n+1}^{\mathbf{a}_0}$.
\end{corollary}
\proof Since $H_n$ is not well defined on $\partial \PP_n$ we need
the assumption that $\partial \PP_{n+1} \subset \PP_n$. Then the
result follows from Lemma~\ref{l:homeo}. \cqfd


\subsection{Graphs and renormalization}~\label{s:graploccon}
In the dynamical plane, the graphs of Definition~\ref{d:puzzle}
are used to prove the following Theorem (see~\cite{Ro1})\,:
\begin{thmi}[\cite{F,Ro1}]\footnote{For
completeness, we will sketch its proof  in
section~\ref{s:annexe2}.} The boundary of $B_\mathbf{a}$, as well
as of any connected component of $\wt
  B_\mathbf{a}$,  is a Jordan curve  for any $\mathbf{a}\in \C$.
\end{thmi}
It follows (by Remark~\ref{r:grotzsch}) from
Proposition~\ref{p:dyn} which is a formulation of Yoccoz' Theorem
in the context of the family $(f_\mathbf{a})$.

\begin{remark}\label{r:grotzsch}
Let $p$ be a point in some region $U\subset \C$. If a sequence of
disjoint annuli $A_{k}$ are homotopic in $U\setminus \{p\}$ and
satisfy $\sum_{k\ge0}\mod A_{k}=\infty$ then the diameter of
$U_k$, the connected component of $U\setminus A_k$ containing $p$,
shrinks to $0$. \end{remark} \proof This is a direct consequence
of the following classical results (see \cite{A})\,:

$\bullet$ Gr\"otzsch inequality\,: $\mod A \ge \sum_i \mod A_i$
when $A_i$ are disjoint sub-annuli of $A$ homotopic to $A$\,;

$\bullet$ for any compact $K$ contained in a disk $D$,  if  the
annulus $D \setminus K$ has
 infinite modulus, then $K$ is just a point.
\cqfd

\begin{proposition}\label{p:dyn}
Let $\mathbf{a}\in  \CC  \setminus  \R^-$. There exist
$\epsilon=\pm1$ and $l_0$ such that for $l\ge l_0$  the puzzle
defined by $I^\mathbf{a}_0(\theta)$ with
$\theta=\frac{\epsilon}{(d-1)^l-1}$ gives a sequence of
non-degenerate
 annuli $A^\mathbf{a}_{n_i}$ satisfying\,:
  \begin{enumerate}
 \item for $i\ge 1$, $A^\mathbf{a}_{n_i}= P^\mathbf{a}_{n_i}
\setminus \ol{ P^\mathbf{a}_{n_i+1}}$\,; so
$A^\mathbf{a}_{n_i}$surrounds the critical value for $i\ge 1$ but
maybe not  for $i=0$\,; \item $f_\mathbf{a}^{n_i-n_0}$ induces a
non-ramified covering map from $\ol{A^\mathbf{a}_{n_i}}$ onto
$\ol{A^\mathbf{a}_{n_0}}$; \item either $\sum_{i\ge0}\mod
A^\mathbf{a}_{n_i}=\infty$ {\rm(}where $\mod A^\mathbf{a}_{n_i}$
denotes the modulus of $A^\mathbf{a}_{n_i}${\rm)} or there exists
$k>1$ such that
 $f_\mathbf{a}^k\from P^\mathbf{a}_{n+k} \to P^\mathbf{a}_n$ is a
quadratic-like map for every large $n$.
  \end{enumerate}
 \end{proposition}
The proof of this proposition  can be found in~\cite{Ro1} as a
consequence of Lemma~2.9, Lemma~2.10 and Theorem~1.10 (Theorem of
Yoccoz) in this article. Similar formulations can be found
in~\cite{F,M3}. It will be used several times latter.

 \begin{definition}
A proper holomorphic map $f\from U\to V$ is {\it
  quadratic-like}
  if $U,V$ are  topological disks with $\ol U \subset V$ and if the degree
  of $f$ is $2$.

A map $f$ is said to be {\it renormalizable} if there exist disks
$U,V$ and some integer $k>1$ such that $f^k \from U \to V$ is
quadratic-like and if the orbit of the unique critical point $x$
of the restriction $f^k_{|U}$  stays in $U$, {\it i.e.}
$f^{kn}(x)\in U$ for all $n\ge0$. The integer $k$ is called the
{\it period}.
\end{definition}

\begin{lemma}~\label{l:dynrenor} The map $f_\mathbf{a}$
is renormalizable if and only if there exists $l_1>l_0$ such that
for $l\ge l_1$  the second case of the alternative of
Proposition~\ref{p:dyn}.3. occurs for the graphs
$I_0^\mathbf{a}(\theta)$ defined by
$\theta=\frac{\epsilon}{(d-1)^l-1}$.
\end{lemma}
\proof First we suppose that $f_\mathbf{a}$ satisfies the second
case of the alternative of Proposition~\ref{p:dyn}.3. Then, there
exists $n_0$ and $k>1$  such that
 $f_\mathbf{a}^k\from P^\mathbf{a}_{n+k} \to P^\mathbf{a}_n$ is a
quadratic-like map for every $n\ge n_0$. Since
$f_\mathbf{a}(-\mathbf{a})$ belongs to $P_n^\mathbf{a}$,  the
critical point $-\mathbf{a}$ belongs to
$f_\mathbf{a}^{k-1}(P_{n+k}^\mathbf{a})=P^\mathbf{a}_{n-1}(f_\mathbf{a}^k(-\mathbf{a}))$
for every $n\ge n_0$
  and the restriction  $f_\mathbf{a}^{k-1}\from P^\mathbf{a}_{n+k} \to P^\mathbf{a}_{n-1}(f_\mathbf{a}^k(-\mathbf{a}))$ is
  a homeomorphism.
  Let denote by $x_\mathbf{a}$ the unique preimage of $-\mathbf{a}$ by the
  restriction. Since the pieces are nested,  this preimage does not
  depend on $n$ and $x_\mathbf{a}\in P_i$ for every $i\ge k$.
  The point $x_\mathbf{a}$ is the critical point of the restriction
   $f_\mathbf{a}^k\from P^\mathbf{a}_{n+k} \to P^\mathbf{a}_n$  for  every $n\ge n_0$.
To prove that $f_\mathbf{a}$ is renormalizable it is enough to see
that if we fix some $n\ge n_0$,  the image
$f_\mathbf{a}^{ki}(x_\mathbf{a})$ belongs to
$P^{\mathbf{a}}_{n+k}$ for all $i\ge 0$. This follows from the
fact that $x_\mathbf{a}\in P^\mathbf{a}_{n+k+ki}$ for every $i\ge
0$ so that $f^{ki}_\mathbf{a}(x_\mathbf{a})\in
P^\mathbf{a}_{n+k}$. \vskip 0.5em

Now we assume that  $f_\mathbf{a}$ is renormalizable. Let
$K_\mathbf{a}$ denote its filled Julia set and $k$ the period. We
prove first that the intersection $K_\mathbf{a}\cap\partial
B_\mathbf{a}$ contains at most one point.
 Assume that there is at least two points in this intersection
 but also that $K_\mathbf{a}$ is not included in $\partial B_\mathbf{a}$.
 Then there is a bounded
connected component in  $\C \setminus (\partial B_\mathbf{a}\cup
K_\mathbf{a})$ so there are points on the boundary  of this
connected component (and also on $\partial B_\mathbf{a}$) which
are not in $\partial B_\mathbf{a}(\infty)$\,; this is not possible
for a polynomial. If $K_\mathbf{a}\subset\partial B_\mathbf{a}$ we
would deduce by iteration that $\partial
B_\mathbf{a}=K_\mathbf{a}$, since $\partial B_\mathbf{a}$ is a
Jordan curve\,; this is not possible for a polynomial (namely for
$f_\mathbf{a}^k$). Now we can prove that $K_\mathbf{a}$ is
included in the puzzle pieces $P_n^\mathbf{a}(-\mathbf{a})$ as
follows. If $K_\mathbf{a}$ is cut by $\partial
P_n^\mathbf{a}(-\mathbf{a})$
 there are some rays in $\tilde B_\mathbf{a}$ converging to points of $K_\mathbf{a}$,
 so by iteration some ray of the graph of depth $0$
 converges to a point of $K_\mathbf{a}\cap \partial B_\mathbf{a}$.
  The intersection point has to be fixed by $f_\mathbf{a}^k$ (else there is
 more than one point in the intersection). For $l_1>l_0$, the
 rays of the graphs  defined for $l\ge l_1$ are not  $k$-periodic
 so cannot converge to $K_\mathbf{a}\cap\partial B_\mathbf{a}$.
Therefore $K_\mathbf{a}$ is included in all the puzzle pieces
$P_i^\mathbf{a}(-\mathbf{a})$, so $f_\mathbf{a}^k\from
P_{n+k}^\mathbf{a}\to P_n^\mathbf{a}$ is quadratic-like and we are
in the second case of the alternative of
Proposition~\ref{p:dyn}.3. for those graphs.\cqfd

\begin{definition}\
A set $\M_0$ is a {\it copy } of $\M$ if there exists a
homeomorphism $\chi$ and an integer $k>1$ (the {\it period}) such
that \begin{enumerate} \item$\M_0=\chi^{-1}(\M)$,\item
$\chi^{-1}(\partial \M)\subset\partial\CC$ and for every
$\mathbf{a}\in \M_0$, \item $f_\mathbf{a}$  is renormalizable near
the critical point $-\mathbf{a}$ with $f_\mathbf{a}^k$
topologically conjugated to $z^2+\chi(\mathbf{a})$ on
neighbourhoods of the filled Julia sets.
\end{enumerate}
\end{definition}
\begin{proposition}\label{p:mandelbrot}
 If $f_{\mathbf{a}_0}$ is renormalizable, $
{\mathbf{M}}_{\mathbf{a}_0}=\Cap_{n\ge 0}\PP_n(\mathbf{a}_0)$ is a
copy of $\M$.
 \end{proposition}
\proof  Since $f_{\mathbf{a}_0}$ is renormalizable, there exist
$l_1>l_0$ such that the graphs defined in Lemma~\ref{l:dynrenor}
satisfy the second case of the alternative of
Proposition~\ref{p:dyn}.3. We prove that $\{f_\mathbf{a}^k\from
P^\mathbf{a}_n \to P^\mathbf{a}_{n-k}, \ \mathbf{a} \in
\PP_n(\mathbf{a}_0)\}$ form a {\it Mandelbrot-like family}.

\noindent For $n\ge n_0$ and $\PP_n=\PP_n(\mathbf{a}_0)$, we
consider the mapping $\mathbf{f}\from \WW'\to \WW$ defined by
 $\WW=\{(\mathbf{a},z) \mid \mathbf{a} \in \PP_n
,\ z\in  P_{n-k}^\mathbf{a}\}$, $\WW'=\{(\mathbf{a},z) \mid
\mathbf{a} \in \PP_n
,\ z\in P_{n}^\mathbf{a}\}$ and
$\mathbf{f}(\mathbf{a},z)=(\mathbf{a},f^k_\mathbf{a}(z))$. They
form an analytic family of quadratic-like maps in the sense of
Douady and Hubbard \cite[p.304]{DH2} since they satisfy the
following three properties\,:

 \begin{itemize}
 \item the map $\mathbf{f}\from \WW' \to \WW$ is holomorphic and
proper\,;

\item the holomorphic motion of the disk $P_n^\mathbf{a}$, resp.
$P_{n-k}^\mathbf{a}$, is a homeomorphism  between $\WW$', resp.
$\WW$, and $\PP_n\times \D$ which  is fibered over $\PP_n$ (since
$\mathbf{a}\in\PP_n$)\,;

\item the projection  $\ol \WW' \cap \WW \to  \PP_n$ ({\it i.e.}
the first coordinate) is proper, since  $\ol \WW' \cap
\WW=\{(\mathbf{a},z) \mid \mathbf{a} \in \PP_n,\ z\in
\ol{P_n^\mathbf{a}}\}$.
\end{itemize}
Let $\M_{\mathbf{f}}=\{\mathbf{a} \mid K(f_\mathbf{a}^k) \hbox{ is
connected}\}$ denote the connectedness locus of $\mathbf{f}$,
where $K(f_{\mathbf{a}_0}^k)=\Cap_{i \ge 0}
(f_{\mathbf{a}_0}^k)^{-i}(P^{\mathbf{a}_0}_n)$ denote its filled
Julia set. Then $\M_{\mathbf{f}}$ coincides with
$\M_{\mathbf{a}_0}$. Indeed, for $\mathbf{a} \in
\M_{\mathbf{a}_0}$, the critical point $-\mathbf{a}$, and its
orbit by $f_\mathbf{a}$, does never cross the graphs
$I_j^\mathbf{a}$ ($j\ge 0$) since $\mathbf{a}$ belongs to every
para-puzzle pieces.
 Therefore the critical point $x_\mathbf{a}$  of
$f_\mathbf{a}^k|_{P_n^\mathbf{a}}$ never escapes the piece
$P_n^\mathbf{a}$ (by iteration by $f_\mathbf{a}^k$). Hence
$K(f_\mathbf{a}^k)$ is connected and
$\mathbf{a}\in\M_{\mathbf{f}}$. Conversely, for $\mathbf{a} \in
\PP_n\setminus\ol\PP_{n+1}$, the critical value belongs to
$A_n^\mathbf{a}$ (by Corollary~\ref{c:critpiece}). Thus
$f_\mathbf{a}^k(f_\mathbf{a}(-\mathbf{a}))$ is not in
$P_{n-k}^\mathbf{a}$ and therefore the critical point of
$f_\mathbf{a}^k$ escapes the domain\,; then
 the filled Julia set is not connected anymore so that
  $\mathbf{a}\notin \M_{\mathbf{f}}$.

Moreover, by Corollary~\ref{c:nondeg} and Proposition~\ref{p:dyn},
there exists a sequence $n_i$
 such that $\ol{\PP_{n_i+1}}\subset \PP_{n_i}$. Thus $\M_{\mathbf{a}_0}$
is also the intersection of the closed pieces\,:
$\M_{\mathbf{a}_0}=\Cap_{n\ge0}\ol \PP_n$ and therefore  is
compact.

 Now, the theory of Mandelbrot-like families  of Douady
 and Hubbard  (see \cite{DH2}, Theorem II.2, Propositions II.14 and IV.21)
 gives a continuous map $\chi \from \PP_n \to \C$ such that
 the maps  $f^k_\mathbf{a}$ and  $z^2+\chi(\mathbf{a})$ are quasi-conformally conjugated on a
 neighbourhood of the filled Julia sets, for every $\mathbf{a} \in \PP_n$.

 Moreover, since $\M_{\mathbf{f}}$ is compact, the map $\chi$ induces a
 homeomorphism between   $\M_{\mathbf{f}}$  and the Mandelbrot
 set~$\M$  if we are in the following situation (see \cite{DH2})\,:
 for a closed disk $\Delta\subset \PP_n$
 containing $\M_{\mathbf{f}}$ in its interior,  the quantity $f^k_\mathbf{a}(x_\mathbf{a})-x_\mathbf{a}$,
 (where   $x_\mathbf{a}$ denotes the critical point of $f^k_\mathbf{a}|_{P_n^\mathbf{a}}$ should turn exactly
 once around $0$ when $\mathbf{a}$ describes $\partial \Delta$. We verify
 this property now.

 Take some piece $\PP_p(\mathbf{a}_0)=\Delta$,  compactly contained   $\PP_n(\mathbf{a}_0)$
 (see  Corollary~\ref{c:nondeg}).
 It is a topological disk  containing  $\MM_{\mathbf{a}_0}$ in its interior.
  To compute the degree on $\partial \Delta$ of
$\gamma(\mathbf{a})=f_\mathbf{a}^k(x_\mathbf{a})-x_\mathbf{a}$ we
make a homotopy of this curve $\gamma$ to the curve
$H_{p-1}(\mathbf{a})-x_{\mathbf{a}_0}$ as follows. The critical
point $x_\mathbf{a}$ of $f_\mathbf{a}^k$ satisfies that
$f_\mathbf{a}^{k-1}(x_\mathbf{a})=-\mathbf{a}$, so
$\gamma(\mathbf{a})=f_\mathbf{a}(-\mathbf{a})-x_\mathbf{a}$. Let
$h(\mathbf{a},z)=h_{p-1}(\mathbf{a},z)-x_\mathbf{a}$, then
$\gamma(\mathbf{a})=h(\mathbf{a}, H_{p-1}(\mathbf{a}))$. Assume
that $\PP_{p-1}$ is a round disk (if not  use a conformal
representation)\,; then the homotopy is simply
$G(t,\mathbf{a})=h(\mathbf{a}_0+t(\mathbf{a}-\mathbf{a}_0),H_{p-1}(\mathbf{a}))$
joining $G(0,\mathbf{a})=H_{p-1}(\mathbf{a})-x_{\mathbf{a}_0}$ and
$G(1,\mathbf{a})= f_\mathbf{a}(-\mathbf{a})-x_\mathbf{a}$.

Since $H_{p-1}$ is a homeomorphism from $\partial \PP_p$ to
$\partial P_p^{\mathbf{a}_0}$ (piece that surrounds
$x_{\mathbf{a}_0}$), the degree of
$H_{p-1}(\mathbf{a})-x_{\mathbf{a}_0}$ around $0$ is exactly $1$,
when $\mathbf{a}$ describes $\partial \PP_p$. \cqfd

\begin{proi}{\rm \bf\ref{p:periodicrenormalizable}.}
If $f_{\mathbf{a}}$ has a periodic point $x\neq 0$ of multiplier
$\rho$ with $|\rho|\le1$, then $f_{\mathbf{a}}$ is renormalizable
near $x$ and $\mathbf{a}$ belongs to a copy of $\M$.
\end{proi}
\proof Since $x\neq 0$, it is not in $\wt B_\mathbf{a}$ and, since
it is not eventually repelling, it is not on any of the graphs. So
we can consider the sequence $(P_n^\mathbf{a}(x))$ of puzzle
pieces containing $x$. Since $x$ is periodic, this  sequence of
pieces is periodic {\it i.e.},
$f_\mathbf{a}^k(P^\mathbf{a}_{n+k}(x))=P^\mathbf{a}_{n}(x)$ for
any large $n$ and for some  $k>1$. Choose the smallest $k$ with
this property. There exists some $i\le k$ such that the critical
point $-\mathbf{a}$ belongs to the piece
$P^\mathbf{a}_{n}(f_\mathbf{a}^i(x))$, for every sufficiently
large $n$. Otherwise the map $f_\mathbf{a}^k\from
P^\mathbf{a}_{n+k}(x)\to P^\mathbf{a}_{n}(x)$ would be invertible
and its inverse $g\from P^\mathbf{a}_{n}(x)\to
P^\mathbf{a}_{n}(x)$ either would be an automorphism or has an
attracting fixed point (by Schwarz' Lemma). This is not possible
since on the one hand $|g'(x)|\ge 1$ and on the other hand the
sequence $P^\mathbf{a}_{n}(x)$ is strictly decreasing ($\exists \
n \mid P^\mathbf{a}_{n+k}(x)=g(P^\mathbf{a}_{n}(x))\neq
P^\mathbf{a}_{n}(x)$). This integer $i$ is independent of $n$
since the pieces $P_n^\mathbf{a}(f_\mathbf{a}^j(-\mathbf{a}))$ are
disjoint for $j<k$. Therefore the map $f^k_\mathbf{a} \from
P^\mathbf{a}_{n+k}(f^i_\mathbf{a}(x))\to
P^\mathbf{a}_{n}(f^i_\mathbf{a}(x))$ is quadratic-like. Hence we
are in the situation of Lemma~\ref{l:dynrenor} where we proved
that $\mathbf{a}$ belongs to a copy of the Mandelbrot set $\M$
(see Proposition~\ref{p:mandelbrot}). \cqfd

\begin{cri}{\rm \bf\ref{c:H2}.}
Any  bounded hyperbolic component either
 is a connected component of $\mathcal H$ or a hyperbolic component of
  a copy of  $\mathbf M$.
\end{cri}
\proof Let $\UU$ be a hyperbolic component  which is not in $\HH$.
For $\mathbf{a}\in\UU$, the map $f_\mathbf{a}$ has an attracting
periodic cycle, which is not the fixed point $0$. Thus we are in
the situation of Proposition~\ref{p:periodicrenormalizable} so
that the parameter $\mathbf{a}$ belongs to a copy of $\M$. \cqfd

 \begin{cri}{\rm \bf\ref{c:brujno}.} If $f_\mathbf{a}$ has a periodic point $x$ of
multiplier $\lambda=e^{2i\pi \theta}$ with $\theta\in \R\setminus
\Q$, then $f_\mathbf{a}$ is linearizable near $x$ if and only if
$\theta \in \mathcal B$. Moreover, if $\theta\notin \mathcal B$
there exist periodic cycles in any neighbourhood of $x$.
\end{cri}
\proof The map  $f_\mathbf{a}$ is renormalizable by
Proposition~\ref{p:periodicrenormalizable}. So there is a
homeomorphism   that conjugates $f_\mathbf{a}^k$ to a quadratic
polynomial $z^2+\chi(\mathbf{a})$ on a neighbourhood of its Julia
set (see~\cite{DH2}). The multiplier at the fixed points are the
same by Na\"{\i}shul' Theorem (see~\cite{N}). So the result
follows from Yoccoz' and Brjuno's work (see~\cite{Y}). \cqfd


\section{Local connectivity}


Fix $\mathbf{a}_0\in \partial\CC\cap \SS$.
 Take $l_0\ge l_1$ given by Proposition~\ref{p:dyn} and
$\theta\in\left\{\pm \frac{1 }{(d-1)^l-1}\right\}$ with $l\ge
l_0$.

Recall that the sequence of graphs $I_n^{\mathbf{a}_0}(\theta)$
and the para-graph $\II_n(\theta)$ associated in
Definition~\ref{d:ppuzzle} satisfy the following properties\,:
\begin{itemize}
\item The sequence of puzzle pieces $P^{\mathbf{a}_0}_n$
containing the critical value is well-defined since the critical
value $f_{\mathbf{a}_0}(-\mathbf{a}_0)$ is on none of the graphs
$I^{\mathbf{a}_0}_n$, $n\ge 0$. \item The sequence of para-puzzle
pieces~$\left(\PP_n\right)_{n\in \N}$ containing~$\mathbf{a}_0$ is
well-defined  by Lemma~\ref{l:homeo} (since the parameter
$\mathbf{a}_0$ also never belongs  to a
 graph $\II_n$).
\item There exists a sequence of (non-degenerate) annuli
$\left(A^{\mathbf{a}_0}_{n_i}\right)_{i\in \N}$ such that, for
$i\ge 1$, $A^{\mathbf{a}_0}_{n_i} = P^{\mathbf{a}_0}_{n_i}
\setminus \ol{ P^{\mathbf{a}_0}_{n_i+1}}$
  (so surrounds the critical value $f_{\mathbf{a}_0}(-\mathbf{a}_0)$) and
 the map $f_{\mathbf{a}_0}^{n_i-n_0}\from \ol{A^{\mathbf{a}_0}_{n_i}}\to
\ol{A^{\mathbf{a}_0}_{n_0}}$ induces  a non-ramified covering map
(Proposition~\ref{p:dyn})\,; \item The  annuli
$\AA_{n_i}=\PP_{n_i}\setminus\ol {\PP_{n_i+1}}$  are
non-degenerate (Corollary~\ref{c:nondeg}) and
surround~$\mathbf{a}_0$.
\end{itemize}

 
\subsection{Tools for proving local connectivity: estimation of moduli and connectivity questions}

The next Proposition follows from Shishikura' trick to compare
moduli of annuli.

  \begin{proposition} \label{p:modules}
There exists a constant $K>1$ such that, for $ i \ge 0$, $$ \frac
1K \mod A_{n_i}^{\mathbf{a}_0} \le \mod \AA_{n_i} \le K \mod
A_{n_i}^{\mathbf{a}_0} \,.$$
  \end{proposition}

\proof The idea is to get a $K$-quasi-conformal homeomorphism
between $\AA_{n_i} $ and $A_{n_i}^{\mathbf{a}_0}$ extending  the
map $H_{n_i}$ (via S\l odkovksi's Theorem and the dynamical
covering).

\noindent Fix $n \in \{n_i,\ i\ge0\}$ and let  $d_n$ be the degree
of $f_{\mathbf{a}_0}^{n-n_0} \from A_n^{\mathbf{a}_0} \to
A_{n_0}^{\mathbf{a}_0}$. For every $\mathbf{a} \in \PP_n$ we
define, for $m\le n$, the dynamical annuli $A_m^\mathbf{a}$,
``motion'' of $A_m^{\mathbf{a}_0}$, by  the connected component of
$P^\mathbf{a}_{m} \setminus {h^\mathbf{a}_{m}(\partial
P_{m+1}^{\mathbf{a}_0})}$ that intersects $I_{m+1}^\mathbf{a}$ in
its interior. By definition (Lemma~\ref{l:motion}) the following
diagram is commutative for $\mathbf{a}\in \PP_n$. $$
\begin{CD}
\partial A_n^{\mathbf{a}_0} @>{\h^\mathbf{a}_n}>> \partial A_n^\mathbf{a} \\
 @V{f_{\mathbf{a}_0}^{n-n_0}}VV                      @VV{f_\mathbf{a}^{n-n_0}}V \\
\partial A_{n_0}^{\mathbf{a}_0}  @>>{\h^\mathbf{a}_{n_0}}> \partial A_{n_0}^\mathbf{a}
  \end{CD}$$
 Thus $f_\mathbf{a}^{n-n_0}$ maps $\ol{A_n^\mathbf{a}}$ to $\ol{A_{n_0}^\mathbf{a}}$ and induces
 a non-ramified covering. Indeed, the critical value $f_\mathbf{a}(-\mathbf{a})$
 remains outside $I_n^\mathbf{a}$ so that the  critical point $-\mathbf{a}$
 and all its preimages cannot enter $\ol{A_n^\mathbf{a}}$.
At depth $n_0$ we extend the holomorphic motion
$h_{n_0}\from\PP_{n_0}\times \partial A^{\mathbf{a}_0}_{n_0} \to
\ol \C$ by S\l odkowski's Theorem \cite{Sl} to  a holomorphic
motion $ \wt\h_{n_0} \from \PP_{n_0} \times \ol\C \longrightarrow
\ol\C $.
 For every $\mathbf{a}\in \PP_{n_0}$, the map $\wt\h^\mathbf{a}_{n_0}$ is a
 $K_\mathbf{a}$-quasi-conformal homeomorphism, with
 $K_\mathbf{a}=\frac{1+|\phi(\mathbf{a})|}{1-|\phi(\mathbf{a})|}$ where $\phi \from \PP_{n_0}\to\D$
 is a  conformal representation sending $\mathbf{a}_0$ to $0$.
 For every $\mathbf{a} \in \PP_n$ the homeomorphism
$\wt\h^\mathbf{a}_{n_0} \from \ol{A_{n_0}^{\mathbf{a}_0}} \to
\ol{A_{n_0}^\mathbf{a}}$ lifts---via the holomorphic covering maps
$f_{\mathbf{a}_0}^{n-n_0}$ and $f_\mathbf{a}^{n-n_0}$---to a
quasi-conformal homeomorphism $ \wt\h^\mathbf{a}_n \from
\ol{A_n^{\mathbf{a}_0}} \longrightarrow \ol{A_n^\mathbf{a}} $ with
the same dilatation  $K_\mathbf{a}$.
 Moreover the identity $f_\mathbf{a}^{n-n_0} \circ \h^\mathbf{a}_n =\h^{\mathbf{a}}_{n_0} \circ f_{\mathbf{a}_0}^{n-n_0}$
 ensures that the map $\wt\h_n
\from \PP_n \times \ol {A_n^{\mathbf{a}_0}} \to  \ol\C$, \
$(\mathbf{a},z) \mapsto \wt\h_n^\mathbf{a}(z) $, \ is a
holomorphic motion that extends $\h_n$. From
Corollary~\ref{c:critpiece} we know that $\mathbf{a}$ belongs to $
 {\AA_n}$ if and only if $f_\mathbf{a}(-\mathbf{a})$
 belongs to $ { A_n^\mathbf{a}}$ so that $\mathbf{a}$ belongs to $
\ol {\AA_n}$ if and only if $f_\mathbf{a}(-\mathbf{a})$
 belongs to $\ol { A_n^\mathbf{a}}$ or
equivalently $\mathbf{a}$ belongs to ${\AA_n} \cup \partial
\PP_{n+1}$ if and only if $(\wt\h^\mathbf{a}_n )^{-1}
(f_\mathbf{a}(-\mathbf{a}) )\in A_n^{\mathbf{a}_0}\cup
\partial P_{n+1}^{\mathbf{a}_0} $. Therefore
the following map $\wt\H_n$ is well-defined
$$ \wt\H_n \from \left\{
  \begin{aligned}
   {\AA_n} \cup \partial \PP_{n+1}&{}\longrightarrow
    {A_n^{\mathbf{a}_0}\cup \partial P_{n+1}^{\mathbf{a}_0}}, \\
   \mathbf{a} \quad&{}\longmapsto
   \wt\H_n(\mathbf{a}) = \left( \wt\h^\mathbf{a}_n \right)^{-1} (f_\mathbf{a}(-\mathbf{a}))
  \end{aligned} \right.  $$
From~\cite[IV.3]{DH2}  the map $\wt H_n$ is $K_n$-quasi-regular
with $ K_n = \sup \bigl\{ K_\mathbf{a},\ \mathbf{a} \in \PP_n
\bigr\}$ (see also~\cite{Ro2}). Moreover $\wt H_n$ is a bijection
since it agrees with $H_n$ on $\partial \PP_{n+1}$
(Lemma~\ref{l:homeo}). Therefore $\wt H_n$ is a
$K$-quasi-conformal homeomorphism from $\AA_n$ to
$A_n^{\mathbf{a}_0}$ with $K=\sup \bigl\{ K_\mathbf{a},\
\mathbf{a} \in \ol\PP_{n_0+1} \bigr\}< +\infty$. The result then
follows.\cqfd

In the rest of this subsection we prove the connectedness of the
intersection of the para-puzzle pieces with $\partial \CC$ and
$\partial \UU$.
\begin{lemma}\label{l:parapiece}
Let $\UU$ be a  connected component of $\HH$. For every $n\ge 0$,
the intersection of $\UU$ and $\ol\PP_n$  is a sector of $\UU$
bounded by $$\
\partial \PP_n\cap \mathcal U=\left(\XX_n\cap (\RR_{\mathcal
U}(t_n)\cup \RR_{\mathcal U}(t'_n))\right)\cup \wt{\mathcal
E}_{\mathcal U}(v_n),$$ where $\wt{\mathcal E}_{\mathcal U}(v_n)$
is a part of the equipotential $\mathcal E_{\mathcal U}(v_n)$ and
$t_n\le t'_n$. Moreover, as $n$ tends to infinity,  $v_n\to 0$ and
$t_n,t'_n$ converge to a common value.
\end{lemma}
\proof Since every ray of $\partial \PP_n\cap \mathcal U$ is
associated  to an external ray, it is not possible to have more
than two rays in $\partial \PP_n\cap \mathcal U$ (the internal and
external rays are connected by equipotentials). So $t_n,t'_n$ are
consecutive angles in $\Theta_n$. We prove that $t_n-t'_n\to 0$
for $\HH_0$\,; the proof is the same for any $\UU$. By definition
of the graphs, any  puzzle piece of depth $0$ intersects $\HH_0$
under an angular sector of width less than $1/d$. Therefore the
puzzle pieces of greater depth have rays in $\HH_0$ whose angles
are consecutive angles  divided  $d-1$. Thus $|t'_{n}-t_{n}|\le
\frac{1}{d(d-1)^n}$ , so $\lim t_n=\lim t_n'=t$ (since both
sequences  are monotone).\cqfd

   \begin{lemma}\label{l:connexe}
For every  connected component $\UU$ of $ \HH$, the intersection
$\ol \PP_n\cap\partial \UU$ is a connected set for all $n\ge 0$.
  \end{lemma}
\proof From Lemma~\ref{l:parapiece}, we know that
 $\partial \UU \cap \ol \PP_n$ is just the decreasing intersection
of the compact connected sets $\ol{\Phi_{\UU}(S^n_k)}$ where
$S^n_k$ is the sector in the disk between the angles $t_n$ and
$t_n'$ and of potential less than $v^n_k$, where
$(v^n_k)_{k\in\N}$ s a sequence  which tends to $0$ with $k$.
Therefore it is compact and connected.\cqfd

 \begin{lemma}\label{l:connexeC}
For every  $n\ge 0$, the intersection $\ol \PP_n\cap\partial \CC$
is  connected.
  \end{lemma}
  \proof The  property of the para-puzzle pieces we use here are  to be a
  disk  whose  boundary is a succession of arcs
   of  the following form\,:
  a part of an equipotential  in $\HH_\infty$ followed by a part of a ray in $\XX_n$
  converging to a point of $\partial \CC$ and another
  part of a ray in $\XX_n$ followed by a part of an equipotential  in $\HH$
  (by Remark~\ref{r:aboutdistrayext} and  Lemma~\ref{l:homeo}).
  We denote the property  by $(*)$.
  Fix $n$ and let $\GG_n$ be the
  bounded connected component of $\ol\C\setminus
  \EE_\infty(1/d^n)$.
Consider $B(k)$ the set of   disks ${\DD}\subset \ol\GG_n$
satisfying Property (*) and such that $\GG_n\setminus  \DD$ has
$k$ connected components.

  We prove  by recurrence on $k$ that for any disk $\DD \in B(k)$,
  $\ol\DD\cap \partial \CC$ is connected.

Let $\DD\in B(1)$, we prove by contradiction that
$\ol\DD\cap\partial \CC$ is connected. Let $V$ be the complement
$\GG_n\setminus\ol\DD$. Since $ \DD$ belongs to $B(1)$  there is
only two parts of equipotentials in its boundary\,: one part of
$\EE_\infty(1/d^n)$ and  one part of an equipotential in a
component $\UU$ of $\HH$. The complement $V$ has the same
property. Therefore the intersection $\partial V\cap  \CC$ is
reduced to the landing points of the external rays so is included
into $\partial \UU\cap \CC$. Assume now that $\ol\DD\cap\partial
\CC$ is not connected\,: $\ol\DD\cap\partial \CC=A\sqcup B$ where
$A$ and $B$ are non empty, closed and disjoint. The intersection
$\ol\DD\cap \partial \UU$ is connected (it is the intersection of
a decreasing sequence of connected compacts as in the previous
Lemma~\ref{l:connexe}) so we can assume that it is contained in
$A$. Therefore $A'=(V\cap\partial\CC) \cup A$ is closed since the
closure of $V\cap\partial\CC$ is included in $A$. Moreover $A'$ is
 disjoint from $B$ and $A'\sqcup B=\partial \CC$. This
contradicts the fact that $\partial \CC$ is connected.

Now fix some integer $k\ge 1$. Assume that we have proved the
result for $B(i)$ with $i\le k$. Take $ \DD$ a  disk of $B(k+1)$.
There exists at least a  connected component, $\VV$, of
$\GG_n\setminus\ol \DD$ whose boundary intersects
$\EE_\infty(1/d^n)$ under exactly  one component. Then $ \DD\cup
\ol\VV$ is a  disk in $B(k)$. So, with the same argument used
before for $k=1$, if $\ol \DD \cap \partial \CC$ is not connected
then $(\ol \DD\cup \ol\VV)\cap
\partial \CC$ is  not connected neither. This gives the
contradiction.  \cqfd

Now we can conclude in the non-renormalizable case.
 \begin{lemma}\label{c:modules}
If the map $f_{\mathbf{a}_0}$ is not renormalizable, then
$\partial \CC$ and $\partial \UU$ are  locally connected at
$\mathbf{a}_0$, where $\UU$ is any connected component of $\HH$.
   \end{lemma}
\proof   By Lemma~\ref{l:dynrenor} if the map $f_{\mathbf{a}_0}$
is not renormalizable it  satisfies the first alternative of
Proposition~\ref{p:dyn}.3.
 The sequence of annuli considered in
 Proposition~\ref{p:dyn} has the property that
 $\sum_{i\ge0}\mod A^{\mathbf{a}_0}_{n_i}=\infty$.
 Hence
 $\sum_{i\ge0}\mod\AA_{n_i}=\infty$ by  Proposition~\ref{p:modules}.
Thus the diameter of  $\PP_{n_i}$  shrinks to $0$ by
Remark~\ref{r:grotzsch}.

 Finally, $\ol{\PP_{n_i}}\cap \partial \UU$, resp.   $\ol{\PP_{n_i}}\cap \partial \CC$,
 form a basis of
 connected neighbourhoods of $\mathbf{a}_0$ in $\partial \UU$, resp. in $\partial \CC$, by
 Lemma~\ref{l:connexe}, resp.  Lemma~\ref{l:connexeC}
 (since $\PP_{n_i}$ is a neighbourhood of
 $\mathbf{a}_0$).
\cqfd

\begin{crr}\label{c:Clc} $\partial \CC$ is locally connected at $\mathbf{a}_0$
as soon as $\mathbf{a}_0$
 does not belong to a copy of $\M$.
\end{crr}
\proof By  Proposition~\ref{p:mandelbrot} we are in the first
alternative of Proposition~\ref{p:dyn}.3., so  the result follows
from Lemma~\ref{c:modules}.\cqfd

\subsection{Local connectivity of $\partial \HH_0$ and Wakes of $\HH_0$}\label{s:H0interM}
\begin{proposition}\label{p:H0interM}
 Let  $\mathbf{a}_0\in
\partial\HH_0\cap\SS$. If $f_{\mathbf{a}_0}$ is renormalizable of period $k$,
there exist $k$-periodic angles $t$ and $\zeta,\zeta'$ {\rm(}by
multiplication by $d-1$ and by $d$ respectively{\rm)} such that
the rays $\RR_0(t)$, $\RR_\infty(\zeta)$ and $ \RR_\infty(\zeta')$
converge to $\mathbf{a}_0$. Moreover the curve $\RR_\infty(\zeta)
\cup \RR_\infty(\zeta')\cup\{\mathbf{a}_0\}$ separates
$\M_{\mathbf{a}_0}\setminus\{\mathbf{a}_0\}$ from $\HH_0$ where  $
{\mathbf{M}}_{\mathbf{a}_0}=\Cap_{n\ge 0}\PP_n(\mathbf{a}_0)$. In
the dynamical plane, the rays $R^0_{\mathbf{a}_0}(t)$,
$R^\infty_{\mathbf{a}_0}(\zeta)$,
$R^\infty_{\mathbf{a}_0}(\zeta')$ converge to the same point which
is  a parabolic periodic point and
$R_{\mathbf{a}_0}^\infty(\zeta)\cup
R_{\mathbf{a}_0}^\infty(\zeta')\cup\{\mathbf{a}_0\}$ separates
$f_{\mathbf{a}_0}(-\mathbf{a}_0)$ from $B_{\mathbf{a}_0}$.
 \end{proposition}
 \begin{figure}[!h]\vskip 0cm
\centerline{\input{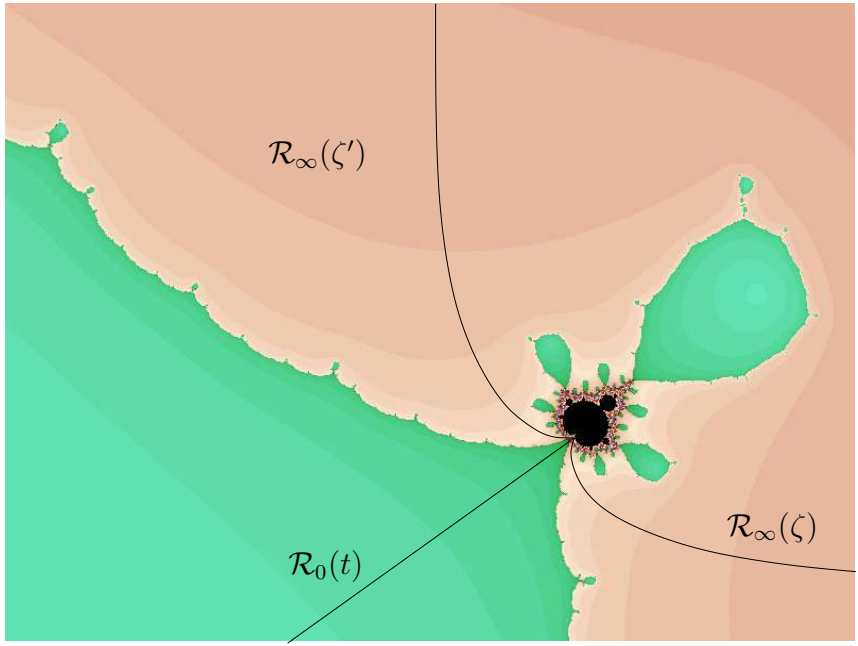}} \caption{Separation of
$\M_{\mathbf{a}_0}$ (copy of $\M$) from $\HH_0$ by
rays.}\label{f:mandelloin}
  \end{figure}

 \proof
 The para-puzzle piece
$\PP_n=\PP_n(\mathbf{a}_0)$ intersects $\HH_0$ since it
 is a neighbourhood of $\mathbf{a}_0\in \partial \HH_0$.
In particular, its boundary $\partial \PP_n$ contains two rays in
$\HH_0$, say  $\RR_0(t_n), \RR_0(t'_n)$, landing at parameters
called $\mathbf{a}_n,\mathbf{a}'_n$ respectively
 and two external rays  say $\RR_\infty(\zeta_n), \RR_\infty(\zeta'_n)$ also
 landing at  $\mathbf{a}_n$ and $\mathbf{a}'_n$ respectively.
 The sequences of angles  $(t_n), (t'_n)$
 converge  to some common value $t$ with $t_n\le t \le t'_n$ (
 see Lemma~\ref{l:parapiece}).
Moreover, the sequence of intervals $(\zeta_n,\zeta'_n)$ is
decreasing, so the angles $\zeta_n,\zeta'_n$ converge to some
values $\zeta,\zeta'$.

The boundary of the dynamical puzzle pieces $P_{n}^\mathbf{a}$ is
given for the  subsequence ${n_i}$ such that $\partial
\PP_{n_i+1}\subset \PP_{n_i}$ by  the bijection described in
Corollary~\ref{c:homeobord}.  It contains parts of the rays
$R^0_{\mathbf{a}_0}(t_{n_i}),R^0_{\mathbf{a}_0}(t'_{n_i})$ and
$R^\infty_{\mathbf{a}_0}(\zeta_{n_i}),
R^\infty_{\mathbf{a}_0}(\zeta'_{n_i})$ landing at points
$z_{n_i}(\mathbf{a}_0), z'_{n_i}(\mathbf{a}_0)$ respectively.
 Since $f_{\mathbf{a}_0}$ is renormalizable, the puzzle pieces are
``periodic'' {\it i.e.},
$f_{\mathbf{a}_0}^k(\ol{P_{n+k}^{\mathbf{a}_0}})=\ol{P_n^{\mathbf{a}_0}}$
by Lemma~\ref{l:dynrenor}. Hence the angles satisfy the relation
 $(d-1)^k\theta_{n+k}=\theta_n\hbox{ mod } 1$ for the rays in
 $B_\mathbf{a}$
 and  $d^k\theta_{n+k}=\theta_n\hbox{ mod } 1$ for the external
 rays.
This is clear for the rays in $B_{\mathbf{a}_0}$ since there  are
only two  rays in $\partial P_n^{\mathbf{a}_0}\cap
B_{\mathbf{a}_0}$\: and it follows for the external rays since
they converge to the same points $z_n,z_n'$.
 Then the angles $t$, $\zeta, \zeta'$ satisfy the same relations
 and therefore are of the form
 $\frac{p}{(d-1)^k-1}$, $\frac{q}{d^k-1}$, $\frac{q'}{d^k-1}$ respectively.

We prove now that $\RR_0(t)$, $\RR_\infty(\zeta)$ and
$\RR_\infty(\zeta')$ converge to $\chi^{-1}_{\mathbf{a}_0}(1/4)$
where $\chi_{\mathbf{a}_0}$  is the homeomorphism that maps
$\M_{\mathbf{a}_0}$ to $\M$. The proof is the same for the three
rays, we do it for $\RR_0(t)$.
The  ray $\RR_0(t)$ 
converges, since $t$ is rational, to some parameter
$\mathbf{a}_1$. For every $n\ge 0$ the part of the ray
$\RR_0(t)\cap\XX_n$ is in $ \PP_n$, because of the bijection
between $\partial \PP_{n_i}$ and $\partial P_{n_i}^{\mathbf{a}_0}$
and of Lemma~\ref{l:parapiece}. Thus $\mathbf{a}_1$ belongs to
$\M_{\mathbf{a}_0}$.
 In the dynamical plane, the  ray $R^0_{\mathbf{a}_1}(t)$
lands at a periodic point $z_1$, since $t$ is periodic. Its
period, say $k$, is the period of $t$ (two rays in
$B_{\mathbf{a}_1}$ cannot land at the same point). The point $z_1$
cannot be eventually critical since $t$ is periodic. It is
parabolic  by Lemma~\ref{l:douady} since  the ray
$R^0_{\mathbf{a}_1}(t)$ is not stable in any neighbourhood of
$\mathbf{a}_1$ by the following remarks. For $\mathbf{a}\in
\RR_0(t)$ near $\mathbf{a}_1$, the critical value is on
$R_\mathbf{a}^0(t)=R_\mathbf{a}^0((d-1)^kt)$. Moreover  the
critical point $-\mathbf{a}$ is in
$f_\mathbf{a}^{k-1}(P_{n+k}^\mathbf{a})$ since $f_{\mathbf{a}_1}$
is renormalizable and $\mathbf{a}$ is very close to
$\mathbf{a}_1$. So $-a$ is on the preimage of $R_\mathbf{a}^0(t)$
that belongs to $f_\mathbf{a}^{k-1}(P_{n+k}^\mathbf{a})$, {\it
i.e.} on $R_\mathbf{a}^0((d-1)^{k-1}t)$. Thus the ray
$R_{\mathbf{a}_1}^0((d-1)^{k-1}t)$, as well as its iterated
preimage
 $R_{\mathbf{a}_1}^0(t)$, is not stable.
Then the point $z_1$ is, for the return map $f_{\mathbf{a}_1}^k$,
a parabolic fixed point. Moreover its multiplier is $1$ since the
ray $R^0_{\mathbf{a}_1}(t)$ is fixed (by  $f_{\mathbf{a}_1}^k$).
Therefore, under the bijection $\chi_{\mathbf{a}_0}$ the parameter
$\mathbf{a}_1$ corresponds to the cusp of $\M$, {\it i.e.}
$\mathbf{a}_1= \chi^{-1}_{\mathbf{a}_0}(1/4)$.

 Finally, the  three rays
$\RR_0(t)$, $\RR_\infty(\zeta)$, $\RR_\infty(\zeta')$ converge to
the same parameter  $\mathbf{a}_1=\chi^{-1}_{\mathbf{a}_0}(1/4)$
of $\M_{\mathbf{a}_0}$ and, the proof above shows that
$\M_{\mathbf{a}_0}\subset \PP_n$ is in the connected component of
$\C \setminus (\ol\RR_0(t_n)\cup \ol\RR_0(t'_n)\cup
\ol\RR_\infty(\zeta_n)\cup \ol\RR_\infty(\zeta'_n))$ that contains
$\mathbf{a}_0$. Therefore, $\M_{\mathbf{a}_0}$ is in  the closure
of the connected component of $\C \setminus
(\ol\RR_\infty(\zeta)\cup \ol\RR_\infty(\zeta'))$   not containing
$\HH_0$. (Misiurewicz parameters in $\M_{\mathbf{a}_0}$ are
accessible by external rays $\RR_\infty(\theta)$ with $\theta\in
(\zeta,\zeta')$.) Then the only possible intersection between
$\M_{\mathbf{a}_0}$ and $\ol \HH_0$ is the cusp $\mathbf{a}_1$.
Therefore $\mathbf{a}_1=\mathbf{a}_0$.

The three rays $R_{\mathbf{a}_0}^0(t)$,
$R_{\mathbf{a}_0}^\infty(\zeta)$,
$R_{\mathbf{a}_0}^\infty(\zeta')$ converge
 to  points  in $\cap P_n^{\mathbf{a}_0}$, fixed by $f_{\mathbf{a}_0}^k$. But $K_{\mathbf{a}_0,k}=\cap P_n^{\mathbf{a}_0}$
 contains only one fixed point with rotation number~$1$ called $\beta$.
 Thus the three rays land at the same point\,: $\beta$.
 They  separate $K_{\mathbf{a}_0,k}$ from $B_{\mathbf{a}_0}$, {\it i.e.} $f_{\mathbf{a}_0}(-\mathbf{a}_0)$
 from  $B_{\mathbf{a}_0}$. Indeed, any eventually repelling periodic point in $K_{\mathbf{a}_0,k}$
 (for instance $\beta'$ the preimage of $\beta$ by $f_{\mathbf{a}_0}^k$) is accessible
  by an external ray whose angle is between $\zeta_n$ and $\zeta'_n$,
  so at the limit between $\zeta$ and $\zeta'$.
 \cqfd

\begin{definition}
A parameter, $\mathbf{a}$, is called  {\it parabolic} (or of {\it
parabolic type}) if $f_\mathbf{a}$ has a parabolic periodic point.
\end{definition}

\begin{corollary}\label{c:H0interM}
Any $\mathbf{a}_0\in \partial\HH_0$, for which $f_{\mathbf{a}_0} $
is renormalizable, is the cusp of a copy of $\M$. More precisely,
the intersection $\M_{\mathbf{a}_0}\cap\partial\HH_0$ reduces to
$\{\mathbf{a}_0\}$ for
$\M_{\mathbf{a}_0}=\cap\PP_n(\mathbf{a}_0)$. Moreover
$\mathbf{a}_0=\chi^{-1}_{\mathbf{a}_0}(1/4)$ where
$\M=\chi_{\mathbf{a}_0}(\M_{\mathbf{a}_0})$, so $\mathbf{a}_0$ is
parabolic.
\end{corollary}
\proof This follows from Proposition~\ref{p:H0interM},
its  proof above and the use of symmetries.

 \begin{corollary}\label{c:lcH} The boundary of
$\HH_0$ is locally connected.
 \end{corollary}
\proof It is locally connected at parameters $\mathbf{a}\in
\partial\HH_0$
 which are not renormalizable (Lemma~\ref{c:modules}). For parameters
$\mathbf{a}_0\in \partial\HH_0$ which are renormalizable,  we
consider the sequence $\QQ_n$ of  subsets
 of $\partial\HH_0$ defined by
$\QQ_n=\ol{\PP_n} \cap \partial \HH_0$. These  subsets are
connected neighbourhoods of $\mathbf{a}_0$ in $\partial \HH_0$
(Lemma~\ref{l:connexe}). Moreover they  form a basis
 since $(\Cap_{n\ge 0}\ol{\PP_n} \cap \partial
\HH_0)\subset (\M_{\mathbf{a}_0}\cap\partial
\HH_0)=\{\mathbf{a}_0\}$ by corollary~\ref{c:H0interM}.\cqfd

The following Lemma is used in Theorem~\ref{th:comp}.

\begin{lemma}\label{l:2rayext} Let $\mathbf{a}\in \CC$. For any
point $z\in \partial B_\mathbf{a}$, there is at most two external
rays converging to $z$. Moreover, if $z$ is not eventually
critical and if  there are two external rays converging to $z$,
then $z$ is (eventually) periodic\,; moreover the two rays define
two connected components, each of them contains at least one
critical point of $f_\mathbf{a}^k$ for some $k\ge 0$.
 \end{lemma}
\proof  If $z$ is eventually critical, there is clearly exactly
two external rays converging to it. Assume now that $z$ is not
eventually critical. Consider the closure of all the external rays
converging to $z$ and let $V$ be a connected component of its
complement.  Assume to get a contradiction that there is no
critical point of $f_\mathbf{a}^k$ in $V$ (for every $k\ge0$).
 Then the  iterates $f_\mathbf{a}^k$ restrict to homeomorphisms on
$V$. This contradicts the fact that the map is doubling the angles
on $B_\mathbf{a}(\infty)$ so that the image of $V$ by some
$f_\mathbf{a}^k$ will contain all $B_\mathbf{a}(\infty)$ and
$f_\mathbf{a}^{k+1}$ will no more be injective on $V$.

Assume now that there are at least three external rays converging
to $z$. Let $V_0$ be the component containing the critical point
$0$ and (up to iterating) one component $V_1$ contains the first
inverse image of $B_\mathbf{a}$. Note that  this  implies that
$-\mathbf{a}\in  V_1$. Take $V_2$ a third component. Since $V_2$
contains a critical point of $f_\mathbf{a}^i$ (for some $i\ge 0$)
it is mapped by some iterate of $f_\mathbf{a}$, to $V_1$ (which
contains $-\mathbf{a}$ and the preimage of $B_\mathbf{a}$).
Indeed, all these sectors are mapped to sectors attached to
$\partial B_\mathbf{a}$ as long as they do not contain a critical
point. Therefore $z$ is a periodic point, of period say $k$. Then
the fact that $V_2$ is mapped to $V_1$  contradicts the fact that
$f_\mathbf{a}^k$ preserves the cyclic order of the rays landing at
$z$ since
 there is a finite number of external rays landing at $z$
(see for instance~\cite{Pe}).
 \cqfd
Now we can describe more precisely the boundary of $\partial
\HH_0$\,:

\begin{thmi}{\bf\ref{th:comp}.} Let $\mathbf{a}\in
\partial\mathcal{H}_0\setminus\R^{-}$\,;
there exists a unique parameter ray  in $\HH_0$   landing at
$\mathbf{a}$, say $\RR^s_0(t)$. The following dichotomy
holds\,{\rm:}
\begin{itemize}

\item either there is a unique external parameter ray converging
to $\mathbf{a}$. In this case $f_\mathbf{a}$ is not renormalizable
so that $\mathbf{a}$ does not belong to a copy of  $\M$. Moreover
in the dynamical plane, the ray $R_\mathbf{a}^0(t)$ lands at the
critical value $f_\mathbf{a}(-\mathbf{a})\in
\partial B_\mathbf{a}$ and there is a unique external ray
converging to $f_\mathbf{a}(-\mathbf{a})$\,;

\item or there are exactly
 two external parameter rays
 converging to $\mathbf{a}$. In this case  $\mathbf{a}$  is
  the cusp of a copy of $\mathbf M$.
  Furthermore, in the dynamical plane, the  ray $R_\mathbf{a}^0(t)$
  lands at a parabolic point
on $\partial B_\mathbf{a}$. The angle $t$ is necessarily periodic
by multiplication by $d-1$.

 \end{itemize}
\end{thmi}

Note that in the first case, the angle $t$ can be periodic by
multiplication by $d-1$ in this case $\mathbf{a}$ is a Misiurewicz
parameter. In Proposition~\ref{p:convzero},  we give the exact
conditions on $t\in \Q$ so that $\mathbf{a}$ is of parabolic or of
Misiurewicz type. If $t\in \R\setminus \Q$ we are clearly in the
first case.

\begin{figure}[!h]\vskip 0cm
\centerline{\input{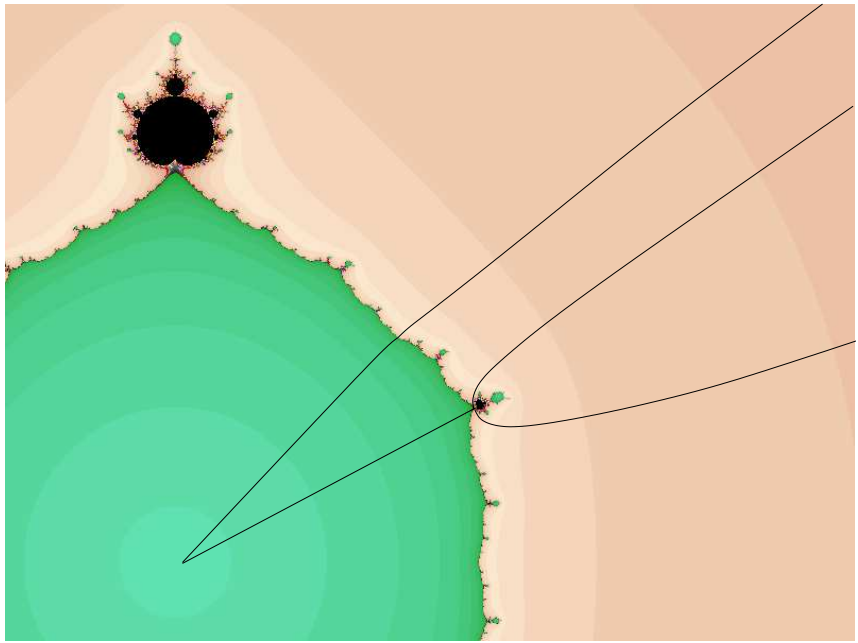}} \caption{Points on $\partial
\HH_0$ and rays converging to them.}\label{f:pararys}
 \end{figure}

\proof We do the proof in several steps.

{\it Any  parameter $\mathbf{a}$ of $\partial \HH_0\setminus \R^-$
is the landing point of a ray in $\HH_0$}\,:

\noindent We consider the fundamental domain $s(\SS)$ containing
$\mathbf{a}$. For $d>3$, the map $\Phi_0$ restricts to a
homeomorphism $\check{\Phi}_0$ from $s(\SS)\cap \HH_0$ onto
$\Delta_d$  (see Proposition~\ref{p:H0}).
 Since the boundary of $\HH_0$ and therefore of $s(\SS)\cap \HH_0$
 (Remark~\ref{r:rayons}) is locally connected,
 the inverse map $\Psi_0=\check{\Phi}_0^{-1}$ extends continuously to
the boundaries\,: $\ol\Psi_0\from \ol{\Delta_d} \to \ol{s(\SS)\cap
\HH_0}$. The analogue statement  for $d=3$ gives a continuous
extension $\ol\Psi_0$  from $\ol\D\setminus\R^+$ into
  $s(\dot{\SS})\cap\ol{\HH_0}$ (one can use a double covering argument to see this).
Therefore the parameter $\mathbf{a}$ on the boundary of $\HH_0$ is
the limit of a ray $\RR^s_0(t)$.

\vskip 0.5em \noindent{\it Now we suppose first that
$f_\mathbf{a}$ is not renormalizable.}

Note that  $\mathbf{a}$ does not  belong to a copy of  $\M$ by the
definition of ``renormalizable''.

 In the dynamical plane, the sequence
of puzzle pieces $(P_n^\mathbf{a})$ shrinks to one point namely
$f_\mathbf{a}(-\mathbf{a})$. Moreover, taking the subsequence
$n_i$ such that $\ol {\PP_{n_i+1}}\subset \PP_{n_i}$, the
existence of the homeomorphism between $\partial \PP_{n_i+1}$ and
$\partial P_{n_i+1}^\mathbf{a}$ preserving angles and potentials
(Corollary~\ref{c:homeobord}) insures that the ray
$R_\mathbf{a}^0(t)$ enters all the puzzle pieces
$P_{n_i+1}^\mathbf{a}$ for $i\ge 0$. Thus $R_\mathbf{a}^0(t)$
converges to $f_\mathbf{a}(-\mathbf{a})$.

By Lemma~\ref{l:2rayext}, there is only one external ray landing
at $z=\va$ since $-\mathbf{a}\in \partial B_\mathbf{a}$. Assume
now that two external rays $\RR_\infty(\xi),\RR_\infty(\xi')$ land
at $\mathbf{a}$.
 These two rays enter any para-puzzle piece $\PP_n(\mathbf{a})$ so by the homeomorphism
of Corollary~\ref{c:homeobord}   the rays
$R^\mathbf{a}_\infty(\xi)$ and $R^\mathbf{a}_\infty(\xi')$ enter
all the pieces $P_n^\mathbf{a}$. Since the intersection $\Cap
_{n\ge 0}P_n^\mathbf{a}$ reduces to $\va$, the rays both converge
to the same point $z=\va$. But we have just seen that this is not
possible.

\vskip 0.5em \noindent{\it Now we consider the second case of the
dichotomy\,: $f_\mathbf{a}$ is renormalizable.}

In this case Proposition~\ref{p:mandelbrot} insures that
$\M_{\mathbf{a}}=\cap\PP_n(\mathbf{a})$ is a copy of $\M$ and
$\mathbf{a}$ is the cusp $\chi^{-1}(1/4)$ where $\chi$ is the
homeomorphism between $\M_{\mathbf{a}}$ and $\M$. There are two
external rays $\RR_\infty(\zeta),\RR_\infty(\zeta')$ converging to
$\mathbf{a}$ by Proposition~\ref{p:H0interM} and in the dynamical
plane the ray $R_\mathbf{a}^0(t)$ converges to a point $z\in
\partial B_\mathbf{a}$ which is a parabolic periodic point. Hence the angle
$t$ is periodic by multiplication by $d-1$.

To prove that there is only two external parameter rays converging
to $\a$ we proceed by contradiction. Assume  that there is a third
ray $\xi$ converging to $\a$. To fixe the ideas assume that the
cyclic order at $\infty$ is $\zeta',\xi,\zeta$. Then the
Mandelbrot copy belongs to one connected component of the
complement of
$\RR_\infty(\xi)\cup\RR_\infty(\zeta)\cup\RR_\infty(\zeta')\cup\{\a\}$,
say the one  containing the rays of angle between $\zeta'$ and
$\xi$. Since the ray $\RR_\infty(\xi)$ enters every para-puzzle
piece, the ray $R_\a^\infty(\xi)$ enters every puzzle piece
$P_n^\a$ by the homeomorphism of Corollary~\ref{c:homeobord}. So
the ray $R_\a^\infty(\xi)$ converges to a point $z$ of $J(f_\a^k)$
(the Julia set of the renormalized map). There exist points  in
$J(f_\a^k)$  that are accessible by external ray of angle
$\xi'\in\Q$ between $\xi $ and $\zeta$ such that neither
 $\frac{t}{d}+\frac{\lfloor\frac{d-1}{2}\rfloor}{d}$
nor $\frac{t}{d}+\frac{\lfloor\frac{d-1}{2}\rfloor+1}{d}$ is
periodic   by multiplication by $d$. Then the ray
$\RR_\infty(\xi')$ lands at a Misiurewicz parameter by
Proposition~\ref{p:convinfini}. This Misiurewicz parameter belongs
to $\M_\a$ since the ray $\RR_\infty(\xi')$ enters every
para-puzzle piece $\PP_n$ (by the homeomorphism of
Corollary~\ref{c:homeobord}). But this contradicts the fact that
$\M_\a$ belongs to connected component  containing the rays of
angle between $\zeta'$ and $\xi$.
\cqfd

\begin{remark} The fact that  $t$ is $k$-periodic by multiplication by $d-1$
 does not imply that $f_\mathbf{a}$ is renormalizable.
\end{remark}
\proof It is possible that the sequence $f_\mathbf{a}^i
(P^\mathbf{a}_{n+k})$ avoids the critical point since there are
other preimages of $f_\mathbf{a}(-\mathbf{a})$ on $\partial
B_\mathbf{a}$ in degree $d>3$. See
Proposition~\ref{p:convzero}.\cqfd

 \begin{definition}\label{d:wake}
We define the {\it wake} $\WW(\mathbf{a}_0)$ of any point
 $\mathbf{a}_0\in \partial \HH_0\cap\SS$ as follows.

If $f_{\mathbf{a}_0}$ is not renormalizable let us take
$\WW(\mathbf{a}_0)=\emptyset$\,; else let $\WW(\mathbf{a}_0)$ be
the connected component of
$$\C\setminus\left(\ol\RR_\infty(\zeta_0)\cup\ol\RR_\infty(\zeta'_0)\right)
\quad \hbox{ containing }\quad \M_{\mathbf{a}_0}\setminus
\{\mathbf{a}_0\}$$ where $\zeta_0$, $\zeta'_0$  are periodic
angles (by multiplication by $d$) such that the rays
$\RR_\infty(\zeta_0)$ and $\RR_\infty(\zeta'_0)$ converge to
$\mathbf{a}_0$. For parameters not in $\SS$ we use the symmetries
to define the wake.

 \end{definition}

 \begin{remark} By Theorem~\ref{th:comp} there are at most  two rays
  converging to a parameter  $\a\in\partial \HH_0$ (those
  defined in Proposition~\ref{p:H0interM}) so that
 the wake is well defined. Moreover since is the landing point of a ray
 $\RR^s_0(t)$, we can also call
 $\WW^s(t)$ the wake $\WW(\mathbf{a})$. Note that the wake of a parameter $\a\in\SS$
 is  not necessarily contained in $\SS$.
 \end{remark}

\begin{lemma}\label{l:wake} For any parameter $\mathbf{a}$ in   $\WW(\mathbf{a}_0)$,
 the rays $R^\infty_\mathbf{a}(\zeta_0)$,
$R^\infty_\mathbf{a}(\zeta_0')$ and $R^0_\mathbf{a}(t_0)$ converge
to the same point which  is repelling of period $k$ {\rm(}the
period of $\M_{\mathbf{a}_0}${\rm)}, where $\zeta_0,\zeta'_0$
define the wake $\WW(a_0)$ and $\RR_0^s(t_0)$ is landing at
$\a_0$. For $\a=\mathbf{a}_0$ these three dynamical  rays also
land at a common  point, which is $k$-periodic parabolic point.
Moreover, for $\mathbf{a}\in
\WW(\mathbf{a}_0)\cup\{\mathbf{a}_0\}$, the critical value is in
the corresponding dynamical wake\,: the connected component of
$\C\setminus(\ol R_\mathbf{a}^\infty(\zeta_0)\cup
 \ol R_\mathbf{a}^\infty(\zeta'_0))$ which does not contain $B_\mathbf{a}$.\end{lemma}

\begin{figure}[!h]\vskip 0cm
\centerline{\input{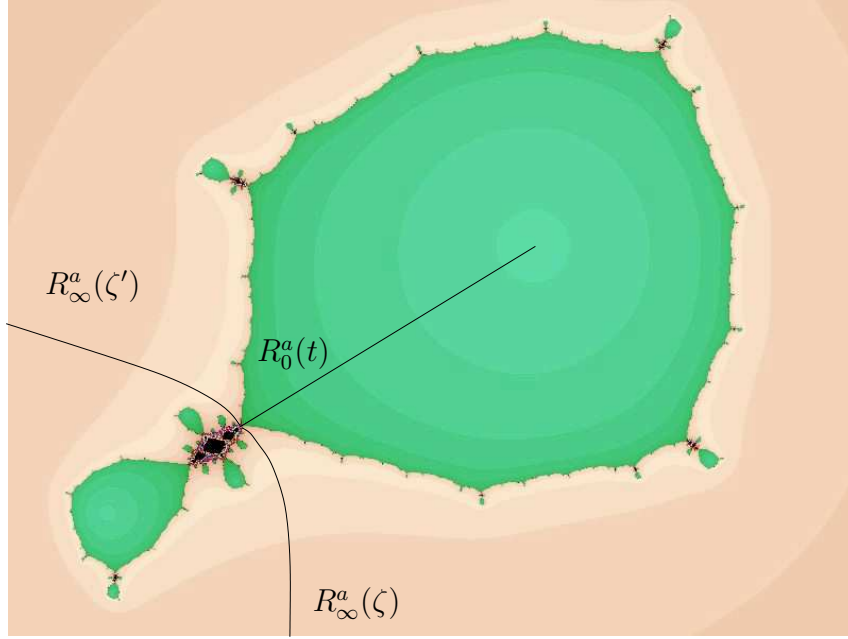}} \caption{Julia set for
$\mathbf{a}$ in the copy $\M_{\a_0}$ of
 figure~\ref{f:mandelloin}.}\label{f:juliam4}
  \end{figure}

\proof  Note first that, for every  parameter $\mathbf{a} \in
\M_{\mathbf{a}_0}$, the three rays $R^\infty_\mathbf{a}(\zeta_0)$,
$R^\infty_\mathbf{a}(\zeta_0')$ and $R^0_\mathbf{a}(t_0)$ converge
to the same $k$-periodic point.  For $\mathbf{a}=\mathbf{a}_0$,
this follows from Proposition~\ref{p:H0interM}. Then it is easy to
check that all the arguments of Proposition~\ref{p:H0interM} go
through for the parameters $\mathbf{a}$ in $\M_{\mathbf{a}_0}$.
Indeed,  the boundary of the puzzle pieces $\partial
P_n^{\mathbf{a}}$ and $\partial P_n^{\mathbf{a}_0}$ are identified
through the holomorphic motion defined on the neighbourhood
$\PP_{n-1}$ of $\M_{\mathbf{a}_0}$.
 
After this remark,  the  proof goes exactly  as point $3$ of
Lemma~\ref{l:holom} and Corollary~\ref{c:stable}, so we give here
just the steps of the argumentation. We consider the set $\Omega$
of parameters $\mathbf{a}$ such that
$R^\infty_\mathbf{a}(\zeta_0)$, $R^\infty_\mathbf{a}(\zeta_0')$
and $R^0_\mathbf{a}(t_0)$ converge to the same point which is a
repelling periodic point  of period $k$. Then $\Omega$ is open and
non empty. Its boundary is included in $PP_k\cup\Cup_{i\ge
0}\left(\RR_0((d-1)^i t_1)\cup U\RR_\infty(d^i\zeta_0)\cup
U\RR_\infty(d^i\zeta'_0)\right)$ (see
definition~\ref{d:pararays}), where $PP_k$ is the set of
parameters $\mathbf{a}$ such that $f_\mathbf{a}$ has a parabolic
point of period $k$ and multiplier $1$.  We claim that, in the
region $\WW(\mathbf{a}_0)$ there is no parameter rays of
$U\RR_\infty(d^i\zeta_0)$, $U\RR_\infty(d^i\zeta'_0)$ and
$\RR_0((d-1)^i t_0)$ with $i\ge 0$. For this  we look first in the
dynamical plane of $f_{\mathbf{a}_0}$. The angles $\zeta_0$,
$\zeta'_0$ are the limits of the sequences $(\zeta_n)_{n\ge 2}$,
$(\zeta'_n)_{n\ge 2}$ defined as follows\,: the two internal rays
$R_{\mathbf{a}_0}^0(t_n)$, $R_{\mathbf{a}_0}^0(t'_n)$ of $\partial
P_n^{\mathbf{a}_0}$ converge to points $z_n,z'_n$  to which are
attached the external rays $R_{\mathbf{a}_0}^\infty(\zeta_n)$,
$R_{\mathbf{a}_0}^\infty(\zeta'_n)$ of $\partial
P_n^{\mathbf{a}_0}$.
There is no iterate of $R_{\mathbf{a}_0}^\infty(\zeta_n),
R_{\mathbf{a}_0}^\infty(\zeta'_n)$ in $Q_n$, the connected
component containing $f_{\mathbf{a}_0}(-\mathbf{a}_0)$ of
$$\C \setminus \left(\ol {R_{\mathbf{a}_0}^0(t_n)} \cup
\ol { R_{\mathbf{a}_0}^0(t'_n)} \cup \ol{
R_{\mathbf{a}_0}^\infty(\zeta_n)}\cup \ol{
R_{\mathbf{a}_0}^\infty(\zeta'_n)} \right).$$ Otherwise, such an
iterate would be attached to an internal ray of some $\partial
P_j^{\mathbf{a}_0}$ for $j\le n$, with angle in $(t_n,t'_n)$. But
this is not possible since $\partial P_n^{\mathbf{a}_0}$ contains
only two rays in $B_{\mathbf{a}_0}$, they are in $\partial Q_n$,
and since $P_n^{\mathbf{a}_0}\subset P_j^{\mathbf{a}_1}$.
 Therefore, there is no element of
$d^i\zeta_0$, $d^i\zeta'_0$ in the segment $(\zeta_n,\zeta'_n)$
``defining'' the region $Q_n$,
 and so neither in  the limit interval  $(\zeta_0,\zeta'_0)$.
Thus there are no rays of $U\RR_\infty(d^i\zeta_0),
U\RR_\infty(d^i\zeta'_0)$, $i\ge 0$ in $\WW(\mathbf{a}_0)$. On the
other hand, there is no point of $PP_k$ in $\WW(\mathbf{a}_0)$\,:
otherwise this would contradict the maximum principle for the
multiplier of the landing point of $R^0_\mathbf{a}(t_0)$ (as in
Lemma~\ref{l:holom}). Therefore, $\WW(\mathbf{a}_0)\subset\Omega$.

 \cqfd


\subsection{Local connectivity of $\partial \UU$
for a component $\UU$  of $\HH\setminus \HH_0$.}
 \begin{proposition}\label{p:UcapM}
Let  $\UU$ be   a connected component of $\HH\setminus \HH_0$. Let
$\mathbf{a}_0\in \partial \UU\cap\SS$ be such that
$f_{\mathbf{a}_0}$ is renormalizable and  denote by
$\chi_{\mathbf{a}_0}$ the homeomorphism between
$\M_{\mathbf{a}_0}$ and $\M$ where $\M_{\mathbf{a}_0}=\cap
\PP_n(\mathbf{a}_0)$ {\rm (}see
Proposition~\ref{p:H0interM}{\rm)}. Then\,:
\begin{itemize} \item $\partial \UU\cap \M_{\mathbf{a}_0}=\{\mathbf{a}_0\}$\,;
 \item $\M_{\mathbf{a}_0}\cap\partial
\HH_0=\chi_{\mathbf{a}_0}^{-1}(1/4):=\mathbf{a}_1$ the cusp of
$\M_{\mathbf{a}_0}$\,; \item $\UU\subset \WW(\mathbf{a}_1)$ the
wake of $\mathbf{a}_1$\,; \item  $\mathbf{a}_0$ is the landing
point of three rays $\RR_\UU(t_1)$, $\RR_\infty(\eta)$ and
$\RR_\infty(\eta')$ where $d^i\eta=\zeta_1$, $d^i\eta'=\zeta_1'$
with $t_1$ and  $\zeta_1,\zeta'_1$ $k$-periodic by multiplication
by $d-1$, and by $d$ respectively. Here $i$ is the depth of $\UU$
{\it i.e.,} $\UU\subset \HH_i\setminus \HH_{i-1}$, and
$t_1,\zeta_1,\zeta'_1$ are associated to $\mathbf{a}_1$
 by Proposition~\ref{p:H0interM}\,; \item The curve $\ol
\RR_\infty(\eta)\cup\ol\RR_\infty(\eta')$ separates
 $\UU$ from $\M_{\mathbf{a}_0}$.
\end{itemize}
  \end{proposition}

\begin{figure}[!h]\vskip 0cm
\centerline{\input{MC5b.pstex_t}} \caption{Intersection of $\UU$
and the copy $\M_0$ of $\M$.}\label{f:mandelcopyinterU}
  \end{figure}

 \proof We will not prove the items in  the order they appear.
 The proof is very similar to that of Proposition~\ref{p:H0interM}.
Using Lemma~\ref{l:defcentre},  we can define  in
$\PP_n=\PP_n(\mathbf{a}_0)$ a holomorphic function $r(\mathbf{a})$
which coincides, for $\mathbf{a} \in \UU$, with the center of
$U(\mathbf{a})$, the connected component of $\wt B_\mathbf{a}$
containing the critical value. Since $\PP_n$ intersects $\UU$,
$\partial \PP_n$ contains two internal rays
$\RR_\UU(\tau_n),\RR_\UU(\tau'_n)$ (Lemma~\ref{l:parapiece}) with
landing points $u_n,u'_n$ respectively, but also external rays
$\RR_\infty(\eta_n),\RR_\infty(\eta'_n)$ landing at $u_n,u'_n$
respectively. Using the homeomorphism of
Corollary~\ref{c:homeobord}, the boundary  $\partial
P_n^{\mathbf{a}_0}$ contains the  part in $\XX_n$ of the rays
$R^{r(\mathbf{a}_0)}_{\mathbf{a}_0}(\tau_n)$,
$R^{r(\mathbf{a}_0)}_{\mathbf{a}_0}(\tau'_n)$ and
$R^{\infty}_{\mathbf{a}_0}(\eta_n)$,
 $R^{\infty}_{\mathbf{a}_0}(\eta'_n)$ (Corollary~\ref{c:corresp}),
 with common landing points, say $z_n,z'_n$ respectively, at
 least for $n$ in  the subsequence $(n_i)_{i\ge0}$ given in
 Proposition~\ref{p:dyn}.

 1. {\it We prove  first that $\partial \PP_n\cap \HH_0\neq \emptyset$ and that
$\M_{\mathbf{a}_0}\cap\partial
\HH_0=\chi_{\mathbf{a}_0}^{-1}(1/4):=\mathbf{a}_1$}\,:

\noindent Since for every $m\ge 0$, $P_m^{\mathbf{a}_0}$
intersects $U(\mathbf{a})$, which is of depth  $i$, the image
$f_{\mathbf{a}_0}^i (P_m^{\mathbf{a}_0})$ intersects
$B_{\mathbf{a}_0}$. Then, for $i\le kj\le m$,  the image
$f_{\mathbf{a}_0}^{kj}(P_m^{\mathbf{a}_0})$ which is the puzzle
piece $P_{m-kj}^{\mathbf{a}_0}$ containing the critical value,
intersects $B_{\mathbf{a}_0}$. Thus, $P_{n}^{\mathbf{a}_0}\cap
B_{\mathbf{a}_0}\neq \emptyset$ for any $n$. By the homeomorphism
of Corollary~\ref{c:homeobord}, we deduce that $\PP_n\cap
\HH_0\neq \emptyset$.  So there is some point $\mathbf{a}_1$ in
the intersection $\Cap_{n\ge 0}\PP_n\cap \ol
\HH_0=\M_{\mathbf{a}_0}\cap \ol\HH_0 $. Applying
Proposition~\ref{p:H0interM} (or Corollary~\ref{c:H0interM}) to
this point $\mathbf{a}_1$, we deduce that
$\M_{\mathbf{a}_0}\cap\partial \HH_0$ is reduced to $\mathbf{a}_1$
and is the cusp of $\chi_{\mathbf{a}_0}^{-1}(1/4)$ of
$\M_{\mathbf{a}_0}$ since $\M_{\mathbf{a}_1}=\M_{\mathbf{a}_0}$.
Indeed, $\mathbf{a}_1\in \M_{\mathbf{a}_0}=\cap\PP_n$, the pieces
$\PP_n(\mathbf{a}_0)$ and $\PP_n(\mathbf{a}_1)$ coincide.

2.  {\it We prove here that  $\PP_{n}$ contains in its boundary
the part in $\XX_n$ of the rays $\RR_{0}(\tau_{n+i})$,
$\RR_{0}(\tau'_{n+i})$, $\RR_{\infty}(\zeta_{n+i})$,
 $\RR_{\infty}(\zeta'_{n+i})$ with $\zeta_n=d^i\eta_n$ and
 $\zeta'_n=d^i\eta'_n$, at least for infinitely many $n\in \N$}\,:

\noindent We prove it for the dynamical puzzle piece
$P_n^{\mathbf{a}_0}$ and then use the homeomorphism of
Corollary~\ref{c:homeobord}. We have seen  in point 1 above that
the puzzle piece $P_{n+i}^{\mathbf{a}_0}$ contains in its boundary
the rays $R^{r(\mathbf{a}_0)}_{\mathbf{a}_0}(\tau_{n+i})$,
$R^{r(\mathbf{a}_0)}_{\mathbf{a}_0}(\tau'_{n+i})$ and also
$R^0_{\mathbf{a}_0}(t_{n+i})$, $R^0_{\mathbf{a}_0}(t'_{n+i})$.
Since $f_{\mathbf{a}_0}^i(P_{n+i}^{\mathbf{a}_0})$ is the piece
$P_n^{\mathbf{a}_0}(f_{\mathbf{a}_0}^i(f_{\mathbf{a}_0}(-\mathbf{a}_0)))$,
it  contains only two rays of $B_{\mathbf{a}_0}$ in its boundary,
 so that the rays $R^{r(\mathbf{a}_0)}_{\mathbf{a}_0}(\tau_{n+i})$ and
$R^{0}_{\mathbf{a}_0}(t'_{n+i})$ have the same image  by
$f^i_{\mathbf{a}_0}$. In particular,
$(d-1)^{i}t'_{n+i}=\tau_{n+i}=t'_{n}$. Therefore the puzzle piece
$f_{\mathbf{a}_0}^{i-1}(P_{n+i}^{\mathbf{a}_0})$
 contains the critical point $-\mathbf{a}_0$ since it is simply connected
 and $f_{\mathbf{a}_0}$ maps its boundary with degree two on its image.
Thus
$f_{\mathbf{a}_0}^i(P_{n+i}^{\mathbf{a}_0})=P_n^{\mathbf{a}_0}$ so
that $i$ is a multiple of $k$. The piece $P_{n}^{\mathbf{a}_0}$
contains in its boundary the rays
$R^{0}_{\mathbf{a}_0}(\tau_{n+i})$,
$R^{0}_{\mathbf{a}_0}(\tau'_{n+i})$, with end points $z_{n+i},
z'_{n+i} $ respectively, to which converge the external rays
$R^{\infty}_{\mathbf{a}_0}(\zeta_{n+i})$,
 $R^{\infty}_{\mathbf{a}_0}(\zeta'_{n+i})$
 with $\zeta_{n+i}=d^i\eta_{n+i}$, $\zeta'_{n+i}=d^i\eta'_{n+i}$.

Using the homeomorphism of
 Corollary~\ref{c:homeobord} we deduce  that
 the para-puzzle piece $\PP_{n}$ contains in its boundary the rays $\RR_{0}(\tau_{n+i})$,
  $\RR_{\infty}(\zeta_{n+i})$ landing at a  commun parameter $\mathbf{a}_{n+i}$ and
$\RR_{0}(\tau'_{n+i})$, $\RR_{\infty}(\zeta'_{n+i})$ landing at
some other  parameter $\mathbf{a}'_{n+i}$, at least for $n$ in the
subsequence $(n_j)_{j\ge 0}$ defined in Proposition~\ref{p:dyn}.

3. {\it We prove now that $\UU\subset \WW(\mathbf{a}_1)$}\,:

\noindent The pieces $\PP_n(\mathbf{a}_0)$ and
$\PP_n(\mathbf{a}_1)$ coincide. As in Proposition~\ref{p:H0interM}
applied to $\mathbf{a}_1$, the sequences of angles
$(\tau_n),(\tau'_n)$ admit a common limit $\tau$ which is, by
point 2 above, equal to  the common limit $t_1$ of the sequences
$(t_n)_{n\ge 2}, (t'_n)_{n\ge 2}$. The monotone sequences
$(\zeta_n)_{n\ge 2}$, $(\zeta'_n)_{n\ge 2}$ converge to limits
$\zeta_1$, $\zeta'_1$ respectively. By
Proposition~\ref{p:H0interM}, the angle  $\tau$ is $k$-periodic by
multiplication by $d-1$, the angles  $\zeta_1$, $\zeta'_1$ are
$k$-periodic by multiplication by $d$ and the rays $\RR_0(t_1)$,
$\RR_\infty(\zeta_1)$, $\RR_\infty(\zeta'_1)$ land at
$\mathbf{a}_1=\chi_{\mathbf{a}_0}^{-1}(1/4)$. The curve
$\ol{\RR_\infty(\zeta_1)} \cup \ol{\RR_\infty(\zeta'_1)}$ bounds
$\WW(\mathbf{a}_1)$.

Let $Q_n$  the connected component containing $\mathbf{a}_0$ of
$$\C \setminus\left( \ol{\RR_0(\tau_{n+i})}\cup \ol{\RR_0(\tau'_{n+i})}
\cup\ol{ \RR_\infty(\zeta_{n+i})}\cup\ol{
\RR_\infty(\zeta'_{n+i})}\right).$$ The para-puzzle piece
$\PP_{n}$ is contained in $Q_n$ and $\RR_\UU(\tau_n)$ is in the
boundary of $\PP_n$. Thus $\UU$ is included in $Q_n$ since
$\UU\cap \partial Q_n=\emptyset$.
  Thus the component $\UU$ is totally included
in $\WW(\mathbf{a}_1)$.

4. { \it Finally, we prove that $\RR_\UU(t_1)$,
$\RR_\infty(\eta)$, $\RR_\infty(\eta')$ land at the same
parameter, which  is $\mathbf{a}_0$, and that the curve
$\gamma=\ol{\RR_\infty(\eta)}\cup\ol{\RR_\infty(\eta')}$ separates
$\UU$ from $M_{\mathbf{a}_0}$}\,:

\noindent
 
Since $t_1, \eta,\eta'$ are rational ($d^i\eta=\zeta_1$ and
$d^i\eta'=\zeta'_1$) the rays $\RR_\UU(t_1)$, $\RR_\infty(\eta)$,
$\RR_\infty(\eta')$  converge to  parameters $\mathbf{a}_{t_1}$,
$\mathbf{a}_\eta$, $\mathbf{a}_{\eta'}$ respectively which are
either parabolic  or Misiurewicz parameters after
Lemma~\ref{l:misiurpreim1}.

If  $\mathbf{a}_{t_1}$ is a Misiurewicz  parameter, the ray
$R^{r(\mathbf{a}_{t_1})}_{\mathbf{a}_{t_1}}(t_1)$
 converges to $f_{\mathbf{a}_{t_1}}(-\mathbf{a}_{t_1})$
 (Lemma~\ref{l:misiurpreim1}). Moreover,
the rays  $R^0_{\mathbf{a}_{t_1}}(t_1)$,
$R^\infty_{\mathbf{a}_{t_1}}(\zeta_1)$ and
$R^\infty_{\mathbf{a}_{t_1}}(\zeta'_1)$ land at the same point
which
 is repelling (Lemma~\ref{l:wake}). Indeed,
 $\mathbf{a}_{t_1}$ is in $\WW(\mathbf{a}_1)$ since $\ol {\RR_\UU(t_1)}\subset
\WW(\mathbf{a}_1)\cup\{\mathbf{a}_1\}$ but $\mathbf{a}_{t_1}\neq
\mathbf{a}_1$ since $\mathbf{a}_1$ is a parabolic parameter.
 Pulling back along the critical orbit we obtain that
 $R^{r(\mathbf{a}_{t_1})}_{\mathbf{a}_{t_1}}(t_1)$,
 $R_{\mathbf{a}_{t_1}}^\infty(\eta)$ and $R_{\mathbf{a}_{t_1}}^\infty(\eta')$
land at the same point {\it i.e.}, at
$f_{\mathbf{a}_{t_1}}(-\mathbf{a}_{t_1})$. Therefore, by
Lemma~\ref{l:misiurewicz} the rays $\RR_\UU(t_1)$,
$\RR_\infty(\eta)$,  $\RR_\infty(\eta')$ land at the same
parameter $\mathbf{a}_{t_1}$. The proof is exactly the same in
case $\mathbf{a}_\eta$, or $\mathbf{a}_{\eta'}$, is of Misiurewicz
type.

Assume now that every parameter $\mathbf{a}\in\{\mathbf{a}_{t_1},
\mathbf{a}_\eta, \mathbf{a}_{\eta'}\}$ is parabolic. Then, the
landing point of $R_\mathbf{a}^0(t_1)$
 is a parabolic periodic point. Thus, the map $f^k_\mathbf{a}$ possesses a
parabolic fixed point of  multiplier $1$. Then,  the only
possibility for $\mathbf{a}\in \M_{\mathbf{a}_0}$ is to be the
cusp of $\M_{\mathbf{a}_0}$, {\it i.e.}
$\mathbf{a}=\chi^{-1}_{\mathbf{a}_0}(1/4)$ so
 $\mathbf{a}_{t_1}= \mathbf{a}_\eta= \mathbf{a}_{\eta'}$.

In both cases, the  curve
$\gamma=\RR_\infty(\eta)\cup\RR_\infty(\eta') \cup
\{\mathbf{a}_{t_1}\}$ separates the plane into two components. Let
$V$ denotes the one containing $\UU$ (since $\UU\cap
\gamma=\emptyset$). The para-puzzle piece $\PP_n$ intersects $V$
and $\UU$. Using Lemma~\ref{l:connexe} one can see that for  any
ray $\RR_\infty(\xi)$ in $V\cap \PP_n$, the angle $\xi$ is either
in $(\eta', \eta'_n)$ or in $(\eta_n,\eta)$. Assume (to get a
contradiction) that $\M_{\mathbf{a}_0}$ intersects $V$. Then let
$\mathbf{a}$ be a Misiurewicz parameter in the intersection
$\M_{\mathbf{a}_0}\cap V$. It is the landing point of an external
ray $\RR_\infty(\xi)$. This ray belongs to $V$ but also enters
every $\PP_n$ since it converges to $\mathbf{a}\in
\M_{\mathbf{a}_0}$. Hence $\xi$ is either in $(\eta', \eta'_n)$ or
in $(\eta_n,\eta)$, so $\xi=\eta$ or $\xi=\eta'$. Then
$\mathbf{a}=\mathbf{a}_{t_1}$ but this contradicts the fact that
$\mathbf{a}\in V$ (since $\mathbf{a}_{t_1}\in
\partial V$). Therefore, the curve $\gamma$ separates $\M_{\mathbf{a}_0}$
from $\UU$ and the unique possible intersection between $\ol \UU$
and $\M_{\mathbf{a}_0}$ is $\mathbf{a}_{t_1}$, so that
$\mathbf{a}_0=\mathbf{a}_{t_1}$ since $\mathbf{a}_0\in
\M_{\mathbf{a}_0}$ and $\mathbf{a}_0\in \ol \UU$. \cqfd

We will see in Theorem~\ref{th:U} that $\mathbf{a}_0$ is always a
Misiurewicz parameter.

\begin{corollary}\label{c:renbordU}
If $\UU$ is a connected component of $\HH\setminus\HH_0$, there
exists at most one parameter $\mathbf{a}$ on $\partial \UU$ such
that $f_\mathbf{a}$ is renormalizable. If it exists it is the
parameter characterized in Proposition~\ref{p:UcapM}.
 \end{corollary}
\proof Let $\mathbf{a}$ be a parameter on $\partial \UU$ such that
$f_\mathbf{a}$ is renormalizable.  By Proposition~\ref{p:UcapM}
$\M_{\mathbf{a}}=\cap\PP_n(\mathbf{a})$ intersects $\ol\HH_0$ and
$\ol\UU$. So if there is another point $\mathbf{a}'$ like this,
one can find a loop in
$\ol\HH_0\cup\ol\UU\cup\M_\mathbf{a}\cup\M_{\mathbf{a}'}$
surrounding points of $\HH_\infty$ and this contradicts the fact
that $\HH_\infty\cup\{\infty\}$ is simply connected. \cqfd

 \begin{corollary}\label{c:lcU} The boundary of
$\UU$ is locally connected.
 \end{corollary}
\proof Let $\mathbf{a}_0\in \partial \UU$. If $f_{\mathbf{a}_0}$
is not renormalizable the result follows from
Lemma~\ref{c:modules}. If $f_{\mathbf{a}_0}$ is  renormalizable,
we are in the situation of Proposition~\ref{p:UcapM}. The sequence
$\ol{\PP_n} \cap
\partial \UU$ of connected neighbourhoods of $\mathbf{a}_0$
(Lemma~\ref{l:connexe}) form a basis since $(\Cap_{n\ge
0}\ol{\PP_n} \cap \partial \UU)\subset
(\M_{\mathbf{a}_0}\cap\partial \UU)=\{\mathbf{a}_0\}$.\cqfd

\begin{thmi} {\bf\ref{th:H}.} For every connected component $\UU$ of $\mathcal H$,
the boundary $\partial \UU$ is a Jordan curve.
\end{thmi}
\proof Every component $\UU$ of $\HH$ is simply connected
(Lemma~\ref{l:1connected}) and its boundary is locally connected
(Corollary~\ref{c:lcU}). Therefore  any conformal map $\Psi \from
\D \to\UU$ extends continuously to a map $\ol\Psi \from \ol\D \to
\ol\UU$ by Caratheodory's Theorem. Thus  the boundary $\partial
\UU$ is the curve\,: $\ol\Psi(\S^1)$. We prove that it is a Jordan
curve by contradiction. If $\ol\Psi(\S^1)$ is not a simple curve,
there is a crossing point $z$ of $\ol\Psi(\S^1)$ and one can find
a simple closed cuve $\gamma$ in $\UU\cup\{z\}$ surrounding points
of $\HH_\infty$ (since $\ol\Psi(\S^1)\subset \partial
\HH_\infty$). This contradicts the maximum principle, exactly as
in Lemma~\ref{l:1connected}, applied to the map $\mathbf{a}
\mapsto f^N_\mathbf{a}(-\mathbf{a})$ for some large $N$. \cqfd

\begin{thmi}{\bf\ref{t:Hyp}.} The boundary of every bounded hyperbolic
component is a Jordan curve.
\end{thmi}
\proof The boundary of a hyperbolic component which is contained
in $\HH$ is a Jordan curve by Theorem~\ref{th:H} above. The other
bounded hyperbolic components are connected components of the
interior of a copy of $\M$ (by Corollary~\ref{c:H2}). Thus the
boundary of such a component is the image by a homeomorphism of
the boundary of a hyperbolic component of the interior of $\M$.
Therefore it is a Jordan curve.
  \cqfd

  The following result is the analogue of
  Theorem~\ref{th:comp} for the captures components.
  We can go out of a capture component $\UU$
  (and stay in $\CC$) only through
  the landing point of rays of angles $\xi$ such that
  $(d-1)^j\xi=t_1$ where $\UU\subset \WW^s(t_1)$.

\begin{thmi}{\bf\ref{th:U}.}
Let $\UU$ be  a connected component of $\HH_i\setminus \HH_0$
{\rm(}with $i \in \N${\rm)}. Any parameter $\mathbf{a}\in\partial
\mathcal U$ is the landing point of a unique ray $\RR_\UU(\xi)$.
In the dynamical plane,  $R^{r(\mathbf{a})}_\mathbf{a}(\xi)$
converges to $f_\mathbf{a}(-\mathbf{a})$, which is not on
$\partial B_\mathbf{a}$ but $f_\mathbf{a}^{i+1}(-\mathbf{a}) \in
\partial B_\mathbf{a}$. No parameter $\mathbf{a}\in
\partial \UU$ can be of parabolic type. If $\xi\in \Q$ the parameter
$\mathbf{a}$ is of Misiurewicz type.
Moreover, let $t_1$ be such that $\UU\subset \WW^s(t_1)$ {\rm(}see
Proposition~\ref{p:UcapM}{\rm)}. Then\,:
\begin{itemize}
 \item if  $(d-1)^j\xi=t_1$,    there are exactly two  external
 parameter rays converging to~$\mathbf{a}$\,;
 \item otherwise there is exactly one external parameter ray converging to $\mathbf{a}$.
 \end{itemize}\end{thmi}

\begin{figure}[!h]\vskip 0cm
\centerline{\input{antenne.pstex_t}} \caption{Zoom on a capture
$\UU\neq \HH_0$.}\label{f:antennes}
 \end{figure}

\proof
 
The proof is similar to that of Theorem~\ref{th:comp}. The
boundary $\partial \UU$ is locally connected by
Theorem~\ref{th:H}. So, the coordinate $\Phi_\UU^{-1} \from \D \to
\UU$ extends to a continuous map from the closure $\ol \D$ to $\ol
\UU$. Hence, any point $\mathbf{a}$ of $\partial \UU$ is the
landing point of a unique ray $\RR_\UU(\xi)$. Let $U_\mathbf{a}$
be the connected component of $\wt B_\mathbf{a}$ containing $\va$
for $\mathbf{a}\in \UU$, let $r(\mathbf{a})$ be its center.  The
dynamical ray $R^{r(\mathbf{a})}_\mathbf{a}(\xi)$ converges to a
point $z(\mathbf{a})\in
\partial U_\mathbf{a}$ since $\partial U_\mathbf{a}$ is  locally connected.

{\it If $f_\mathbf{a}$ is not renormalizable
 then $z(\mathbf{a})=\va$}\,:

\noindent  As in the proof of  Theorem~\ref{th:comp}, the point
$z(\mathbf{a})$ belongs to every $P_n^\mathbf{a}$ for  $n\ge 0$ so
that $\cap P_n^\mathbf{a}=f_\mathbf{a}(-\mathbf{a})$. Thus
$z(\mathbf{a})=\va$.

{\it If $f_\mathbf{a}$ is renormalizable then $z(\mathbf{a})$ is
{\rm(}eventually{\rm)} parabolic or $z(\mathbf{a})=\va$ \,:}

\noindent If $f_\mathbf{a}$ is renormalizable,
Proposition~\ref{p:UcapM} insures that $\xi$ is rational (with the
unicity proved just above). Then $z(\mathbf{a})$ is either a
(eventually) parabolic point
 or $z(\mathbf{a})=\va$ (see Lemma~\ref{l:misiurpreim1}).

{\it The critical value is not on $\partial B_\mathbf{a}$}\,:

\noindent Assume, to get a contradiction, that
$f_\mathbf{a}(-\mathbf{a})\in
\partial B_\mathbf{a}$. Then, it is the landing point of
 exactly one ray of $B_\mathbf{a}$, say  $R_\mathbf{a}^0(\xi')$.
 Note that the parameter $\mathbf{a}$ cannot be of parabolic  type since
$\va$ is not in a Fatou component. Then $z(\mathbf{a})=\va$ by the
two remarks above.
 Thus the two rays   $f_\mathbf{a}^i(R_\mathbf{a}^{r(\mathbf{a})}(\xi))=R_\mathbf{a}^0(\xi)$ and
$f_\mathbf{a}^i(R_\mathbf{a}^{0}(\xi'))=R_\mathbf{a}^0((d-1)^{i}\xi')$
converge to $f_\mathbf{a}^{i}(z(\mathbf{a}))$. By
Lemma~\ref{l:aboutdistray0}, this implies that
$(d-1)^{i}\xi'=\xi$. Then, the two different rays
$R_\mathbf{a}^0((d-1)^{i-1}\xi')$ and
$f_\mathbf{a}^{i-1}(R_\mathbf{a}^{r(\mathbf{a})}(\xi))$ landing at
$f_\mathbf{a}^{i-1}(z(\mathbf{a}))$ have the same image by
$f_\mathbf{a}$. Thus,
$f_\mathbf{a}^{i-1}(z(\mathbf{a}))=-\mathbf{a}$, so
$f_\mathbf{a}^{i}(-\mathbf{a})=-\mathbf{a}$ but this is not
possible  since the critical point  would be periodic on the Julia
set.

{\it The parameter $\mathbf{a}$ is not parabolic }\,:

\noindent Assume to get a contradiction that $\mathbf{a}$ is of
parabolic type. Then, by
Proposition~\ref{p:periodicrenormalizable}, $f_\mathbf{a}$ is
renormalizable, so that  $\xi$ is in $\Q$ and $z(\mathbf{a})$ is a
(eventually) parabolic periodic point. Proposition~\ref{p:UcapM}
insures that $\UU$  is in some wake $\WW(\mathbf{a}_1)$ where
$\mathbf{a}_1\in
\partial \HH_0$.
 Then $\mathbf{a}\in \ol{\WW(\mathbf{a}_1)}$, but in
$\WW(\mathbf{a}_1)$ the point $z(\mathbf{a})$ should be eventually
repelling, so $\mathbf{a}=\mathbf{a}_1$. Hence, the four parameter
rays (of Proposition~\ref{p:UcapM}) $\RR_\infty(\zeta_1)$,
$\RR_\infty(\zeta'_1)$, $\RR_\infty(\eta)$ and $\RR_\infty(\eta')$
land at the same point $\mathbf{a}_1$. But this is not possible by
Theorem~\ref{th:comp} since $\a_1\in\partial\HH_0$.
 
{\it If $\xi $ is rational $\mathbf{a}$ is of Misiurewicz type}\,:

\noindent  The conclusion of the points before is that
$z(\mathbf{a})=\va$.
Then $\va$
is eventually repelling (since $\xi$ is eventually periodic by
multiplication by $d$). Thus $\mathbf{a}$ is of Misiurewicz type.

 {\it The iterate $f_\mathbf{a}^{i+1}(-\mathbf{a})$ belongs to $\partial B_\mathbf{a}$}\,:

 Indeed, $z(\mathbf{a})=\va\in \partial U_\mathbf{a}$ so that $f_\mathbf{a}^i(\va) \in \partial B_\mathbf{a}$.

{\it Finally, we consider the number of external parameter rays
converging to $\mathbf{a}$}\,:

\noindent  We have seen above that $\mathbf{a}$ is of Misiurewicz
type and that the ray $R_\mathbf{a}^{r(\mathbf{a})}(\xi)$ is
landing at $\va$.

First case\,: $(d-1)^j\xi=t_1$. Since   $\mathbf{a}\in
\WW(\mathbf{a}_1)$, the rays $R_\mathbf{a}^0(t_1)$,
$R_\mathbf{a}^\infty(\zeta_1)$ and $R_\mathbf{a}^\infty(\zeta'_1)$
land at a common point $y$. We can pull back these rays and get
that two external rays $R_\mathbf{a}^\infty(\eta)$ and
$R_\mathbf{a}^\infty(\eta')$ landing at the same point as
$R_\mathbf{a}^{r(\mathbf{a})}(\xi)$, that is at $\va$. Since
$\mathbf{a}$ is of Misiurewicz type, the external parameter rays
$\RR_\infty(\eta),\RR_\infty(\eta')$ land at $\mathbf{a}$.
 Assume (to get a
contradiction) that there is a third external parameter ray
$\RR_\infty(\xi')$ converging to $\mathbf{a}$,
$\xi'\notin\{\eta,\eta'\}$. Then, since $\mathbf{a}$ is of
Misiurewicz type, the ray $R_\mathbf{a}^\infty(\xi')$ lands at
$\va$ (Lemma~\ref{l:misiurpreim1}). Thus, since $(d-1)^j\xi=t_1$,
$f_\mathbf{a}^i(R_\mathbf{a}^\infty(\xi'))$ gives a third external
ray converging to $y$, since $d^i\xi'\notin\{\zeta_1,\zeta_1'\}$.
But this contradicts Lemma~\ref{l:2rayext} which insures that at
most two external rays land at a point of $\partial B_\mathbf{a}$.

Now we do the other cases. Assume (to get a contradiction) that
there are (at least) two external dynamical rays converging to
$\va$. Since $\va$ is not on $\partial B_\mathbf{a}$ but
$y=f_\mathbf{a}^i(\va)\in
\partial B_\mathbf{a}$, we iterate the three rays $i$-times and get
$R_\mathbf{a}^0(\xi)$ and two external rays, say
$R_\mathbf{a}^\infty(\theta_1)$ and
$R_\mathbf{a}^\infty(\theta_2)$, converging to $y$. The two
external rays separate $\C$ into two connected components. Let $V$
be the one
 containing $B_\mathbf{a}$ and $V_0$  be the other one.  Note that, as long as the
 $V_j=f^j(V_0)$  does not contain a  critical point of $f_\mathbf{a}$, its image by
 $f_\mathbf{a}$ is still a sector not containing $B_\mathbf{a}$ and  attached by  $f^{i+1}(y)$
 to $\partial B_\mathbf{a}$ (since $-\mathbf{a}$ is not in the orbit of $y$).
 From
Lemma~\ref{l:2rayext}, some iterate of $V_0$ has to contain a
critical point. If $V_j$ contains $0$, then $V_{j-1}$ contains
$B'_\mathbf{a}$ and so the preimage of $f^j(y)$ on $\partial
B_\mathbf{a}'$. Then  $f_\mathbf{a} \from V_{j-1} \to V_j$ is not
a homeomorphism so $-\mathbf{a}\in V_{j-1}$. We consider the first
$j$ such that $-\mathbf{a}\in V_{j-1}$.
  The sector  $V_{j-1}$ contains a preimage of
$\ol{R_\mathbf{a}^\infty(d^{j}\theta_1)}$ and of
$\ol{R_\mathbf{a}^\infty(d^{j}\theta_2)}$ landing at a preimage of
$f^{j}_\mathbf{a}(y)$.
 Let  $V'_{j}$ be the connected component of
 $\C \setminus (\ol{R_\mathbf{a}^\infty(d^{j}\theta_1)}\cup
 \ol{R_\mathbf{a}^\infty(d^{j}\theta_2)}$) containing the critical value.
 Denote by $V'_{j-1}$ the connected component of $f_\mathbf{a}^{-1}(V'_{j})$
 containing $-\mathbf{a}$. Then $V'_{j-1}\subset V_{j-1}$.
 Applying Lemma~\ref{l:2rayext} to $V'_{j}$, we get that
 some iterate $f_\mathbf{a}^r(V'_{j})$ contains $-\mathbf{a}$,  so if $r$ is
 the smallest iterate to satisfy this condition  one has\,: $f_\mathbf{a}^r(V'_{j})=V_{j-1}$.
 Hence the  landing point $f_\mathbf{a}^{j-1}(y)=f_\mathbf{a}^{j+r}(y)$ is a repelling periodic
 point. We  thicken the sector $V_{j-1}$ by adding a small disk
 $D$  around $f_\mathbf{a}^{j-1}(y)$ which satisfies $f_\mathbf{a}^{r+1}(D)\supset \ol D$
and  by taking external rays closed to
$R_\mathbf{a}^\infty(d^{j-1}\theta_1)$ and
$R_\mathbf{a}^\infty(d^{j-1}\theta_2)$ landing at points in $D$.
Denote by $V_{j+r}''$ this new domain\,: the domain $V''_{j+r}$ is
bounded by these two rays union a  part of $\partial D$ and
containing $D \cup V_{j-1}$. Take the inverse image of $V''_{j+r}$
by $f_\mathbf{a}^{r+1}$ along the previous orbit ({\it i.e.}
 backward along the orbit $\{f_\mathbf{a}^{j-1}(y),\cdots,
f_\mathbf{a}^{j+r}(y)\}$). We obtain a domain $V_{j-1}''$ with
$\ol{V_{j-1}''}\subset V_{j+r-1}''$ and $f_\mathbf{a}^r \from
V_{j-1}'' \to V_{j+r-1}''$ is proper of degree two.  Hence, since
$y=f_\mathbf{a}^i(f_\mathbf{a}(-\mathbf{a}))$, the forward orbit
by $f_\mathbf{a}^{r+1}$ of $-\mathbf{a}$ will stay in $V_{j-1}''$
so that $f_\mathbf{a}$ is renormalizable. This contradicts the
fact that two external rays  converge to $\va$ in the non
renormalizable case.

Assume now that two external rays
$\RR_\infty(\xi'),\RR_\infty(\xi'')$ land at $\mathbf{a}$.
 These two rays enter any para-puzzle piece $\PP_n(\mathbf{a})$ so by the homeomorphism
of Corollary~\ref{c:homeobord}   the rays
$R^\mathbf{a}_\infty(\xi')$ and $R^\mathbf{a}_\infty(\xi'')$ enter
all the pieces $P_n^\mathbf{a}$. Since the intersection $\Cap
_{n\ge 0}P_n^\mathbf{a}$ reduces to $\va$, the rays both converge
to the same point $z=\va$. But we have just seen that this is not
possible. \cqfd

\begin{cri}{\bf\ref{c:indiff}.} There is no parameter $\mathbf{a}$ on
the boundary of a connected component of $\HH$
 such that $f_\mathbf{a}$ has
 an irrational indifferent periodic point.
\end{cri}
\proof Let  $\UU$ be a connected component of $\HH$ and let
$\mathbf{a}$ be a parameter of $\partial \UU$. Assume to get a
contradiction that
 $f_\mathbf{a}$ admits  an irrational indifferent periodic point denoted by $x$.
 If $\UU=\HH_0$, Proposition~\ref{p:periodicrenormalizable} asserts
 that $f_\mathbf{a}$ is renormalizable, then Proposition~\ref{p:H0interM}
 implies that $\mathbf{a}$ is of parabolic type. If $\UU\in  \HH_i\setminus\HH_0$
 for some $i>0$,
 Theorem~\ref{th:U} gives that  $f_\mathbf{a}^{i+1}(-\mathbf{a})\in\partial B_\mathbf{a}$.
 Therefore $x\in\partial B_\mathbf{a}$ since a  subsequence of $(f_\mathbf{a}^{n}(-\mathbf{a}))$
 accumulates $x$.
Denote by  $R_\mathbf{a}^0(\xi)$ the ray landing at $x$ ($\partial
B_\mathbf{a}$ is a Jordan curve). Since $x$ is periodic, $\xi \in
\Q$ (Lemma~\ref{l:aboutdistray0}). This contradicts the Snail
Lemma (see~\cite{M1}) which asserts that the landing point of a
periodic ray  either is a repelling periodic point or has
multiplier equal to $1$. \cqfd

\section{Description of $\CC$ and size of the limbs of
$\HH_0$.}\label{s:limbs}
 
\subsection{Connections in $\CC$}

\begin{definition}\ Let $\M_0$ be a copy  of $\M$, let $\UU$ be
 a hyperbolic component. We say that  $\M_0$ and $\UU$ are attached
 if $\M_0$ intersects $\ol\UU$.
\end{definition}
\begin{definition}A {\it tip} of $\M$ is the landing point
of an external ray of the form $\RR_\M(p/2^n)$ with $0<p<2^n$.
 Equivalently it is a  parameter $c$ such that  $P_c^j(0)=\beta_c$
for some $j>0$  where $\beta_c$ is the non-separating fixed point
of $K(P_c)$ (here $P_c(z)=z^2+c$).

A {\it tip} of  a copy $\M_0$ of $\M$ is the image $\chi^{-1}(c)$
of a tip $c$ of $\M$, where $\M_0=\chi^{-1}(\M)$.
\end{definition}

\begin{proi}{\bf\ref{th:intercopy}.}
Let $\M_0$ be  a copy  of $\M$ and  $\UU$ be any  connected
component of $\HH$. If the intersection $\M_0\cap\ol\UU\neq
\emptyset$, then it reduces to exactly one point. Moreover this
intersection point is\,:
\begin{itemize} \item either the cusp of $\M_0$, if $\UU=\HH_0$\,;
\item or  a tip of $\M_0$, if $\UU\neq \HH_0$. Moreover, in this
case $\M_0$ is also attached to $\HH_0$ (by its cusp).
\end{itemize}
\end{proi}

\begin{figure}[!h]\vskip 0cm
\centerline{\input{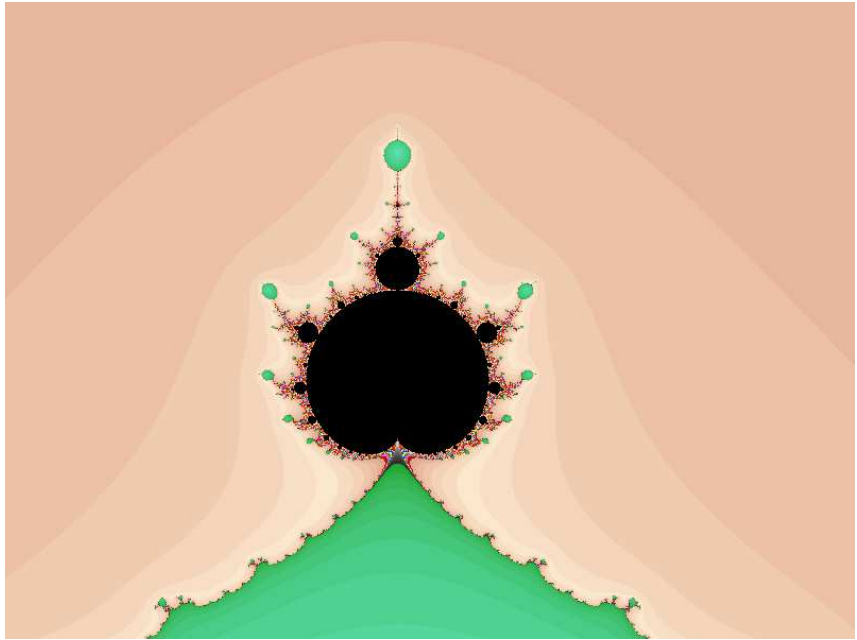}} \caption{Components $\UU$
attached to the tips of the copy $\M_0$ of $\M$.}\label{f:tips}
  \end{figure}

\proof Let  $\mathbf{a}$ be in $\M_0\cap \partial \UU$. Since
$f_\mathbf{a}$ is renormalizable, the intersection
$\M_\mathbf{a}=\cap\PP_n(\mathbf{a})$ is a copy of $\M$
(Proposition~\ref{p:UcapM}). We first prove that
$\M_0=\M_\mathbf{a}$.

If $\UU=\HH_0$,  the puzzle pieces $\PP_n(\mathbf{a})$ intersect
$\HH_0$ along a sector that contains $\mathbf{a}$ in its boundary.
Thus, they cannot cut $Card_0$, the main cardioid of $\M_0$. Then
$Card_0$ is contained in every $\PP_n(\mathbf{a})$ so in  the
intersection $\M_\mathbf{a}$. Thus $\M_\mathbf{a}=\M_0$. The proof
of Theorem~\ref{th:comp} insures that $\mathbf{a}$ is the cusp of
$\M_\mathbf{a}$. Hence, $\M_0\cap\partial \HH_0=\{\mathbf{a}\}$
the cusp of $\M_0$.

 If $\UU\neq \HH_0$, $\M_\mathbf{a}$ is attached to $\HH_0$
 by Proposition~\ref{p:UcapM}. Hence
$\M_\mathbf{a}=\M_0$, otherwise one can find a loop in
$\M_\mathbf{a}\cup\M_0\cup\ol\HH_0$ surrounding points of
$\HH_\infty$, which contradicts the fact that
$\HH_\infty\cup\{\infty\}$ is connected. Then, from
Proposition~\ref{p:UcapM} we get that
$\M_\mathbf{a}\cap\UU=\{\mathbf{a}\}$. Now we prove that
$\mathbf{a}$ is a tip of $\M_0$. From Theorem~\ref{th:U}, the
parameter $\mathbf{a}$ is of Misiurewicz type since
$x:=f^{i+1}_{\mathbf{a}}(-\mathbf{a})\in
\partial B_{\mathbf{a}}$ is an (eventually) repelling periodic point.
Some iterate $z=f_\mathbf{a}^r(x)$ belongs to
$K(f_\mathbf{a}^k)=\cap P_n^\mathbf{a}$. Then $z$ is a repelling
fixed point of $f^k_\mathbf{a}$ which does not separate
$K(f_\mathbf{a}^k)$. So $z$ corresponds  to the $\beta$ fixed
point of $z^2+\chi(\mathbf{a})$, the quadratic polynomial to which
$f_\mathbf{a}^k$ is conjugated.
 Therefore, $\mathbf{a}$ is a tip of $\M_{0}$.
 \cqfd
\begin{corollary}
Let  $\M_0$ be a copy of $\M$ attached to $\HH_0$ and contained in
a wake $\WW^s(t)$. Let $\UU$ be a connected component of
$\HH\setminus \HH_0$ attached to $\M_0$. The landing point of an
internal ray of $\UU$, $\RR_\UU(\xi)$, is a tip of $\M_0$ if and
only if $\xi=t$.
\end{corollary}
\proof Let $\mathbf{a}$ be the landing point of $\RR_\UU(\xi)$ and
let $\mathbf{a}_0$ be the intersection point $\ol \UU\cap \M_0$
(there is only one point by Proposition~\ref{th:intercopy}). we
prove first that $\a_0$ is the landing point of $\RR_\UU(t)$.

Since  $f_{\mathbf{a}_0}$ is renormalizable  we  can apply
Proposition~\ref{p:UcapM}.  Namely, the intersection
$\M_{\mathbf{a}_0}=\cap\PP_n(\mathbf{a}_0)$ is attached to $\HH_0$
and the intersection point is the landing point of $\RR_0(t)$,
since $\mathbf{a}_0\in \WW^s(t)$ and  since the curve
$\ol{\RR_\infty(\zeta)}\cup\ol{\RR_\infty(\zeta')}$ separates
$\ol\HH_0$ from the points of $\WW^s(t)$. Thus,
$\M_0=\M_{\mathbf{a}_0}$\,; otherwise one can find a loop in
$\M_0\cup\M_{\mathbf{a}_0}\cup\ol\HH_0$ surrounding points  of
$\HH_\infty$ which contradicts the fact that
$\HH_\infty\cup\{\infty\}$ is connected. Then
Proposition~\ref{p:UcapM} insures that $\mathbf{a}_0$ is the
landing point of $\RR_\UU(t)$.

Now, if $\mathbf{a}$ is a tip of $\M_0$,
$\mathbf{a}=\mathbf{a}_0$ since it is the only intersection point
between $\ol\UU$ and $\M_{\mathbf{a}_0}$
(Proposition~\ref{th:intercopy}). Therefore $\mathbf{a}$ is the
landing point of $\RR_\UU(t)$ and  $t=\xi$ by the unicity of
Theorem~\ref{th:U}.

Conversely, if $\xi=t$ since $\mathbf{a}$ is the landing point of
$\RR_\UU(\xi)$ and $\mathbf{a}_0$ is the landing point of
$\RR_\UU(t)$, $\mathbf{a}=\mathbf{a}_0$ so that $\mathbf{a}$ is a
tip of $\M_0$ (Proposition~\ref{th:intercopy}). \cqfd

\begin{proi}{\bf\ref{p:tips}.} If a copy $\M_0$ of $\M$ is attached to $\HH_0$, then at any
of its tips,  there is a connected component $\UU$ of
$\HH\setminus \HH_0$ which is attached.
\end{proi}

\proof Let  $\mathbf{a}$ be a tip of $\M_0$ and
$\mathbf{a}_0=\M_0\cap \ol\HH_0$\,; $\mathbf{a}_0$ is the cusp of
$\M_0$ by Proposition~\ref{th:intercopy}. Moreover,
$\M_0=\M_{\mathbf{a}_0}$ and
$\PP_n(\mathbf{a})=\PP_n(\mathbf{a}_0)$ for all $n\ge 0$, also
from the proof of Proposition~\ref{th:intercopy} above . Since
$\M_0$ is attached to $\HH_0$ every para-puzzle piece
$\PP_n(\mathbf{a})$ intersects $\HH_0$. Therefore, the puzzle
pieces $P_n^\mathbf{a}$ intersect the basin $B_\mathbf{a}$
 (applying the homeomorphism of Corollary~\ref{c:homeobord}).
 Then, the  intersection $K_{\mathbf{a}}=\cap P_n^\mathbf{a}$
 intersects $\ol B_\mathbf{a}$
 and since $K(f_\mathbf{a})$ is full,
 $K_{\mathbf{a}}\cap \ol B_\mathbf{a}$ reduces
 to one point, say $x$. Since $f_\mathbf{a}$ is renormalizable,
 there exists $k>0$ such that $f_\mathbf{a}^k$ maps $
 P_{n+k}^\mathbf{a}$ onto $P_n^\mathbf{a}$,
 so the point $x$ is $k$-periodic.
 Since $x\in\partial B_\mathbf{a}$, $x$ is the landing point
 of a unique  ray, say $R_\mathbf{a}^0(\tau)$. The unicity implies that
 $\tau $ is $k$-periodic by multiplication by $d-1$, so that
  $x$ is the non separating  fixed point of the renormalized map
  $f_\mathbf{a}^k$.
Since $\mathbf{a}$ is a tip of $\M_0$, $f_\mathbf{a}(-\mathbf{a})$
is mapped by some iterate of $f_\mathbf{a}^k$ to $x$. Thus some
iterated preimage $R^{r(\mathbf{a})}_\mathbf{a}(\tau)$ of
$R_\mathbf{a}^0(\tau)$ is landing at
$f_\mathbf{a}(-\mathbf{a})\neq x$ with $r(\mathbf{a})$ the center
of some connected component $U_\mathbf{a}$ of $\wt B_\mathbf{a}
\setminus B_\mathbf{a}$ (since $K(f_\mathbf{a})$ is full). The
puzzle pieces $P_n^\mathbf{a}$ intersect $U_\mathbf{a}$ for all
$n\ge 0$. Let denote by $R_\mathbf{a}^{r(\mathbf{a})}(\tau_n)$ and
$R_\mathbf{a}^{r(\mathbf{a})}(\tau'_n)$ the rays involved in
$U_\mathbf{a}\cap\partial P_n^\mathbf{a}$ and by
$R_\mathbf{a}^{\infty}(\eta_n)$, respectively
$R_\mathbf{a}^{\infty}(\eta'_n)$, the external rays of $\partial
P_n^\mathbf{a}$ converging to the landing points of the rays in
$U_\mathbf{a}$. By the homeomorphism of
Corollary~\ref{c:homeobord}, the para-puzzle pieces
$\PP_n(\mathbf{a})$ should intersect some component $\UU$ of
$\HH\setminus \HH_0$ and the rays involved in $\partial
\PP_n(\mathbf{a})\cap\UU$ are $R_\UU(\tau_n)$ and $R_\UU(\tau'_n)$
converging to the same points as the external rays
$R_{\infty}(\eta_n)$ and $R_{\infty}(\eta'_n)$ respectively (at
least for infinitely many $n\in \N$). The sequences
$(\tau_n),(\tau'_n)$ converge to $\tau$ and  the sequence
$(\eta_n)$, resp. $(\eta'_n)$ converges to $\eta$, resp. $\eta'$
by Propositions~\ref{p:H0interM} and~\ref{p:UcapM}. In the
dynamical plane $R_\mathbf{a}^{r(\mathbf{a})}(\tau)$,
$R_\mathbf{a}^{\infty}(\eta)$ and $R_\mathbf{a}^{\infty}(\eta')$
converge to the same point $\va$. Since $\mathbf{a}$ is a
Misiurewicz parameter, $\RR_\infty(\eta)$ lands at $\mathbf{a}$
(Lemma~\ref{l:misiurewicz}). Assume that the internal ray
$\RR_\UU(\tau)$ lands at a parameter $\mathbf{a}'\neq \mathbf{a}$.
Since $\RR_\UU(\tau)$ enters every puzzle pieces
$\PP_n(\mathbf{a})$, the parameter $\mathbf{a}'$ belongs to
$\M_\mathbf{a}=\M_0$. This contradicts the fact that
$\RR_\infty(\eta)$ will  have to land at $\mathbf{a}'$ by
Proposition~\ref{p:UcapM}.
  So $\UU$ is attached to $\M_0$ at $\mathbf{a}$.
\cqfd

\begin{lemi}{\bf\ref{l:inter}.}
Any two distinct components of $\HH$ have disjoint closures.
\end{lemi}
\proof Assume, to get a contradiction, that there exist $\UU_1$,
$\UU_2$ two distinct components of $\HH$  and $\mathbf{a}\in \ol
\UU_1\cap \ol \UU_2$. Since one connected component is distinct
from $\HH_0$, the parameter $\mathbf{a}$ is a Misiurewicz point by
Theorem~\ref{th:U}. Moreover Theorem~\ref{th:U} and
Theorem~\ref{th:comp} insure that there exist $\xi_1$, $\xi_2$
such that $\RR_\UU(\xi_1)$, $\RR_\UU(\xi_2)$ land at $\mathbf{a}$
and $R^{r_1(\mathbf{a})}_\mathbf{a}(\xi_1)$,
$R^{r_2(\mathbf{a})}_\mathbf{a}(\xi_2)$ land at $\va$ for some
adapted centers $r_1(\mathbf{a})$, $r_2(\mathbf{a})$. This is
impossible since $\va$ would be eventually critical by
Lemma~\ref{l:aboutdistray1}, although the critical point cannot be
periodic on the Julia set. \cqfd

\begin{thmi}{\bf\ref{th:conections}.} The only intersections
between closures of hyperbolic components, and also copies of $\M$
are the following\,:
\begin{itemize}
\item The central component $\HH_0$ has Mandelbrot copies $\M_t$
attached to it at angles $t$ which are $(d-1)$-periodic (a full
characterization of these values is given in
Proposition~\ref{p:convzero})\,:

\item At every tip  of such a satellite $\M_t$, a capture
component $\UU$ of $\HH\setminus \HH_0$ is attached.
\end{itemize}
Nevertheless, there are infinitely many  copies of $\M$ in $\CC$
and infinitely many captures  components not contained  in the
category described above.
\end{thmi}

\begin{figure}[!h]\vskip 0cm
\centerline{\input{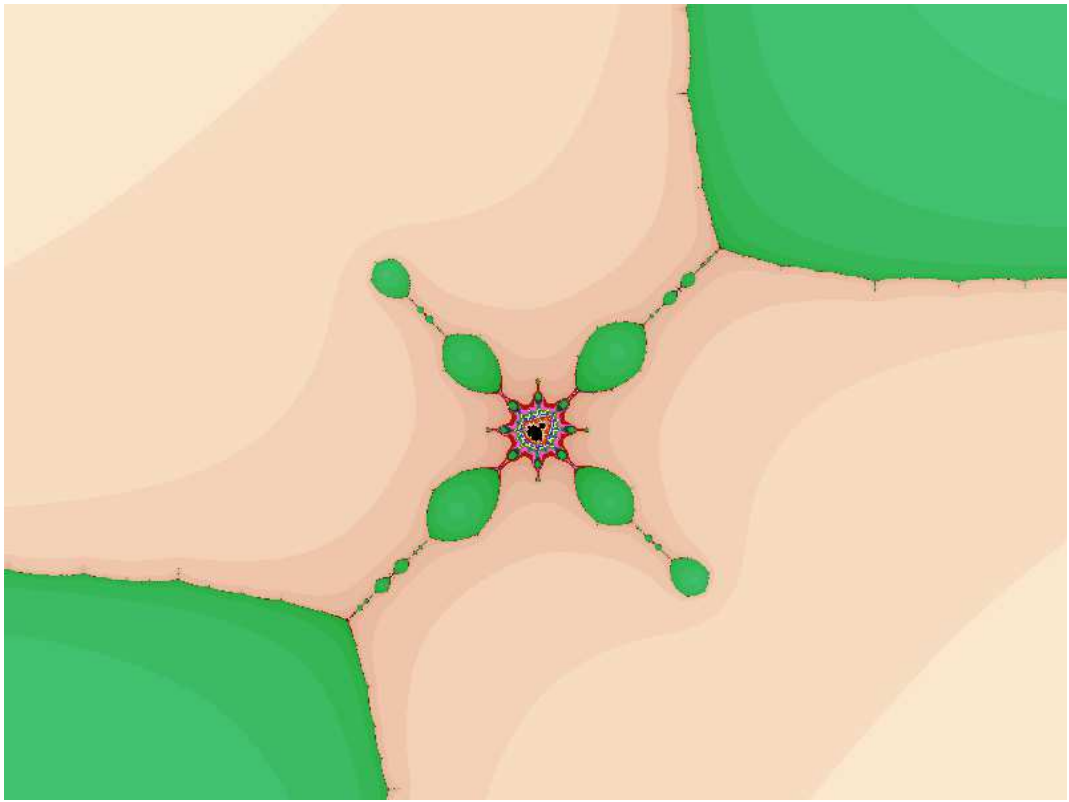}} \caption{Mandelbrot set
without connections with $\HH$.}\label{f:Misole}
  \end{figure}

\proof The description of the intersections of copies of $\M$ with
components of $\HH$ follows from Proposition~\ref{th:intercopy}
and~\ref{p:tips}. The intersections between components of $\HH$
follows from Lemma~\ref{l:inter}.

 There are other copies of $\M$ in $\CC$  (see figure~\ref{f:Misole}) which
  are not attached to components of $\HH$  since  Mandelbrot
  copies are dense in $\CC$ by~\cite{McM}.

  We prove now that there are capture components $\VV$ of
  $\HH\setminus \HH_0$ which are not attached to Mandelbrot copies
  (and also  not attached to hyperbolic components by
  Lemma~\ref{l:inter}).

  We start with a capture component $\UU$ attached to some
  satellite copy $\M_t$. Let $\mathbf{a}$ be the landing point of
  $\RR_\UU(\xi)$ with $(d-1)\xi=t$.
  In the dynamical plane, the critical value is on the boundary of
  some component $U_\mathbf{a}$ which is mapped by some
  $f_\mathbf{a}^i$ to $B_\mathbf{a}$.
  Since $\UU$ is in $\WW^s(t)$, the critical value $\va$ is
  separated from $B_\mathbf{a}$ by two external rays $R^\infty_\mathbf{a}(\theta_1),
  R_\mathbf{a}^\infty(\theta_2)$ landing at $f_\mathbf{a}^i(\va)$ a repelling
  periodic point (Proposition~\ref{p:UcapM} and
  Lemma~\ref{l:wake}).
  Then, considering the puzzle pieces, there is some $n$ such that
  $\partial P_n^\mathbf{a}$ involves rays in $U_\mathbf{a}$ and other
   rays in some iterated
  inverse image of $U_\mathbf{a}$ denoted by $V_\mathbf{a}$.
  Using the homeomorphism of Corollary~\ref{c:homeobord} we get
  that $\partial \PP_n(\mathbf{a})\cap\PP_{n-1}(\mathbf{a})$ contains rays in some
  capture component $\VV$ (because of the external parameter rays).
This para-puzzle piece is contained in $\WW^s(t)$ so that $\VV$ is
contained in $\WW^s(t)$. Therefore, $\VV$ cannot be attached to a
Mandelbrot copy. Otherwise this copy is  attached to $\HH_0$
(Proposition~\ref{p:UcapM}) so would coincide with $\M_0$.
Therefore $\UU$ and $\VV$ cannot be connected by an arc avoiding
$\M_0$, else one can complete it by an  arc  in $\UU\cup
\M_0\cup\VV$ to get a loop in $\CC$
 surrounding points of $\HH_\infty$ which contradicts
 the fact that $\HH_\infty\cup\{\infty\}$ is  connected.
But since $\CC\cap \PP_n(\mathbf{a})$ is connected, we get that
$\UU$ and $\VV$ are connected outside $\M_0$. Hence, $\VV$ is not
attached to a Mandelbrot copy. \cqfd
 
\subsection{Size of the limbs}
\begin{definition}
For any $\mathbf{a}\in \partial \HH_0$, we denote by {\it limb}
containing $\mathbf{a}$  the intersection
$\LL(\mathbf{a})=\ol{\WW(\mathbf{a})}\cap \CC$  if
$\WW(\mathbf{a})\neq \emptyset$ and else
$\LL(\mathbf{a})=\{\mathbf{a}\}$.
 \end{definition}

\begin{remark}By Lemma~\ref{l:wake}, the definition of the limbs
 coincides with the one given by Milnor in~\cite{M4}.
 \end{remark}

\begin{remark} Note that $\LL(\mathbf{a})\cap\partial \HH_0=\{\mathbf{a}\}$.
We will adopt sometimes the notation  $\LL^s(t)$ for a limb
$\LL(\mathbf{a})$ where $\mathbf{a}$ is the landing point of
$\RR_0^s(t)$.
\end{remark}

 The following result  was conjectured in~\cite{M4} (with another
parametrization of the rays). It follows from
 Theorem~\ref{th:comp}.
\begin{corollary}
A limb $\LL^s(t)$ contains more than one point if and only if the
angle $\frac{t}{d-1}+\frac{\lfloor\frac{d-1}{2}\rfloor}{d-1}$ is
periodic by multiplication by $d-1$.
\end{corollary}
\proof Let $\mathbf{a}_0$ be the point of $\partial \HH_0\cap
\LL^s(t)$, so $ \LL^s(t)=\LL(\mathbf{a}_0)$. From
Theorem~\ref{th:comp} there are two external rays converging to
$\mathbf{a}_0$ if and only if $\mathbf{a}_0$ is of parabolic type.
This corresponds to the statement about periodicity of
proposition~\ref{p:convzero}. \cqfd

\begin{thmi}{\bf\ref{th:limbs}.}\ For any $\epsilon>0$,
 there exists only a finite number
of limbs with diameter greater than $\epsilon$.
\end{thmi}
\proof This proof is inspired by Milnor's one in the quadratic
case. It proceeds in two steps.

\noindent {\it 1. Every point of $\CC \setminus \HH_0$ belongs to
a unique limb.}

 Let $\mathbf{a}\in \CC\cap\SS
\setminus \HH_0$. For every $n$, there exists a sector
$S_n(\mathbf{a})$ containing $\mathbf{a}$ such that
$\XX_n\cap\partial S_n(\mathbf{a})\subset \II_n$  is the union of
two rays in $\HH_0$ and two external rays. This sector contains a
para-puzzle piece, say $\PP_n$, which intersects $\HH_0$ along two
internal rays $\RR_0(t_n)$, $\RR_0(t'_n)$ and $\HH_\infty$ along
two external rays $\RR_\infty(\zeta_n)$, $\RR_\infty(\zeta'_n)$.

The study of the para-puzzle pieces (section~\ref{s:graphes})
gives the following informations. The sequences of angles $t_n$,
$t'_n$ converge to the commun limit $t$. Let $\mathbf{a}_0$ be the
landing point of $\RR_0(t)$. Either $\mathbf{a}_0$ is parabolic.
Then the sequences $(\zeta_n)$, $(\zeta'_n)$ converge to different
angles\,: $\zeta$ and  $\zeta'$ respectively, with
$\RR_\infty(\zeta)$, $\RR_\infty(\zeta')$ landing at
$\mathbf{a}_0$ (Theorem~\ref{th:comp}). Thus $\mathbf{a}$ belongs
to the wake $\WW(\mathbf{a}_0)$ so to the limb $\LL(t)$. Or,
$f_{\mathbf{a}_0}$ is not renormalizable, $\a=\a_0$ since
$\a\in\cap\PP_n=\{\a_0\}$, so $\cap P_n^{\mathbf{a}_0}=\va$ and
$\cap\PP_n=\{\mathbf{a}\}$. If two external parameter rays enter
the para-puzzle pieces $\PP_n$, by the homeomorphism of
Corollary~\ref{c:homeobord} the two corresponding dynamical rays
enter the puzzle pieces $P_n^{\mathbf{a}_0}$, so converge to
$f_{\mathbf{a}_0}(-\mathbf{a}_0)$. This contradicts
Theorem~\ref{th:comp}. Therefore the two external rays
$\RR_\infty(\zeta_n)$ and $\RR_\infty(\zeta'_n)$ converge to an
external ray $\RR_\infty(\zeta)$ which lands at $\mathbf{a}_0$.
Thus $\WW(\mathbf{a}_0)=\emptyset$ and
$\mathbf{a}=\mathbf{a}_0=\LL(t)$.
 
\noindent{\it 2. We assume (to get a contradiction) that there
exists a sequence $\LL(t_n)$ of limbs which accumulates  two
points $\mathbf{a}_1\neq \mathbf{a}_2$.}

 First suppose that
$\mathbf{a}_1,\mathbf{a}_2$ belong two different limbs $\LL(t_1)$
and $\LL(t_2)$ respectively. There exists $t\in (t_1,t_2)$ such
that $\RR_0(t)$ lands at a Misiurewicz parameter $\mathbf{a}$.
Then this point is the  landing point of an external ray
$\RR_\infty(\zeta)$. The curve $\RR_0(t)\cup\RR_\infty(\zeta)$
separates $\mathbf{a}_1$ from $\mathbf{a}_2$ so that the sequence
$\LL(t_n)$ cannot accumulate on both $\mathbf{a}_1$ and
$\mathbf{a}_2$.

Now suppose that
$\mathbf{a}_1$ and $\mathbf{a}_2$ belong to the
same limb $\LL(t)$,  we can assume that
$\mathbf{a}_1\notin\partial\HH_0$ (since $\partial\HH_0$ is
locally connected), so there
 exists $r>0$ such that  $B(\mathbf{a}_1,r) \subset \WW(\mathbf{a}_1)$. Therefore,
 if a sequence $\mathbf{a}_n\in \LL(t_n)$ accumulates $\mathbf{a}_1$, then for large
 $n$ the parameter $\mathbf{a}_n$ belongs to the ball $B(\a,r)$ so to $\WW(\a_1)$
 and therefore to $\LL(t)$ so $t_n=t$.
 \cqfd

\subsection{Local connectivity of $\partial \CC$ }
\begin{thmi}{\bf\ref{th:Clc}.}
$\partial \CC$ is locally connected at every point which
 is not in a  copy of $\M$ and at any point of $\partial \UU$ for
every connected components $\UU$ of $\HH$.
\end{thmi}
\proof Corollary~\ref{c:Clc} gives the local connectivity for
parameters which are not  in a  copy of $\M$.  Let
$\mathbf{a}_0\in
\partial\UU$ where $\UU$ is a connected  component of $\HH_i$ and $\PP_n$ be the para-puzzle piece containing
$\mathbf{a}_0$. If $\mathbf{a}_0$ is not parabolic
$\cap\PP_n=\{\mathbf{a}_0\}$ so $\partial \CC$ is locally
connected at this point by Lemma~\ref{l:connexeC}. If
$\mathbf{a}_0$ is parabolic, then the intersection $\cap\PP_n$ is
a copy $\M_0$ of $\M$. So to get the local connectivity at
$\mathbf{a}_0$ we consider the restricted puzzle pieces $\PP'_n$
as follows. Let $\gamma$ be the parametrization by the  internal
angle of the boundary of the main cardioid of $\M_0$ with
$\gamma(0)$ at the cusp. We consider an external rays converging
to $\gamma(1/n)$ and  another to $\gamma(-1/n)$ then we complete
their union to get with $\infty$ a closed curve $\delta_n$ by
adding some segment of curve  $c_n$ inside the main cardioid. Let
$\PP'_n$ be the connected component of $\PP_n\setminus \delta_n$
containing $\mathbf{a}_0$. The choice of $c_n$ is such that the
intersection of  $\PP'_n$ with the cardioid of $\M_0$ is a
sequence decreasing to the cusp of $\M_0$ (the boundary of
cardioid is locally connected). We prove now that the sequence
$\PP'_n$ is a basis of neighbourhoods of $\mathbf{a}_0$. Assume
that there is a point $\mathbf{a}\neq \mathbf{a}_0$ in
$\cap\PP_n'$. The parameter $\mathbf{a}$ is not on the cardioid of
$\M_0$, otherwise it would be some $\gamma(t)$ with
$\frac{-1}{n}\le t\le \frac{1}{n}$ so can only be $\a_0$. Then,
since the intersection of $\PP_n$ is $\M_0$,  the parameter $\a$
belongs to a Limb of $\M_0$ (the image of a limb of $\M$ by the
homeomorphism) say $\LL_\M(t)$. Thus for $n$ such that $1/n<|t|$
the point $\mathbf{a}$ does not belong to $\PP_n'$. This gives the
contradiction. The intersection $\CC\cap\PP_n'$ is clearly
connected by the same arguments as in Lemma~\ref{l:connexeC} since
we cut nicely the piece $\PP_n$ with the curve $\delta_n$. \cqfd


\section{Appendix}\label{s:annexe2}
For the completeness of the article, we recall  here the proof of
the following result given  in~\cite {F,Ro1}
\begin{thmi}[\cite{F,Ro1}]\ref{th:lcdyn} The boundary of every connected
component of $\wt B_\mathbf{a}$ is a Jordan curve.
\end{thmi}

\subsection{Yoccoz' Theorem for rational-like maps}\label{s:annexe2b}

\begin{definition}
A map $f\from X' \to X$ is {\it rational-like}  if\,:
\begin{itemize}
\item $X,X'$ are connected open sets of $\ol \C$ with smooth
boundary, such that $X\supset \ol{X'}$ and  $\partial X$ has a
finite number of connected components\,; \item $f\from X'\to X$ is
a holomorphic proper map with a finite number of
  critical points and extends to a continuous map from $\ol{X'}$ to $\ol X$.
\end{itemize}
\end{definition}

For a rational-like map $f\from X' \to X$, a graph $\Gamma$ is
{\it admissible} if\,:
\begin{itemize}
\item $\Gamma$ is connected, finite,  included in $\ol X$ and
contains
  $\partial X$\,;
\item $\Gamma$ is stable, {\it i.e.} $f^{-1}(\Gamma)\supset
\Gamma\cap X'$\.; \item the forward orbits of the critical points
are disjoint from $\Gamma$.
\end{itemize}

\noindent To an admissible graph, $\Gamma$, a puzzle is
associated\,:

\begin{definition}
The {\it puzzle pieces of depth $n$} are by definition the
connected components of $f^{-n}(X\setminus \Gamma)$.

The {\it end } of a point $x$ is  the nested sequence $(P_0(x),
P_1(x),\cdots, P_n(x),\cdots)$ of  the puzzle pieces containing
$x$.

The end of $x$ is {\it periodic} if there exists $k$ such that
$f^k(P_{k+n}(x))=P_n(x)$ for every~$n$ large enough.

The {\it impression} of $x$ is the intersection $\cap_{n\ge 0}
P_n(x)$ of the puzzle pieces containing $x$.

The point $x$ is {\it surrounded} at depth $i$  if  the annulus
$A_i=P_i\setminus\ol P_{i+1}$ surrounds $x$, {\it i.e.} if $x\in
P_{i+1}$ and $\ol {P_{i+1}}\subset P_i$.
\end{definition}

Yoccoz' Theorem can be stated in the context of rational-like maps
as follows\,:
\begin{thmi}Let $f\from X' \to X$ be a rational-like map with a unique critical
  point $x_0$ of multiplicity $2$ and $x$ be a point of $K(f)$.
If $\Gamma$ is an admissible graph that surrounds $x_0$ and surrounds infinitely many
times
$x$ then\,:
\begin{itemize}
\item if the  end of $x_0$ is not periodic, then the impression of
$x$ is equal to $\{x\}$\,; \item otherwise, let $k$ be the period
of the end of $x_0$, the map $f^k\from P_{l+k}(x_0)\to P_{l}(x_0)$
is quadratic-like for $l$ large enough and the impression of $x_0$
is the filled Julia set $K(f^k|_{P_{l+k}(x_0)})$ of the
renormalized map. Moreover the impression of $x$ reduces to $x$ or
to a preimage of the impression of $x_0$ if some iterate of $x$
falls in the impression of $x_0$.
\end{itemize}
\end{thmi}

\begin{remark}\label{r:7.3}
Let  $C$ be a forward invariant set under a rational-like map $f$.
A  compactness argument shows that instead of finding one
admissible graph and infinitely many annuli surrounding $x\in C$
with this graph, it is enough to find a finite number of
admissible graphs $\Gamma_0,\cdots,\Gamma_l$ such that every point
of $C$ is surrounded at bounded depth by one of these graphs which
surrounds also the critical point $x_0$.
\end{remark}

\subsection{Application to the family $f_\mathbf{a}$ }
It is enough to prove that $B_\mathbf{a}$ is locally connected\,:
one gets then  the result for  every connected component of $\wt
B_\mathbf{a}$ by pulling back .

\begin{remark}
If $\mathbf{a}\notin\CC$ the connected components of the  Julia
set are  locally connected since $f_\mathbf{a}$ is hyperbolic.
Thus we consider only parameters $\mathbf{a}\in \CC$.
\end{remark}

Let  $X$ be the connected  component of $\ol\C \setminus
(E_\mathbf{a}^\infty(1)\cup E_\mathbf{a}^0(1))$ containing
$J(f_\mathbf{a})$ and $X'=f_\mathbf{a}^{-1}(X)$. The map
 $f_\mathbf{a}\from X' \to X $ is  a rational-like map. We consider
the graphs given in section~\ref{s:graphes}. They are clearly
admissible. We  prove now that they  satisfy the conditions of
Yoccoz' Theorem\,: using  Remark~\ref{r:7.3} it suffices to show
that every point $x\in
\partial B_\mathbf{a}$ is surrounded, at bounded depth, by one of the
graphs which also surrounds  the critical point $-\mathbf{a}$.

\begin{lemma} For $\theta=\frac{1}{(d-1)^l-1}$ and $\theta'=\frac{1}{(d-1)^{l'}-1}$
with $l'>l+1$ and $l$ large enough, every point of $\partial
B_\mathbf{a}$ is surrounded  by one of the graphs at bounded
depth.
\end{lemma}
\proof Let $U(\theta)$ be the connected component of
$\ol\C\setminus \gamma$ containing $R^0_\mathbf{a}(0)$ where
$\gamma$ denotes the curve  in
$I_1^\mathbf{a}(\theta)=f_\mathbf{a}^{-1}(I_0^\mathbf{a}(\theta))$
formed by the internal rays $R_\mathbf{a}^0(\theta+\frac{1}{d-1})$
and $R_\mathbf{a}^0(\frac{\theta}{d-1})$ and the corresponding
external rays. One sees that every point of $\partial
B_\mathbf{a}$ which is in $U(\theta)$ but not on the graph
$I_1^\mathbf{a}(\theta)$ is surrounded at depth $0$ by the graph
$I_0(\theta)$.  Now the points on $\partial B_\mathbf{a}$ of
$I_1(\theta)$ and of $I_1(\theta')$ are distinct. Moreover the
union $V=U(\theta)\cup U(\theta')$ is the connected component of
the complement of $\delta$ containing $R^0_\mathbf{a}(0)$, where
$\delta$ is the curve formed by the internal rays
$R_\mathbf{a}^0(\theta'+\frac{1}{d-1})$ and
$R_\mathbf{a}^0(\frac{\theta}{d-1})$ and the corresponding
external rays. Then for $n$ such that
$\frac{1}{n}(\theta'+\frac{1}{d-1})<\frac{\theta}{d-1}$, the union
$\Cup_{i\le n}f_\mathbf{a}^{-i}(V)$ covers $\ol \C \setminus
\{0\}$. Thus every point of $\partial B_\mathbf{a}$ is surrounded
at depth less than $n$ by $I_0(\theta)$ or by $I_0(\theta')$.
\cqfd
\begin{remark}
This result  clearly holds  also  for
$\theta=\frac{-1}{(d-1)^l-1}$ and
$\theta'=\frac{-1}{(d-1)^{l'}-1}$.
\end{remark}
\begin{lemma}
There exists $k_0$ and $\epsilon\in\{\pm1\}$ (depending on
$\mathbf{a}$) such that for every $k>k_0$ the critical point
$-\mathbf{a}$ is surrounded at depth $1$ by the graph
$I_0(\epsilon\theta)$ where $\theta=\frac{1}{(d-1)^k-1}$.
\end{lemma}
\proof We take the open set $U(\theta)$ defined in previous Lemma.
For $d>3$ it is clear that the union $U(\theta)\cup U(-\theta)$
covers $B_\mathbf{a}\setminus\{0\}$ and therefore all $\C
\setminus\{0\}$. This solves the question for $d>3$.

For $d=3$,  the union $U(\theta)\cup U(-\theta)\cup
f_\mathbf{a}^{-1}(U(\theta)\cup U(-\theta))$ covers all
$B_\mathbf{a}\setminus\{0\}$ and therefore all $\C
\setminus\{0\}$.
 \cqfd

 The proof that the intersection  $\ol{P_n(x)}\cap\partial B_\mathbf{a}$ is
 a connected set is exactly the same as the proof of Lemma~\ref{l:connexe}
(in the parameter plane).

Yoccoz' Theorem (stated in section~\ref{s:annexe2b}) and previous
Lemmas allow to conclude that if the end of the critical point
$-\mathbf{a}$ is not periodic, the boundary $\partial
B_\mathbf{a}$ is locally connected. Then Caratheodory's Theorem
together with Lemma~\ref{l:aboutdistray0} insures that $\partial
B_\mathbf{a}$ is a Jordan curve.

We now consider the case  where the end of the critical point
$-\mathbf{a}$ is periodic of period $k$.

The map $f_\mathbf{a}^k\from P_{n+k}(-\mathbf{a})\to
P_n(-\mathbf{a})$ is quadratic-like and the orbit of the critical
point never escapes the puzzle piece $P_n(-\mathbf{a})$. So by the
straightening theorem of~\cite{DH2}, the restriction of
$f^k_\mathbf{a}$ is conjugated to a unique quadratic polynomial
$z^2+c$. Let $K_\mathbf{a}=K(f_\mathbf{a}^k)$ denotes the filled
Julia set of the restriction and $K$ the filled Julia set of the
quadratic polynomial.

We assume that $K_\mathbf{a}\cap \partial B_\mathbf{a}\neq
\emptyset$.
\begin{lemma}
There exists an internal ray $R_\mathbf{a}^0(\eta)$ of period $k$
converging to the non separating fixed point $\beta_\mathbf{a}$ of
$f_\mathbf{a}^k$ in $K_\mathbf{a}$ and two external rays
$R_\mathbf{a}^\infty(\zeta),R_\mathbf{a}^\infty(\zeta')$
converging to $\beta_\mathbf{a}$ and separating $K_\mathbf{a}$
from $B_\mathbf{a}$.
\end{lemma}
\proof The angles of the  internal rays  that bound the puzzle
piece $P_n(-\mathbf{a})$ are of the form $\eta_n<\eta'_n$ with
$(d-1)^k\eta_{n+k}=\eta_n \mod 1$ (and the same for $\eta'_n$)
with $|\eta'_n-\eta_n|<\frac{1}{(d-1)^{n-1}}$. Therefore they
converge to a common limit  $\eta$ which is periodic of period
$k$. Moreover the ray $R_\mathbf{a}^0(\eta)$ lies in all the
puzzles pieces $P_n(-\mathbf{a})$ so its landing point is in
$\partial B_\mathbf{a}\cap K_\mathbf{a}$. Since it is a fixed
point of $f_\mathbf{a}^k$ with rotation number $0$ it is the
non-separating fixed point of $K_\mathbf{a}$, {\it i.e.}
$\beta_\mathbf{a}$.

For the external rays the proof is the same. The external rays
attached to $R_\mathbf{a}^0(\eta_n)$ and $R_\mathbf{a}^0(\eta'_n)$
which are in the boundary of $P_n(-\mathbf{a})$  are of the form
$R_\mathbf{a}^0(\zeta_n)$ and $R_\mathbf{a}^0(\zeta'_n)$ and the
angle satisfies the equation $d^k\zeta_{n+k}^{(')}=\zeta_n^{(')}$.
Thus they converges to periodic angles $\zeta, \zeta'$. The rays
$R_\mathbf{a}^\infty(\zeta)$ and $R_\mathbf{a}(\zeta')$ converge
to $\beta_\mathbf{a}$ by the same argument as before. To see that
 the curve $R_\mathbf{a}^\infty(\zeta)\cup R_\mathbf{a}^\infty(\zeta')\cup\beta_\mathbf{a}$
 separates $B_\mathbf{a}$ from $K_\mathbf{a}\setminus\{\beta_\mathbf{a}\}$ it is enough to
 note that the preimage $\beta'_\mathbf{a}$ of $\beta_\mathbf{a}$  in $K_\mathbf{a}$  is the
 landing point of a ray of the form
 $R_\mathbf{a}^\infty(\frac{\zeta}{d}+\frac{i}{d})$ which is always
 contained in $P_n(-\mathbf{a})$ so  converges to $K_\mathbf{a}$ and separates
$R_\mathbf{a}^\infty(\zeta)$ from $ R_\mathbf{a}^\infty(\zeta')$.
\cqfd

\begin{corollary} The boundary of $\partial B_\mathbf{a}$ is locally
connected.
\end{corollary}
\proof When the end of the critical point $-\mathbf{a}$ is
periodic, the impression of the end is $K_\mathbf{a}$. If
$K_\mathbf{a}\cap\partial B_\mathbf{a}$  is empty the end of any
point forms a basis of connected neighbourhoods of that point. If
$K_\mathbf{a}\cap\partial B_\mathbf{a}$  is not empty, we take for
sequence of neighbourhoods of $\beta_\mathbf{a}$ in $\partial
B_\mathbf{a}$ the intersection $U_n=\ol{V_n}\cap\partial
B_\mathbf{a}$ where $V_n$ is the connected component of
$P_n(\beta_\mathbf{a})\setminus(R_\mathbf{a}^\infty(\zeta)\cup
R_\mathbf{a}^\infty(\zeta')\cup\beta_\mathbf{a})$ which intersects
$B_\mathbf{a}$. It is easy to see that the sequence $(U_n)$ forms
a basis of connected neighbourhoods of $\beta_\mathbf{a}$ in
$\partial B_\mathbf{a}$ since the intersection  $\cap U_n$ reduces
to $\partial B_\mathbf{a}\cap K_\mathbf{a}=\beta_\mathbf{a}$. Then
we pull back those neighbourhoods along the backward orbit of
$\beta_\mathbf{a}$. \cqfd

  \end{document}